\catcode`\@=11

\font\tenmsa=msam10
\font\sevenmsa=msam7
\font\fivemsa=msam5
\font\tenmsb=msbm10
\font\sevenmsb=msbm7
\font\fivemsb=msbm5
\newfam\msafam
\newfam\msbfam
\textfont\msafam=\tenmsa  \scriptfont\msafam=\sevenmsa
  \scriptscriptfont\msafam=\fivemsa
\textfont\msbfam=\tenmsb  \scriptfont\msbfam=\sevenmsb
  \scriptscriptfont\msbfam=\fivemsb

\def\hexnumber@#1{\ifnum#1<10 \number#1\else
 \ifnum#1=10 A\else\ifnum#1=11 B\else\ifnum#1=12 C\else
 \ifnum#1=13 D\else\ifnum#1=14 E\else\ifnum#1=15 F\fi\fi\fi\fi\fi\fi\fi}

\def\msa@{\hexnumber@\msafam}
\def\msb@{\hexnumber@\msbfam}
\mathchardef\boxdot="2\msa@00
\mathchardef\boxplus="2\msa@01
\mathchardef\boxtimes="2\msa@02
\mathchardef\square="0\msa@03
\mathchardef\blacksquare="0\msa@04
\mathchardef\centerdot="2\msa@05
\mathchardef\lozenge="0\msa@06
\mathchardef\blacklozenge="0\msa@07
\mathchardef\circlearrowright="3\msa@08
\mathchardef\circlearrowleft="3\msa@09
\mathchardef\rightleftharpoons="3\msa@0A
\mathchardef\leftrightharpoons="3\msa@0B
\mathchardef\boxminus="2\msa@0C
\mathchardef\Vdash="3\msa@0D
\mathchardef\Vvdash="3\msa@0E
\mathchardef\vDash="3\msa@0F
\mathchardef\twoheadrightarrow="3\msa@10
\mathchardef\twoheadleftarrow="3\msa@11
\mathchardef\leftleftarrows="3\msa@12
\mathchardef\rightrightarrows="3\msa@13
\mathchardef\upuparrows="3\msa@14
\mathchardef\downdownarrows="3\msa@15
\mathchardef\upharpoonright="3\msa@16

\mathchardef\downharpoonright="3\msa@17
\mathchardef\upharpoonleft="3\msa@18
\mathchardef\downharpoonleft="3\msa@19
\mathchardef\rightarrowtail="3\msa@1A
\mathchardef\leftarrowtail="3\msa@1B
\mathchardef\leftrightarrows="3\msa@1C
\mathchardef\rightleftarrows="3\msa@1D
\mathchardef\Lsh="3\msa@1E
\mathchardef\Rsh="3\msa@1F
\mathchardef\rightsquigarrow="3\msa@20
\mathchardef\leftrightsquigarrow="3\msa@21
\mathchardef\looparrowleft="3\msa@22
\mathchardef\looparrowright="3\msa@23
\mathchardef\circeq="3\msa@24
\mathchardef\succsim="3\msa@25
\mathchardef\gtrsim="3\msa@26
\mathchardef\gtrapprox="3\msa@27
\mathchardef\multimap="3\msa@28
\mathchardef\therefore="3\msa@29
\mathchardef\because="3\msa@2A
\mathchardef\doteqdot="3\msa@2B

\mathchardef\triangleq="3\msa@2C
\mathchardef\precsim="3\msa@2D
\mathchardef\lesssim="3\msa@2E
\mathchardef\lessapprox="3\msa@2F
\mathchardef\eqslantless="3\msa@30
\mathchardef\eqslantgtr="3\msa@31
\mathchardef\curlyeqprec="3\msa@32
\mathchardef\curlyeqsucc="3\msa@33
\mathchardef\preccurlyeq="3\msa@34
\mathchardef\leqq="3\msa@35
\mathchardef\leqslant="3\msa@36
\mathchardef\lessgtr="3\msa@37
\mathchardef\backprime="0\msa@38
\mathchardef\risingdotseq="3\msa@3A
\mathchardef\fallingdotseq="3\msa@3B
\mathchardef\succcurlyeq="3\msa@3C
\mathchardef\geqq="3\msa@3D
\mathchardef\geqslant="3\msa@3E
\mathchardef\gtrless="3\msa@3F
\mathchardef\sqsubset="3\msa@40
\mathchardef\sqsupset="3\msa@41
\mathchardef\trianglerighteq="3\msa@44
\mathchardef\trianglelefteq="3\msa@45
\mathchardef\bigstar="0\msa@46
\mathchardef\between="3\msa@47
\mathchardef\blacktriangledown="0\msa@48
\mathchardef\blacktriangleright="3\msa@49
\mathchardef\blacktriangleleft="3\msa@4A
\mathchardef\blacktriangle="0\msa@4E
\mathchardef\triangledown="0\msa@4F
\mathchardef\eqcirc="3\msa@50
\mathchardef\lesseqgtr="3\msa@51
\mathchardef\gtreqless="3\msa@52
\mathchardef\lesseqqgtr="3\msa@53
\mathchardef\gtreqqless="3\msa@54
\mathchardef\Rrightarrow="3\msa@56
\mathchardef\Lleftarrow="3\msa@57
\mathchardef\veebar="2\msa@59
\mathchardef\barwedge="2\msa@5A
\mathchardef\doublebarwedge="2\msa@5B
\mathchardef\angle="0\msa@5C
\mathchardef\measuredangle="0\msa@5D
\mathchardef\sphericalangle="0\msa@5E
\mathchardef\varpropto="3\msa@5F
\mathchardef\smallsmile="3\msa@60
\mathchardef\smallfrown="3\msa@61
\mathchardef\Subset="3\msa@62
\mathchardef\Supset="3\msa@63
\mathchardef\Cup="2\msa@64

\mathchardef\Cap="2\msa@65

\mathchardef\curlywedge="2\msa@66
\mathchardef\curlyvee="2\msa@67
\mathchardef\leftthreetimes="2\msa@68
\mathchardef\rightthreetimes="2\msa@69
\mathchardef\subseteqq="3\msa@6A
\mathchardef\supseteqq="3\msa@6B
\mathchardef\bumpeq="3\msa@6C
\mathchardef\Bumpeq="3\msa@6D
\mathchardef\lll="3\msa@6E

\mathchardef\ggg="3\msa@6F

\mathchardef\circledS="0\msa@73
\mathchardef\pitchfork="3\msa@74
\mathchardef\dotplus="2\msa@75
\mathchardef\backsim="3\msa@76
\mathchardef\backsimeq="3\msa@77
\mathchardef\complement="0\msa@7B
\mathchardef\intercal="2\msa@7C
\mathchardef\circledcirc="2\msa@7D
\mathchardef\circledast="2\msa@7E
\mathchardef\circleddash="2\msa@7F
\def\ulcorner{\delimiter"4\msa@70\msa@70 }
\def\urcorner{\delimiter"5\msa@71\msa@71 }
\def\llcorner{\delimiter"4\msa@78\msa@78 }
\def\lrcorner{\delimiter"5\msa@79\msa@79 }
\def\yen{\mathhexbox\msa@55 }
\def\checkmark{\mathhexbox\msa@58 }
\def\circledR{\mathhexbox\msa@72 }
\def\maltese{\mathhexbox\msa@7A }
\mathchardef\lvertneqq="3\msb@00
\mathchardef\gvertneqq="3\msb@01
\mathchardef\nleq="3\msb@02
\mathchardef\ngeq="3\msb@03
\mathchardef\nless="3\msb@04
\mathchardef\ngtr="3\msb@05
\mathchardef\nprec="3\msb@06
\mathchardef\nsucc="3\msb@07
\mathchardef\lneqq="3\msb@08
\mathchardef\gneqq="3\msb@09
\mathchardef\nleqslant="3\msb@0A
\mathchardef\ngeqslant="3\msb@0B
\mathchardef\lneq="3\msb@0C
\mathchardef\gneq="3\msb@0D
\mathchardef\npreceq="3\msb@0E
\mathchardef\nsucceq="3\msb@0F
\mathchardef\precnsim="3\msb@10
\mathchardef\succnsim="3\msb@11
\mathchardef\lnsim="3\msb@12
\mathchardef\gnsim="3\msb@13
\mathchardef\nleqq="3\msb@14
\mathchardef\ngeqq="3\msb@15
\mathchardef\precneqq="3\msb@16
\mathchardef\succneqq="3\msb@17
\mathchardef\precnapprox="3\msb@18
\mathchardef\succnapprox="3\msb@19
\mathchardef\lnapprox="3\msb@1A
\mathchardef\gnapprox="3\msb@1B
\mathchardef\nsim="3\msb@1C
\mathchardef\napprox="3\msb@1D
\mathchardef\nsubseteqq="3\msb@22
\mathchardef\nsupseteqq="3\msb@23
\mathchardef\subsetneqq="3\msb@24
\mathchardef\supsetneqq="3\msb@25
\mathchardef\subsetneq="3\msb@28
\mathchardef\supsetneq="3\msb@29
\mathchardef\nsubseteq="3\msb@2A
\mathchardef\nsupseteq="3\msb@2B
\mathchardef\nparallel="3\msb@2C
\mathchardef\nmid="3\msb@2D
\mathchardef\nshortmid="3\msb@2E
\mathchardef\nshortparallel="3\msb@2F
\mathchardef\nvdash="3\msb@30
\mathchardef\nVdash="3\msb@31
\mathchardef\nvDash="3\msb@32
\mathchardef\nVDash="3\msb@33
\mathchardef\ntrianglerighteq="3\msb@34
\mathchardef\ntrianglelefteq="3\msb@35
\mathchardef\ntriangleleft="3\msb@36
\mathchardef\ntriangleright="3\msb@37
\mathchardef\nleftarrow="3\msb@38
\mathchardef\nrightarrow="3\msb@39
\mathchardef\nLeftarrow="3\msb@3A
\mathchardef\nRightarrow="3\msb@3B
\mathchardef\nLeftrightarrow="3\msb@3C
\mathchardef\nleftrightarrow="3\msb@3D
\mathchardef\divideontimes="2\msb@3E
\mathchardef\varnothing="0\msb@3F
\mathchardef\nexists="0\msb@40
\mathchardef\mho="0\msb@66
\mathchardef\thorn="0\msb@67
\mathchardef\beth="0\msb@69
\mathchardef\gimel="0\msb@6A
\mathchardef\daleth="0\msb@6B
\mathchardef\lessdot="3\msb@6C
\mathchardef\gtrdot="3\msb@6D
\mathchardef\ltimes="2\msb@6E
\mathchardef\rtimes="2\msb@6F
\mathchardef\shortmid="3\msb@70
\mathchardef\shortparallel="3\msb@71
\mathchardef\smallsetminus="2\msb@72
\mathchardef\thicksim="3\msb@73
\mathchardef\thickapprox="3\msb@74
\mathchardef\approxeq="3\msb@75
\mathchardef\succapprox="3\msb@76
\mathchardef\precapprox="3\msb@77
\mathchardef\curvearrowleft="3\msb@78
\mathchardef\curvearrowright="3\msb@79
\mathchardef\digamma="0\msb@7A
\mathchardef\varkappa="0\msb@7B
\mathchardef\hslash="0\msb@7D
\mathchardef\hbar="0\msb@7E
\mathchardef\backepsilon="3\msb@7F
\def\Bbb{\ifmmode\let\next\Bbb@\else
 \def\next{\errmessage{Use \string\Bbb\space only in math mode}}\fi\next}
\def\Bbb@#1{{\Bbb@@{#1}}}
\def\Bbb@@#1{\fam\msbfam#1}

\catcode`\@=\active

\def\del{\partial}
 \def\CC{\hbox{{$\cal C$}}}
 
 \def\CH{\hbox{{$\cal H$}}}
 \def\CU{\hbox{{$\cal U$}}}
 
\def\CR{\hbox{{$\cal R$}}} 
 
\def\CM{\hbox{{$\cal M$}}}

\def\cu{\hbox{\sl u}} % used for special element u
\def\cv{\hbox{\sl v}} % used for special element v
\def\cg{\hbox{{\sl g}}} % used for Lie algebra 'gothic g'
\def\Vec{{\rm Vec}}

\def\lform{\hbox{$\sqcup$}\llap{\hbox{$\sqcap$}}}
\def\h{{{1\over2}}}

\def\R{{\Bbb R}}
\def\C{{\Bbb C}}
\def\Z{{\Bbb Z}}

\def\Q{{\Bbb Q}}

\def\A{{\Bbb A}}

\def\eps{{\epsilon}}

\def\rcross{{\triangleright\!\!\!<}}

\def\dcross{{\bowtie}}

\def\rbiprod{{\cdot\kern-.33em\triangleright\!\!\!<}}
\def\lbiprod{{>\!\!\!\triangleleft\kern-.33em\cdot\, }}

\def\tens{\mathop{\otimes}}

\def\la{{\triangleright}}\def\ra{{\triangleleft}}

\def\Ad{{\rm Ad}}

\def\ev{{\rm ev}}
\def\coev{{\rm coev}}

\def\id{{\rm id}}

\def\<{\langle}
\def\>{\rangle}

\def\equad{\kern -1.7em}
\def\qqquad{\qquad\quad}

\def\eqn#1#2{\begin{equation}#2\label{#1}\end{equation}}

\def\haj#1{{\mathaccent20 {#1}}}
\def\Vhaj{{V\haj{\ }}}

\def\o{{}_{\scriptscriptstyle(1)}}
\def\t{{}_{\scriptscriptstyle(2)}}
\def\th{{}_{\scriptscriptstyle(3)}}
\def\fo{{}_{\scriptscriptstyle(4)}}
\def\fiv{{}_{\scriptscriptstyle(5)}}
\def\six{{}_{\scriptscriptstyle(6)}}
\def\sev{{}_{\scriptscriptstyle(7)}}

\def\bo{{}^{\bar{\scriptscriptstyle(1)}}}
\def\bt{{}^{\bar{\scriptscriptstyle(2)}}}
\def\Ro{{\CR^{\scriptscriptstyle(1)}}}
\def\Rt{{\CR^{\scriptscriptstyle(2)}}}

\def\und#1{{\underline {#1}}}

\def\uo{{{}^{\scriptscriptstyle(1)}}}
\def\ut{{{}^{\scriptscriptstyle(2)}}}

\def\baro{{}_{\bar{\scriptscriptstyle(1)}}}
\def\bart{{}_{\bar{\scriptscriptstyle(2)}}}
\def\barth{{}_{\bar{\scriptscriptstyle(3)}}}
\def\barfo{{}_{\bar{\scriptscriptstyle(4)}}}
\def\barfiv{{}_{\bar{\scriptscriptstyle(5)}}}
\def\barsix{{}_{\bar{\scriptscriptstyle(6)}}}

\def\umo{{{}^{\scriptscriptstyle-(1)}}}
\def\umt{{{}^{\scriptscriptstyle-(2)}}}

\def\Bo{{{}_{\und{\scriptscriptstyle(1)}}}}
\def\Bt{{{}_{\und{\scriptscriptstyle(2)}}}}
\def\Bth{{{}_{\und{\scriptscriptstyle(3)}}}}
\def\Bfo{{{}_{\und{\scriptscriptstyle(4)}}}}
\def\Bfiv{{{}_{\und{\scriptscriptstyle(5)}}}}
\def\Bsix{{{}_{\und{\scriptscriptstyle(6)}}}}

\def\opBo{{{}_{\und{\scriptscriptstyle(1)^{\rm op}}}}}
\def\opBt{{{}_{\und{\scriptscriptstyle(2)^{\rm op}}}}}

\def\text#1{\mbox{\rm #1}}
\def\note#1{}

\def\blacksquare{{\lform}}%AMS Tex Fakes
\def\frac#1#2{{{#1\over#2}}}

\def\proof{\goodbreak\noindent{\bf Proof\quad}}

\def\endproof{{\ $\lform$}\bigskip }

\def\align#1{\begin{eqnarray*}#1\end{eqnarray*}}
\def\alignn#1#2{\begin{eqnarray}\label{#1}#2
\end{eqnarray}}
\def\cmath#1{\[\begin{array}{c} #1 \end{array}\]}
\def\ceqn#1#2{\begin{equation}\label{#1}\begin{array}{c}#2
\end{array}\end{equation}}

\def\vect{{\bf t}}

\def\vecl{{\bf l}}

\def\vecm{{\bf m}}\def\vece{{\bf e}}
\def\vecf{{\bf f}}

\def\Ro#1{{\CR_{#1}\uo}}
\def\Rt#1{{\CR_{#1}\ut}}
\def\Rmo#1{{\CR_{#1}\umo}}
\def\Rmt#1{{\CR_{#1}\umt}}
\def\CH{{\bar C\lbiprod \bar H}}
\def\HB{{H\rbiprod B}}

\documentstyle[11pt,epsf]{article}
\def\nosum{}

\newtheorem{lemma}{Lemma}[section] \newtheorem{propos}[lemma]{Proposition}
\newtheorem{example}[lemma]{Example} \newtheorem{theorem}[lemma]{Theorem}
 \newtheorem{corol}[lemma]{Corollary}
\newtheorem{defin}[lemma]{Definition} 
\newtheorem{rem}[lemma]{Remark}

\begin{document}\baselineskip 20pt

{\ } \hskip 5in DAMTP/95-57
\vspace{.2in}

\begin{center} {\LARGE  DOUBLE-BOSONISATION OF BRAIDED GROUPS AND THE
CONSTRUCTION OF {$U_q(\cg)$}}
\\
\baselineskip 12pt
{\ } {\ }\\ S.  Majid\footnote{Royal Society University Research Fellow and
Fellow of Pembroke College, Cambridge}\\{\ }\\
Department of Mathematics, Harvard University\\
1 Oxford Street, Cambridge, MA02138, USA\footnote{During the calendar years
1995+1996}\\
+\\
Department of Applied
Mathematics \& Theoretical Physics\\ University of Cambridge, Cambridge CB3
9EW, UK
\end{center}
\begin{center}
October, 1995
\end{center}

\vspace{5pt}

\begin{quote}\baselineskip 12pt
\noindent{\bf Abstract} We introduce a quasitriangular Hopf algebra or `quantum
group' $U(B)$, the {\em double-bosonisation}, associated to every braided group
$B$ in the category of $H$-modules over a quasitriangular Hopf algebra $H$,
such that $B$ appears as the `positive root space', $H$ as the `Cartan
subalgebra' and the dual braided group $B^*$ as the `negative root space' of
$U(B)$. The choice $B=f$ recovers Lusztig's
construction of $U_q(\cg)$, where $f$ is Lusztig's algebra associated to a
Cartan datum; other choices give more novel quantum groups.  As an application,
our construction  provides a canonical way of building up quantum groups from
smaller ones by repeatedly extending their positive and negative root spaces by
linear braided groups; we explicitly construct $U_q(sl_3)$ from $U_q(sl_2)$ by
this method, extending it by the quantum-braided plane $\A_q^2$. We provide a
fundamental representation of $U(B)$ in $B$. A projection from the quantum
double, a theory of double biproducts and a Tannaka-Krein reconstruction point
of view are also provided.

\bigskip
\end{quote}
\baselineskip 5pt
\tableofcontents

\section{Introduction}
\baselineskip 22pt

The theory of braided groups or Hopf algebras in braided categories has been
introduced by the author in 1989-1990 \cite{Ma:bra}\cite{Ma:bg} as
a more fundamental object underlying the theory of quantum groups. Braided
planes, lines, matrices, Lie algebras, differentials and other constructions
are now known in this braided-geometrical setting, developed in a series of
papers by the author and collaborators; see \cite{Ma:introm}\cite{Ma:introp}
for reviews. The main idea of braided groups is that they are like Hopf
algebras,
with a diagonal or coproduct map $B\to B\tens B$, but the tensor product here
is
not the usual commutative one; rather, it is a braided non-commutative tensor
product. This `outer noncommutativity' between two algebras is quite a
different
foundation for `braided geometry' from the usual conception of non-commutative
geometry based on the idea of a single `co-ordinate algebra' becoming
non-commutative. It is intended instead as a generalisation of supergeometry.

The reason that braided groups provide the foundation of a kind of $q$-deformed
algebraic geometry is quite fundamental; It can be expected that
the starting point of such a geometry should be the additive properties of
$\R^n$, which means an additive or linear coproduct $\Delta b=b\tens 1+1\tens
b$ on suitable generators. Such a coproduct is not interesting for an ordinary
Hopf algebra since, being cocommutative, it belongs essentially to an
enveloping algebra and not to a quantum group. However, for braided groups such
a coproduct {\em is} compatible with $q$-deformation and, indeed, familiar
examples
such as the so-called quantum plane with relations $yx=qxy$ have exactly such a
linear coproduct or `coaddition'\cite{Ma:poi}.

Apart from large classes of examples, there are also theorems (due to the
author)  which relate braided groups in certain braided categories to ordinary
quantum groups\cite{Ma:tra}\cite{Ma:bos}. This allows braided-group
constructions to be used to obtain results about ordinary quantum groups. So
far, an important application has been to the
construction\cite{Ma:poi}\cite{Ma:qsta} of inhomogeneous quantum groups by
`bosonisation' of linear braided groups. The braided group appears as the
`linear' part of the inhomogeneous quantum group. The bosonisation
construction\cite{Ma:bos}\cite{Ma:skl} associates to every braided group $B$ in
the category of representations of a quasitriangular Hopf algebra $H$ an
ordinary Hopf algebra $B\lbiprod H$. We recall the basic theory in the
Preliminaries section~2.

The present paper extends this close connection between quantum groups and
braided groups
with a new construction $U(B)$ associated to the same data. This is a
quasitriangular
Hopf algebra in which our previous bosonisation $H\rbiprod B$ appears as
`positive Borel subalgebra' and the bosonisation $B^*\lbiprod H$ of the dual of
$B$ appears as `negative subalgebra'. The new part is the nontrivial cross
relations between the these two sub-Hopf algebras within $U(B)$. The background
quantum group $H$ (in the category of representations of which $B,B^*$ live)
plays the role of `Cartan subalgebra'.
We introduce this construction in Section~3 and develop some of its basic
properties. The triangular decomposition
of $U(B)$ into $B^*, H, B$ is an intrinsic feature of our constructive
definition. The most novel aspect is that apart from the general nature of our
construction, the braided groups $B^*$ and $B$ need not be isomorphic; the
positive and negative `roots' are in general dual to each other rather than
isomorphic.

Independently of the author's development of braided groups and bosonisation,
G.~Lusztig in \cite{Lus} introduced a novel construction of the
quantum enveloping algebras $U_q(\cg)$ of V.G.~Drinfeld and
M.~Jimbo\cite{Dri}\cite{Jim:dif} associated
to complex semisimple Lie algebras $\cg$. Although
Lusztig does not use the formalism of braided groups, it is obvious that his
algebra $f$ with `coproduct' $r:f\to f\tens f$ could be viewed as an example of
the $q$-braided type associated to a bilinear form, and that the resulting
quantum Borel subalgebra $U_q(b_+)$ could be viewed as its bosonisation. This
is clear by comparison with the corresponding physics literature\cite{Ma:csta}
where such $q$-braided groups and their bosonisation were studied in a quite
different (physical) context. Such a view on Lusztig's approach to $U_q(b_+)$
has been pointed out most recently by M. Schauenberg at the Chicago AMS meeting
in March 1995, and is one of the motivations for our new $U(B)$ construction.
In fact, we need something stronger, namely that $f$ lives in the category
generated by a certain {\em weakly quasitriangular} Hopf algebra, which we
introduce. We are then able to cast Lusztig's construction for all of
$U_q(\cg)$ into a braided setting.

In fact, we still require Lusztig's elegant construction of $f$ as the
coradical of a bilinear pairing induced by the
Cartan matrix datum, which provides the $q$-Serre relations of $U_q(\cg)$ in
his approach.  But once we are given this as a (self-dual) braided group in a
certain braided category, we can simply feed $B=f=B^*$ into the abstract $U(B)$
construction in the present paper and recover $U_q(\cg)$
directly without the explicit proofs and calculations in \cite{Lus}. This is
demonstrated in detail in Section~4 and provides, we believe, a useful abstract
setting for Lusztig's approach. The fundamental Verma module representation in
\cite{Lus} is recovered now (in an adjoint form) as a natural construction for
an action of $U(B)$ on $B$ by `braided differentiation action'\cite{Ma:fre} of
$B^*$ and the `braided adjoint action'\cite{Ma:lie} of $B$.

The $U(B)$ construction is also more general. We demonstrate an application of
this in Section~5,
where we begin with a central extension of the quantum group $U_q(sl_2)$ in the
role of `Cartan subalgebra'
and adjoin the so-called quantum $B=\A_q^2$ to the positive root elements (and
another copy to the negative roots) by means of our $U(B)$ construction. The
result is $U_q(sl_3)$, but constructed now in a novel way. We see that it is
naturally represented on the quantum plane $\A_q^2$ by braided differentiation,
the braided adjoint action and $U_q(su_2)$ rotations. The same construction
works for $U(V(R))$, where $V(R)$ is the linear braided group associated to
suitable R-matrix data in $M_n\tens M_n$. It adds one to the rank of the
quantum group and $n$ to the positive and negative roots. Some other, more
physical, examples will be computed elsewhere as a construction of
$q$-deformed conformal groups.

A second motivation for the $U(B)$ construction is Drinfeld's quantum double
$D(H)$ in \cite{Dri}. Many constructions for Hopf algebras generalise easily
(in a diagrammatic notation) to the setting of braided groups, so it is natural
to ask for a braided-group version of Drinfeld's double. While a braided group
double cross product theory does exist, the example of the Drinfeld double
(based on mutual coadjoint actions) appears to become `tangled
up'; i.e., it works fine in a symmetric monoidal category but it encounters
problems in a truly braided one. This is also true for even some of the
simplest quantum group constructions, such as tensor products of braided
groups. To overcome this problem we use the bosonisation procedure; since the
bosonisation $B\lbiprod H$ is equivalent in a certain categorical sense to $B$
(the ordinary modules of the former are in correspondence with the braided
modules of the latter\cite{Ma:bos}), one can expect that the double of a
braided group $B$, if it exists, should be closely related to the Drinfeld
double of the bosonisation. We compute the latter in Appendix~A and show that
it projects onto $U(B)$. Hence $U(B)$ could be regarded as some kind of bosonic
(i.e. not braided) version of the `braided double' of $B$; it reduces to
Drinfeld's quantum double when $H=k$. The projection is also the means by which
the quasitriangular structure of $U(B)$ verified directly in Section~3, can be
deduced, which is in analogy with the way that the quasitriangular structure of
$U_q(\cg)$ is obtained from the quantum double of the Borel
subalgebra\cite{Dri}.

In Appendix~B  we show that  $U(B)$ can be viewed as special cases of a more
general `double-biproduct' construction, which we also introduce. Single
bosonisations can be viewed as examples or single biproducts in the sense of
\cite{Rad:str}, so this generalisation is a natural question. However, the
double bosonisations remain the main examples of interest and their key
properties do not come from this point of view. Finally, Appendix~C provides a
still different way of thinking about the complicated relations of the Hopf
algebra $U$, namely as obtained by Tannaka-Krein reconstruction from a suitable
category of braided crossed $B-B^*$-bimodules. One can also think of the latter
as braided crossed $B$-modules in the sense of \cite{Bes:cro}\cite{Dra:bos}.
These are sufficiently complicated, however, that this is not a very convenient
way to prove that $U$ is a Hopf algebra, but provides an alternative viewpoint.

We work over a ground field $k$. With a little care one can work over a
commutative ring just as well. We also note that our constructions will not be
limited to finite-dimensional Hopf algebras.

\subsection*{Acknowledgements} I would like to thank Arkadiy Berenstein for
some useful discussions. The main part of the research was completed during my
visit to Kyoto with funding under the
joint Research Institute of Mathematical Sciences -- Isaac Newton Institute
programme; I thank my hosts there for support.

\section{Preliminaries}

Here we collect basic facts and notation from the theory of braided groups and
their bosonisation, needed in Section~3 and the Appendix. For a more detailed
review, see \cite{Ma:introm}. We also recall Lusztig's algebra $f$ which is
needed in Section~4. We begin with quantum groups in the sense of V.G.
Drinfeld.

1. A {\em quasitriangular} Hopf algebra is $(H,\Delta,\eps,S,\CR)$ where $H$ is
a unital algebra, $\Delta:H\to H\tens H$ and $\eps:H\to k$ are algebra
homomorphisms
forming a coalgebra. This defines a bialgebra. In addition, $S$ is the
convolution-inverse of the identity $H\to H$, i.e. characterised by $\nosum h\o
S h\t=\eps(h)=\nosum (Sh\o)h\t$, where we use the Sweedler
notation\cite{Swe:hop} $\Delta h=\nosum h\o\tens h\t$ (summation understood).
This defines a Hopf algebra. Finally, $\CR\in H\tens H$ is invertible and
obeys\cite{Dri}
\ceqn{qua}{(\Delta\tens\id)(\CR)=\CR_{13}\CR_{23},\quad
(\id\tens\Delta)(\CR)=\CR_{13}\CR_{12}\\
 \tau\circ\Delta=\CR(\Delta\ )\CR^{-1},}
where $\CR_{12}=\CR\tens 1\in H^{\tens 3}$ etc., and where $\tau$ is the usual
transposition map.

A {\em dual quasitriangular} bialgebra or Hopf algebra is $(A,\CR)$ where $A$
is a bialgebra or Hopf algebra and $\CR:A\tens A\to k$ is a
convolution-invertible linear map obeying the obvious dualisation of
(\ref{qua}), namely
\ceqn{dqua}{\CR\circ(\cdot\tens\id)=\CR_{13}*\CR_{23},\quad
\CR\circ(\id\tens\cdot)=\CR_{13}*\CR_{12}\\
\cdot\circ\tau=\CR*\cdot*\CR^{-1}}
in the convolution algebras $\hom(A\tens A\tens A,k)$ and $\hom(A\tens A,A)$
respectively\cite{Ma:pro}\cite{Ma:tan}.

In between these two formulations of Drinfeld's ideas is an intermediate
one\cite{Ma:mor} called a {\em weakly quasitriangular dual pair}. This is a
pair $(H,A)$ of Hopf algebras equipped with a duality pairing $\<\ ,\ \>:H\tens
A\to k$ and convolution-invertible algebra/anti-coalgebra maps
$\CR,\bar{\CR}:A\to H$ obeying
\ceqn{wqua}{ \<\bar{\CR}(a),b\>=\<{\CR}^{-1}(b),a\>,\quad\forall a,b\in A,\quad
\del^R h=\CR*(\del^L h)*\CR^{-1},\quad \del^R h=\bar{\CR}*(\del^L
h)*\bar{\CR}^{-1}}
for all $h\in H$. Here $*$ is the convolution product in $\hom(A,H)$ and
$(\del^L h)(a)=\nosum \<h\o,a\>h\t$, $(\del^R h)(a)=\nosum h\o\<h\t,a\>$ are
left and right `differentiation operators' regarded as maps $A\to H$ for fixed
$h$. It is evident that given a dual pair of bialgebras or Hopf algebras,
\[ {\rm quasitriangularity}\Rightarrow{\rm weak\
quasitriangularity}\Rightarrow{\rm dual\ quasitriangularity}\]
by the appropriate evaluation using the duality pairing. It was shown by the
author in 1989\cite{Ma:qua}\cite{Ma:mor} that the dually paired bialgebras
$(A(R),\widetilde{U(R)})$ associated to a matrix solution $R\in M_n\tens M_n$
of the {\em Yang-Baxter equations} $R_{12}R_{13}R_{23}=R_{23}R_{13}R_{12}$ are
weakly quasitriangular, which means that $A(R)$ is dual quasitriangular. Here
$M_n$ denotes $n\times n$ matrices and $R_{12}=R\tens\id\in M_n^{\tens 3}$,
etc.

2. If $H$ is a quasitriangular bialgebra or Hopf algebra then the categories
${}_H\CM$, $\CM_H$ of left modules and right modules are each braided. This
means that for every two objects $V,W$ there are functorial isomorphisms
$\Psi_{V,W}:V\tens W\to W\tens V$ which behave appropriately under $\tens$ of
modules. Explicitly, the braidings for left, right modules are
\eqn{lrbraid}{ \Psi_{V,W}(v\tens w)=\nosum \Rt{}\la w\tens \Ro{}\la v, \quad
\Psi_{V,W}(v\tens w)=\nosum w\ra\Ro{}\tens v\ra \Rt{}}
for all $v\in V, w\in W$. Here $\CR=\nosum\Ro{}\tens\Rt{}$ is a notation for
explicit components of $\CR\in H\tens H$ (summation understood) and $\la,\ra$
refer to left, right actions respectively.

It is a trivial matter to recast these formulae for the cases of weakly
quasitriangular and dual quasitriangular bialgebras or Hopf algebras.  Then the
categories ${}^A\CM$, $\CM^A$ of left, right comodules become
braided\cite{Ma:pro}\cite{Ma:tan}.

3. An algebra in the category $\CM_H$ means a (right) $H$-module algebra, i.e.
an algebra for which the structure maps intertwine (are covariant under) the
action of $H$. A first result of the theory of braided groups is the
observation (due to the author) that if $B,C$ are such module algebras in a
braided category (i.e. if $H$ is quasitriangular) then there is an associative
algebra $B\und\tens C$, the {\em braided tensor product algebra}\cite{Ma:bg}
again in the category. The product rule is
\[ (b\tens c)(d\tens e)=b\cdot\Psi_{C,B}(c\tens d)\cdot e\]
where the output of $\Psi$ is multiplied from the left by $b$ and from the
right by $e$. A bialgebra in the category means a unital algebra $B$ equipped
with algebra homomorphisms $\und\Delta:B\to B\und\tens B$, $\und\eps:B\to k$
forming a coalgebra and intertwining the action of $H$. A Hopf algebra means
that in addition there is an intertwiner $\und S:B\to B$ which is the
convolution inverse of the identity, i.e. $\nosum (\und
Sb\Bo)b\Bt=\und\eps(b)=\nosum b\Bo \und S b\Bt$, where $\und \Delta=\nosum
b\Bo\tens b\Bt$ denotes the braided coproduct. We use the term {\em braided
group} to denote bialgebras or Hopf algebras in a braided category.  They have
been introduced and studied by the
author \cite{Ma:bg}\cite{Ma:tra}\cite{Ma:bos}, where basic
properties such as
\eqn{bantimul}{ \und S(bc)=\cdot\circ \Psi_{B,C}(\und S b\tens \und S c)}
are proven. There is a also a further theory of quasitriangular braided groups
(or braided quantum groups) which we do not
need here; see \cite{Ma:tra}\cite{Ma:bg}. The theory works in any braided
category and we can easily read off the particular formulae for the other cases
${}_H\CM$, ${}^A\CM$ and $\CM_A$ of interest.

Two braided groups $C,B$ are said to be dually paired if there is an
intertwiner $\ev:C\tens B\to k$ such that
$\ev(cd,b)=\nosum\ev(d,b\Bo)\ev(c,b\Bt)$ and
$\ev(c,ab)=\nosum\ev(c\Bt,a)\ev(c\Bo,b)$
hold for all $a,b\in B$ and $c,d\in C$. This is the natural {\em categorical
duality} pairing\cite{Ma:introm}. In the finite-dimensional non-degenerate case
we write $C=B^\star$. In applications where we are finally interested in
ordinary Hopf algebras, it is also useful to consider an {\em ordinary duality
pairing} $\<\ ,\ \>$ between braided groups $C,B$ defined in the more usual way
(without reversing the product or coproduct as for the categorical $\ev$). In
this case, if $B$ lives in $\CM_H$ then $C$ lives naturally in ${}_H\CM$ and
$\<h\la c,b\>=\<c,b\ra h\>$ for all $h\in H$ is the appropriate `covariance'
condition. In the finite-dimensional non-degenerate case we write $C=B^*$.

4. If $B$ is a braided group in ${}_H\CM$ then its {\em bosonisation} is the
Hopf algebra $B\lbiprod H$ defined as $B\tens H$ with product, coproduct and
antipode\cite{Ma:bos}
\ceqn{lbos}{ (b\tens h)(c\tens g)=\nosum b h\o\la c\tens h\t g,\quad
\Delta(b\tens h)=\nosum b\Bo\tens \CR\ut h\o\tens \CR\uo\la b\Bt\tens h\t\\
S(b\tens h)=\nosum (Sh\t)\cu\CR\uo\la\und Sb\tens S(\CR\ut h\o);\quad
\cu\equiv\nosum(S\CR\ut)\CR\uo}
and tensor product unit and counit. The algebra structure is a smash product
while the coalgebra is a smash coproduct by a particular coaction $b\mapsto
\CR_{21}\la b$ induced by the quasitriangular structure\cite{Ma:dou}.
Bosonisations can be viewed\cite{Ma:skl} as a particular class of
biproducts\cite{Rad:str} (this class was not considered in any sense in
\cite{Rad:str}, however). The right-handed version for $B\in \CM_H$ is
$H\rbiprod B$ defined by
\ceqn{rbos}{ (h\tens b)(g\tens c)=\nosum h g\o\tens (b\ra g\t)c,\quad
\Delta(h\tens b)=\nosum h\o\tens b\Bo\ra\CR\uo\tens h\t\CR\ut\tens b\Bt\\
S(h\tens b)=\nosum S(h\t \Rt{})\tens \und Sb\ra\Ro{}\cv Sh\o,\quad
\cv=\nosum\Ro{}S\Rt{}.}

The corresponding formulae for bosonisation in the comodule categories
${}^A\CM$, $\CM^A$ are trivially obtained by the usual conversions of the
module formulae; $A\rbiprod B$ and $B\lbiprod A$ have the smash coproduct
coalgebra by the given coaction of $A$ and smash product by the coaction
induced by $\CR$. For example, the first case (for $B\in{}^A\CM$) and its
duality pairing with $C\lbiprod H$ is explicitly\cite{Ma:mec}\cite{Ma:qsta}
\ceqn{rcobos}{ (a\tens b)(d\tens c)=ad\o\tens b\bo c\<\CR(b\bt),d\t\>,\quad
\Delta(a\tens b)=a\o\tens b\Bo\bo\tens a\t b\Bo\bt\tens b\Bt\\
 \<c\tens h,a\tens b\>=\<\Rmt{}h\o,a\>\ev(\und S^{-1}c,\Rmo{}h\t\la b),\quad
\<S(c\tens h),a\tens b\>=\<Sh,a\>\ev(c,b).}
Moreover, if $(H,A)$ is weakly quasitriangular then we can adapt the
bosonisation formulae (\ref{lbos})--(\ref{rbos}) to $B\in \CM^A$, ${}^A\CM$
without dualisation; we  define the action of $H$ on $B$ by evaluation against
the given coaction of $A$ and also replace $\CR_{21}\la b$ by $\CR$ evaluated
against the given coaction. These are all variants of the bosonisation
construction in \cite{Ma:bos} in one form or another.

5. Let $R\in M_n\tens M_n$ obey the Yang-Baxter equations. We
define\cite{Ma:fre} the {\em free braided vector} algebra $V(R)$ to be the
braided group with
\ceqn{freevec}{  B=k\<e^i|i=1,\cdots ,n\>,\quad \und\Delta e^i=e^i\tens
1+1\tens e^i,\quad \eps e^i=0\\
\und S e^i=-e^i,\quad \Psi(e^i\tens e^j)=\sum_{a,b} R^j{}_a{}^i{}_b e^a\tens
e^b.}
in the braided category of (left) $A(R)$-comodules. Here $A(R)$ can be replaced
by
any dual-quasitriangular bialgebra or Hopf algebra, provided it induces the
same braiding.

According to the braided-geometrical point of view (where $e^i$ are like
co-ordinates on a vector space), there is also a dually-paired {\em braided
covector algebra} $\Vhaj(R)$\cite{Ma:fre}. The version of it which has an
ordinary duality pairing with $V(R)$ is
\ceqn{freecov}{ D=k\<f_i|i=1,\cdots, n\>,\quad \und\Delta f_i=f_i\tens 1+1\tens
f_i,\quad \und\eps f_i=0\\
\und S f_i=-f_i,\quad \Psi(f_i\tens f_j)=\sum_{a,b} f_b\tens f_a
R^a{}_i{}^b{}_j}
in the braided category of right $A(R)$-comodules. The pairing is defined by
$\<f_j,e^i\>=\delta_j{}^i$, where $\delta$ is the Kronecker delta function. The
categorical dual $C$ has the opposite product and coproduct; it looks the same
on generators but has $\Psi(f_i\tens f_j)=\sum_{a,b} f_a\tens f_b
R^a{}_j{}^b{}_i$. According to (a version of) the theory in \cite{Ma:fre} there
is a left action of $C$ on $V(R)$ by {\em braided differentiation}
\eqn{bdiff}{\del_i:V(R)\to V(R),\quad \del_i (e^{i_1}\cdots e^{i_m})=
[m;R]^{i_1\cdots i_m}_{ij_2\cdots j_m}e^{j_2}\cdots e^{j_m},}
where
\eqn{braint}{ [m;R]=\id+(PR)_{12}+(PR)_{23}(PR)_{12}+\cdots (PR)_{m-1\
m}\cdots(PR)_{12}\in M_n^{\tens m}}
is the {\em braided integer matrix}\cite{Ma:fre}. Here $P$ denotes the
permutation matrix in $M_n\tens M_n$ and $(PR)_{12}=PR\tens\id$, etc. Then the
categorical pairing between $C$ and $B$ is $\ev(c(f),b)=\eps(c(\del) b)$.
Equivalently, the ordinary duality pairing between $D$ and $B$ is
\eqn{DBpair}{ \<f_{j_1}\cdots f_{j_m},e^{i_1}\cdots e^{i_r}\>=\delta^r_m
[m;R]!^{i_1\cdots i_m}_{j_1\cdots j_m},\quad [m;R]!=[m;R]_{1\cdots
m}[m-1;R]_{2\cdots m}\cdots[2;R]_{m-1\ m}}
where the numerical suffices denote as usual the embedding of our matrices in
the corresponding positions in $M_n^{\tens m}$. Equivalently, $\del_i b$ is
characterised by $\und\Delta b=e^i\tens\del_i b+$ terms which do not consist of
$e^i$ in the first tensor factor. There are also right braided differentials
$\overleftarrow{\del_i}$.

Various properties are proven about such braided differential operators
$V(R)\to V(R)$ in \cite{Ma:fre}. In particular, we introduce  a corresponding
exponential
\eqn{BDexp}{ \exp_R=\sum_m e^{i_1}\cdots e^{i_m}\tens f_{j_1}\cdots f_{j_m}
([m;R]!^{-1})^{j_1\cdots j_m}_{i_1\cdots i_m} }
as a formal power-series in the matrix entries of $R$  with coefficients in
$B\tens D$.  This assumes that the pairing is non-degenerate, i.e. that the
$[m,R]!$ matrices are all invertible. In many applications, such as the braided
Taylor's theorem\cite{Ma:fre}, only a finite number of terms from $\exp_R$
contribute.

When the pairing $\<\ ,\ \>$ is degenerate, one can add to further relations to
both sides until it becomes non-degenerate.  In
particular\cite{Ma:fre} if there is a second matrix $R'\in M_n\tens M_n$
obeying certain relations with $R$, one can add to $B,D$ the quadratic
relations
\eqn{BDrel}{ e^ie^j=\sum_{a,b} R'{}^j{}_a{}^i{}_b e^a e^b,\quad
f_if_j=\sum_{a,b} f_b  f_a R'{}^a{}_i{}^b{}_j}
to $B$ and $D$. If both $R,R'$ are $q$-deformations of the identity matrix then
by genericity arguments we will know that the resulting quotients $V(R',R)$ and
$\Vhaj(R',R)$ are non-degenerately paired. This is the case for the $sl_n$
quantum planes $\A_q^n$ in Section~5, where $R'=q^{-2}R$. The $\del_i$ descend
to these quotients\cite{Ma:fre}.

Lusztig's algebra $f$ in \cite{Lus} can  be viewed as a braided group quotient
of the free braided plane $V(R)$ with $R^i{}_j{}^k{}_l=q^{i\cdot
k}\delta^i{}_j\delta^k{}_l$, where $0\ne q\in k$ and $i\cdot j$ are the
components of a bilinear form over $\Z$. In Lusztig's example the braided
differentiation operators are denoted ${}_ir$ and $r_i$, and the radical of the
pairing is far from quadratic. The braided coproduct $\und\Delta$ is denoted
$r$ in \cite{Lus}. The braided antipode $\und S$ does not appear explicitly.
See Section~4 for further details on how to view Lusztig's algebra $f$ as a
braided group.

\section{Quasitriangular Hopf algebra {$U(B)$} associated to a braided group}

Let $H$ be a quasitriangular Hopf algebra and $B$ a braided group in the
category $\CM_H$. Let $C$ be a braided group in $\CM_H$ which is dual to $B$ in
the sense of an categorical duality pairing $\ev$. Equivalently,   $D=C^{\rm
op/op}$ is a Hopf algebra in ${}_H\CM$ which is dual to $B$ in the sense of an
ordinary duality pairing $\<\ ,\ \>$ which is $H$-bicovariant as explained
above. We construct in this section an associated quasitriangular Hopf algebra
$U(B)$. We suppose that $C$ (or, equivalently, $D$) has invertible
braided-antipode. We denote by $\bar H$ the same Hopf algebra as $H$ but
equipped with the quasitriangular structure $\bar{\CR}=\CR_{21}^{-1}$.

\begin{lemma}\cite{Ma:introm} If $D$ is a braided group in ${}_H\CM$ then
$D^{\und{\rm cop}}$ defined by the same algebra, unit and counit as $D$ and the
opposite coproduct $\Psi_{D,D}^{-1}\circ\und\Delta$ and antipode $\und S^{-1}$,
is a braided group in ${}_{\bar H}\CM$.
\end{lemma}
\proof A diagrammatic proof is given in \cite{Ma:introm}. The lemma is also
easily checked directly from the axioms (and is in fact the reason that the
naive concept of braided-cocommutativity for braided groups
$\und\Delta=\Psi^{-1}\circ\und\Delta$ does not make sense; see \cite{Ma:bra}).
\endproof

We see that $\bar C\equiv D^{\und{\rm cop}}=(C^{\rm op/op})^{\und{\rm cop}}$ is
a Hopf algebra in ${}_{\bar H}\CM$. We denote its coproduct explicitly by
$\bar\Delta=\nosum c\baro\tens c\bart$ and its antipode by $\bar S$. According
to the bosonisation theory recalled in the preliminaries, we immediately have
two Hopf algebras $H\rbiprod B$ and $\bar C\lbiprod \bar H$. $H$ itself is a
sub-Hopf algebra of each.

\begin{theorem} There is a unique Hopf algebra structure $U=U(\bar C,H,B)$ on
$\bar C\tens H\tens B$ such that $H\rbiprod B$ and $\bar C\lbiprod\bar H$ are
sub-Hopf algebras by the canonical inclusions and
\[ bc=\nosum ( \Rt1\la c\bart) \Rt2 \Rmo1 (b\Bt\ra\Rmo2)\<\Ro1\la
c\baro,b\Bo\ra\Ro2\>\<\Rmt1\la \bar S c\barth,b\Bth\ra\Rmt2\>\]
for all $b\in B, c\in \bar C$ viewed in $U$. Here $\CR_1,\CR_2$ etc. are
distinct copies of the quasitriangular structure $\CR$ of $H$.
\end{theorem}

We will prove this via a series of lemmas. We begin by proving associativity.
Note first that if there is an associative product as stated, then it is given
uniquely by
\alignn{Uprod}{ (c\tens h\tens b)\cdot(d\tens g\tens a)\equad&&=\nosum c(
h\o\Rt1\la d\bart)\tens h\t \Rt2 \Rmo1 g\o\tens (b\Bt\ra\Rmo2 g\t)a\nonumber\\
&&\qquad\qquad \<\Ro1\la d\baro,b\Bo\ra\Ro2\>\<\Rmt1\la\bar S
d\barth,b\Bth\ra\Rmt2\>}
for all $c,d\in \bar C$, $a,b\in B$ and $h,g\in H$: Because $\bar C\rbiprod
\bar H$ and $H\lbiprod B$ are subalgebras, we know that a general product has
the form
\eqn{Uprodform}{(chb)\cdot(dga)=\sum c(hd_i)R_i(b_ig)a= \sum c(h\o\la d_i)
h\t R_i g\o (b_i\ra g\t)a}
if $bd=\sum d_i R_i b_i$ say, where $d_i\in \bar C$, $R_i\in H$ and $b_i\in B$
are all viewed in $U$ in the canonical way. We take the right hand side as the
definition of the product of general elements.

\begin{lemma}  The map (\ref{Uprod}) is an associative product on $U=\bar
C\tens H\tens B$.
\end{lemma}
\proof It is enough to prove associativity in the special case
$(a\cdot(chb))\cdot d=a\cdot((chb)\cdot d)$ for all $a,b\in B$, $c,d\in \bar C$
and $h\in H$, viewed in $\bar C\tens H\tens B$ in the canonical way. One can
then deduce the general case by breaking down products in the form
(\ref{Uprodform}) and using that $C$ is a left $H$-module algebra and $B$ a
right $H$-module algebra.

To prove the special case, we compute:
\align{&&\equad (a\cdot(chb))\cdot d=\left(\Rt{1}\la c\bart\tens \Rt{2}\Rmo{1}
h\o\tens (a\Bt\ra\Rmt{2}h\t)b\right)\cdot\left(d\tens1\tens1\right)\\
 &&\qquad\quad \<\Ro{1}\la c\baro,a\Bo\ra\Ro{2}\>\<\Rmt{1}\la \bar S
c\barth,a\Bth\ra\Rmt{2}\>\\
&&=(\Rt{1}\la c\bart)\left((\Rt{2}\Rmo{1}h\o)\o\Rt{3}\la d\bart\right)\tens
(\Rt{2}\Rmo{1}h\o)\t \Rt{4}\Rmo{3}\\
&&\ \tens \left((a\Bt\ra\Rmo{2} h\t)b\right)\Bt\ra\Rmo{4} \<\Ro{1}\la
c\baro,a\Bo\ra\Ro{2}\>\<\Rmt{1}\la \bar S c\barth,a\Bth\ra\Rmt{2}\>\\
&&\ \<\Ro{3}\la d\baro, \left((a\Bt\ra\Rmo{2}
h\t)b\right)\Bo\ra\Ro{4}\>\<\Rmt{3}\la\bar S d\barth,\left((a\Bt\ra\Rmo{2}
h\t)b\right)\Bth\ra\Rmt{4}\>}
from the definition (\ref{Uprod}).  Next, we use that $\und\Delta$ and hence
$(\id\tens\und\Delta):B\to B\und\tens B\und\tens B$ is an algebra homomorphism
to the braided tensor product algebra, i.e.
\[(\id\tens\und\Delta)\circ\und\Delta(ab)=a\Bo(b\Bo\ra\Ro{8}\Ro{9})
\tens(a\Bt\ra\Rt{9})(b\Bt\ra\Ro{10})\tens (a\Bth\ra\Rt{8}\Rt{10})b\Bth\]
where $\CR_{8},\cdots,\CR_{10}$ are fresh copies of $\CR$ distinct from others
to be used below. We also use that the product and coproduct of $B$ are
covariant under $H$. The first gives us the action $\ra\Rmt{4}$ etc on products
of elements of $B$, and the second gives us the coproduct
$(\id\tens\und\Delta)\circ\und\Delta(a\ra h)$ etc. We then use the axioms
(\ref{qua}) to convert all coproducts of $\CR$ to  products of $\CR$, suitably
numbered. We use covariance of $\<\ ,\ \>$ to move $\Ro{1}\la$ etc. in its
first input to $\ra\Ro{1}$ etc. in its second input. We arrive by these steps
at the expression:
\align{&&\equad (a.(chb)).d=(\Rt{1}\la c\bart)(\Rt{2}\Rmo{1}h\o
\Rt{3}\Rt{11}\la d\bart)\\
&&\tens \Rt{5}\Rmo{5}h\t \Rt{4}\Rt{10}\Rmo{3}\Rmo{8}\tens (a\Bth\ra \Rmo{6}
h\fo \Rt{9}\Rmo{4}\Rmo{9})(b\Bt\ra\Rmo{10})\\
&&\<c\baro,a\Bo\ra\Ro{5}\Ro{2}\Ro{1}\>\<\bar S
c\barth,a\Bfiv\ra\Rmt{7}\Rmt{6}\Rmt{2}\Rmt{5}\Rmt{1}\>\\
&&\<d\baro,(a\Bt\ra\Rmo{2} h\th
\Ro{4}\Ro{3})(b\Bo\ra\Ro{8}\Ro{9}\Ro{10}\Ro{11})\>\\
&&\<\bar S d\barth,(a\Bfo\ra \Rmo{7} h\fiv
\Rt{8}\Rmt{4}\Rmt{3})(b\Bth\ra\Rmt{10}\Rmt{9}\Rmt{8})\>.}
We next use the braided duality pairing between $B,C$, which between $B,\bar C$
takes the form
\eqn{pairab}{ \<d,ab\>=\<d\bart,a\ra\Rt{}\>\<d\baro,b\ra\Ro{}\>,\quad \<\bar S
d,ab\>=\<\bar S d\baro,a\>\<\bar S d\bart,b\>,\quad\forall d\in\bar C,\ a,b\in
B.}
Hence,
\align{&&\equad (a.(chb)).d=(\Rt{1}\la c\bart)(\Rt{2}\Rmo{1}h\o
\Rt{3}\Rt{11}\la d\barth)\\
&&\tens \Rt{5}\Rmo{5}h\t \Rt{4}\Rt{10}\Rmo{3}\Rmo{8}\tens (a\Bth\ra \Rmo{6}
h\fo \Rt{9}\Rmo{4}\Rmo{9})(b\Bt\ra\Rmo{10})\\
&&\<c\baro,a\Bo\ra\Ro{5}\Ro{2}\Ro{1}\>\<\bar S
c\barth,a\Bfiv\ra\Rmt{7}\Rmt{6}\Rmt{2}\Rmt{5}\Rmt{1}\>\\
&&\<d\bart,a\Bt\ra\Rmo{2} h\th
\Ro{4}\Ro{3}\Rt{12}\>\<d\baro,b\Bo\ra\Ro{8}\Ro{9}\Ro{10}\Ro{11}\Ro{12}\>\\
&&\<\bar S d\barfo,a\Bfo\ra \Rmo{7} h\fiv \Rt{8}\Rmt{4}\Rmt{3}\>\<\bar S
d\barfiv,b\Bth\ra\Rmt{10}\Rmt{9}\Rmt{8}\>.}

On the other side, we similarly compute
\align{&&\equad a\cdot((chb)\cdot d)=\left(1\tens 1\tens
a\right)\cdot\left(c(h\o\Rt{1}\la d\bart)\tens h\t\Rt{2}\Rmo{1}\tens
b\Bt\ra\Rmo{2}\right)\\
&&\ \<\Ro{1}\la d\baro, b\Bo\ra\Ro{2}\>\<\Rmt{1}\la\bar S d\barth,
b\Bth\ra\Rmt{2}\>\\
&&=\Rt{3}\la\left(c(h\o \Rt{1}\la d\bart)\right)\bart\tens \Rt{4}\Rmo{3}(h\t
\Rt{2}\Rmo{1})\o\\
&&\ \tens \left(a\Bt\ra \Rmo{4}(h\t \Rt{2}
\Rmo{1})\t\right)(b\Bt\ra\Rmo{2})\<\Ro{1}\la d\baro,
b\Bo\ra\Ro{2}\>\<\Rmt{1}\la\bar S d\barth, b\Bth\ra\Rmt{2}\>\\
&&\ \<\Ro{3}\la\left(c(h\o \Rt{1}\la
d\bart)\right)\baro,a\Bo\ra\Ro{4}\>\<\Rmt{3}\la\bar S\left(c(h\o \Rt{1}\la
d\bart)\right)\barth,a\Bth\ra\Rmt{4}\>\\
&&=(\Rt{3}\la c\bart)(\Rt{5}\Rt{10}\Rmo{7}h\t \Rt{6}\la d\barth)\tens
\Rt{4}\Rt{9}\Rmo{3}\Rmo{9}h\fo \Rt{2}\Rmo{1}\\
&&\ \tens (a\Bth\ra \Rmo{4}\Rmo{8} h\fiv \Rt{8}\Rmo{5})(b\Bt\ra\Rmo{2})\\
&&\ \<d\baro,b\Bo\ra\Ro{8}\Ro{2}\Ro{7}\Ro{6}\Ro{1}\>\<\bar S
d\barfiv,b\Bth\ra\Rmt{2}\Rmt{5}\Rmt{1}\>\\
&&\ \<c\baro,a\Bo\ra \Ro{4}\Ro{5}\Ro{3}\>\<d\bart,a\Bt\ra
\Ro{9}\Ro{10}\Rmo{6}h\o\Rt{1}\>\\
&&\ \<\bar S c\barth,a\Bfiv\ra \Rmt{8}\Rmt{9}\Rmt{10}\Rmt{7}\Rmt{6}\>\<\bar S
d\barfo,a\Bfo\ra\Rmt{4}\Rmt{3}
\Rmo{10} h\th \Rt{7}\>,}
using in order: the definition (\ref{Uprod}), the homomorphism property of the
iterated braided-coproduct of $\bar C$ in the form
\[ (\id\tens\bar\Delta)\circ\bar\Delta(cd)=c\baro(\Rmo{11}\Rmo{6}\la
d\baro)\tens (\Rmt{11}\la c\bart)(\Rmo{7}\la d\bart)\tens (\Rmt{7}\Rmt{6}\la
c\barth)d\barth,\]
the covariance of the product and coproduct of $\bar C$ and of $\<\ ,\ \>$, the
quasitriangularity axiom (\ref{qua}) and the evaluation pairing between $B,\bar
C$ in the form
\eqn{paircd}{ \<cd,b\>=\<c,b\Bo\>\<d,b\Bt\>,\quad \<\bar S(cd),b\>=\<\bar S
d,b\Bo\ra \Rmo{}\>\<\bar S c,b\Bt\ra\Rmt{}\>,\quad\forall c,d\in\bar C,\ b\in
B.}

These steps are similar to the proof of associativity of the usual Drinfeld
quantum double, except that now there are several copies of $\CR$ inserted at
various points (arising from the braiding in the categories in which $B,\bar C$
live). It remains to show that these are correctly placed. First, we use the
quantum Yang-Baxter equations for $\CR$ in $H^{\tens 3}$, applied to
$\CR_{10},\CR^{-1}_6,\CR^{-1}_7$. Then we are able to use the
quasicocommutativity axiom (\ref{qua}) in the form $\CR_{10}(h\o\tens
h\t)=(h\t\tens h\o)\CR_{10}$. We make this kind of rearrangement three more
times: we use the QYBE applied to $\CR^{-1}_3,\CR^{-1}_9,\CR^{-1}_{10}$ and
reverse the order of $\CR^{-1}_{3}(h\fo\tens h\th)$; we then use the QYBE
applied to $\CR_9,\CR^{-1}_9,\CR^{-1}_6$ and reverse the order of
$\CR_9(h\t\tens h\th)$; we use the QYBE applied to
$\CR^{-1}_4,\CR^{-1}_8,\CR^{-1}_{10}$ and reverse the order of
$\CR^{-1}_4(h\fiv\tens h\fo)$. The result is
\align{&&\equad a\cdot((chb)\cdot d)=(\Rt{3}\la c\bart)(\Rt{5}\Rmo{7}h\o
\Rt{10}\Rt{6}\la d\barth)\\
&&\ \tens \Rt{4}\Rmo{9}h\th \Rt{9}\Rmo{3}\Rt{2}\Rmo{1} \tens (a\Bth\ra
\Rmo{8}h\fo\Rmo{4}\Rt{8}\Rmo{5})(b\Bt\ra\Rmo{2})\\
&&\ \<d\baro,b\Bo\ra\Ro{8}\Ro{2}\Ro{7}\Ro{6}\Ro{1}\>\<\bar S
d\barfiv,b\Bth\ra\Rmt{2}\Rmt{5}\Rmt{1}\>\\
&&\ \<c\baro,a\o\ra \Ro{4}\Ro{5}\Ro{3}\>\<d\bart,a\Bt\ra
\Rmo{6}h\th\Ro{9}\Ro{10}\Rt{1}\>\\
&&\ \<\bar S c\barth,a\Bfiv\ra \Rmt{10}\Rmt{8}\Rmt{6}\Rmt{9}\Rmt{7}\>\<\bar S
d\barfo,a\Bfo\ra\Rmo{10}h\fiv \Rmt{4}\Rmt{3}\Rt{7}\>.}

This coincides with our previous expression for $(a\cdot(chb))\cdot d$ in all
respects except the order of some copies of $\CR$. We use the QYBE several
times more (i.e. we identify the corresponding braids as generated by the
appropriate braidings $\Psi$) to see that the two expressions are in fact
equal. Hence we have an associative product on $U=\bar C\tens H\tens B$. It is
clear that $(1\tens 1\tens 1)$ a unit element for it. \endproof

Next we note that if the coproducts of $\bar C\lbiprod \bar H$ and $H\rbiprod
B$ extend to a coproduct on $U$, the latter must be given by
\eqn{Ucoprod}{ \Delta(c\tens h\tens b)=\nosum c\baro\tens \Rmo{}h\o\tens
b\Bo\ra\Ro{}\tens \Rmt{}\la c\bart\tens h\t\Rt{}\tens b\Bt}
for all $c\in \bar C$, $b\in B$, $h\in H$. This is $(\Delta(ch))\Delta b$ or
$(\Delta c)\Delta(hb)$ in $U\tens U$ as computed from (\ref{Uprod}).

\begin{lemma} The  map (\ref{Ucoprod}) makes the algebra $U$ into a bialgebra.
\end{lemma}
\proof It is enough to prove that $\Delta(bc)=(\Delta b)\cdot(\Delta c)$ for
all $b\in B$ and $c\in \bar C$, where the product in $U\tens U$ is the usual
tensor product one. After proving this, we compute $\Delta((chb)\cdot(dga))$
using the definition (\ref{Uprodform}) of a general product, and then use that
$\bar C$ is a braided group in the category of left $\bar H$-modules, $B$ a
braided group in the category of right $H$-modules and the
quasi-cocommutativity axiom (\ref{qua}) to obtain $(\Delta c)(\Delta
h)(\Delta(bd))(\Delta g)\Delta(a)$. Since $\bar C\lbiprod \bar H$ and
$H\rbiprod B$ are subalgebras and the product restricted to them is already
known (by the bosonisation theorem) to form a bialgebra, we obtain
$(\Delta(ch))(\Delta b)(\Delta d)(\Delta (ga))=(\Delta(chb))(\Delta(dga))$.

It remains to prove the special case. We compute
\align{&&\equad \Delta(b\cdot c)=\Delta\left( \Rt1\la c\bart \tens \Rt2 \Rmo1
\tens b\Bt\ra\Rmo2\right)\\
&&\qqquad\qquad \<\Ro1\la c\baro,b\Bo\ra\Ro2\>\<\Rmt1\la \bar S
c\barth,b\Bth\ra\Rmt2\>\\
&&=(\Rt{1}\la c\bart)\baro\tens \Rmo{3}(\Rt{2}\Rmo{1})\o\tens
(b\Bt\ra\Rmo{2})\Bo\ra\Ro{3}\\
&&\ \tens \Rmt{3}\la(\Rt{1}\la c\bart)\bart\tens (\Rt{2}\Rmo{1})\t
\Rt{3}\tens(b\Bt\ra\Rmo{2})\Bt\\
&&\qqquad\qquad \<\Ro1\la c\baro,b\Bo\ra\Ro2\>\<\Rmt1\la \bar S
c\barth,b\Bth\ra\Rmt2\>\\
&&=\Rt{1}\la c\bart\tens \Rmo{3}\Rt{2}\Rmo{1}\tens b\Bt\ra\Rmo{2}\Ro{3}\\
&&\quad\tens \Rmt{3}\Rt{4}\la c\barth\tens \Rt{5}\Rmo{4}\Rt{3}\tens
b\Bth\ra\Rmo{5}\\
&&\ \<c\baro,b\Bo\ra\Ro{5}\Ro{2}\Ro{4}\Ro{1}\>\<\bar S c\barfo,
b\Bfo\ra\Rmt{5}\Rmt{2}\Rmt{4}\Rmt{1}\>}
using the product (\ref{Uprod}) and the definition (\ref{Ucoprod}) of the map
$\Delta$. We used the covariance of the coproducts of $\bar C$ and $B$ and
their pairing $\<\ ,\ \>$, and then wrote all coproducts of $\CR$ as products
using (\ref{qua}).

On the other side, we compute
\align{&&\equad (\Delta b)\cdot(\Delta c)=(1\tens 1\tens
b\Bo\ra\Ro{1})\cdot(c\baro\tens\Rmo{1}\tens 1)\tens (1\tens \Rt{1}\tens
b\Bt)\cdot(\Rmt{1}\la c\bart\tens1\tens1)\\
&&=\Rt{2}\la c\baro\bart\tens
\Rt{3}\Rmo{2}\Rmo{1}\o\tens(b\Bo\ra\Ro{1})\Bt\ra\Rmo{3}\Rmo{1}\t\\
&&\ \<\Ro{2}\la c\baro\baro,(b\Bo\ra\Ro{1})\Bo\ra\Ro{3}\>\<\Rmt{2}\la\bar S
c\baro\barth,(b\Bo\ra\Ro{1})\Bth\ra\Rmt{3}\>\\
&&\ \tens \Rt{1}\o \Rt{4}\la(\Rmt{1}\la c\bart)\bart\tens
\Rt{1}\t\Rt{5}\Rmo{4}\tens b\Bt\Bt\ra\Rmo{5}\\
&&\ \<\Ro{4}\la(\Rmt{1}\la c\bart)\baro,b\Bt\Bo\ra\Ro{5}\>\<\Rmt{4}\la\bar
S(\Rmt{1}\la c\bart)\barth, b\Bt\Bth\ra\Rmt{5}\>\\
&&=\Rt{2}\la c\bart\tens \Rt{3}\Rmt{2}\Rmo{1}\Rmo{7}\Rmo{8}\tens
b\Bt\ra\Ro{9}\Ro{7}\Rmo{3}\Rmo{6}\Rmo{9}\Rmo{10}\\
&&\ \tens \Rt{1}\Rt{7}\Rt{8}\Rt{4}\Rmt{9}\Rmt{7}\la c\barfiv\tens
\Rt{6}\Rt{9}\Rt{10}\Rt{5}\Rmo{4}\tens b\Bfiv\ra\Rmo{5}\\
&&\ \<c\baro,b\Bo\ra\Ro{6}\Ro{1}\Ro{3}\Ro{2}\>\<\bar S
c\barsix,b\Bsix\ra\Rmt{5}\Rmt{4}\Rmt{10}\Rmt{8}\>\\
&&\ \<\bar S c\barth,b\Bth\ra\Ro{10}\Ro{8}\Rmt{3}\Rmt{2}\>
\<c\barfo,b\Bfo\ra\Ro{5}\Ro{4}\Rmt{6}\Rmt{1}\>}
using the usual induced product in $U\tens U$, the definition (\ref{Ucoprod})
of the map $\Delta$ and the covariance of the coproducts of $\bar C,B$ and
their pairing $\<\ ,\ \>$. As usual, we write coproducts of $\CR$ as products
via (\ref{qua}). We next use (\ref{qua}) in reverse to recognise the element of
$H\tens H$ acting on $b\Bth\tens b\Bfo$ as in the image of the coproduct if
$H$. Then (\ref{paircd}) tells us that
\align{&&\equad (\Delta b)\cdot(\Delta c)=\Rt{2}\la c\bart\tens
\Rt{3}\Rmo{1}\Rmo{7}\Rmo{8}\tens b\Bt\ra\Ro{9}\Ro{7}\Rmo{6}\Rmo{9}\Rmo{10}\\
&&\ \tens \Rt{1}\Rt{7}\Rt{4}\Rmt{9}\Rmt{7}\la c\barfo\tens
\Rt{6}\Rt{9}\Rt{5}\Rmo{4}\tens b\Bfo\ra\Rmo{5}\\
&&\ \<c\baro,b\Bo\ra\Ro{6}\Ro{1}\Ro{3}\Ro{2}\>\<\bar S
c\barsix,b\Bfiv\ra\Rmt{5}\Rmt{4}\Rmt{10}\Rmt{8}\>\\
&&\ \<(\bar S c\barth)c\barfo,b\Bth\ra\Ro{5}\Ro{4}\Rmt{6}\Rmt{1}\>\\
&&=\Rt{2}\la c\bart\tens \Rt{3}\Rmo{7}\Rmo{8}\tens b\Bt\ra\Ro{9}\Rmo{10}\tens
\Rt{1}\Rmt{7}\la c\barth\tens \Rt{6}\Rt{9}\Rmt{4}\\
&&\ \tens b\Bth\ra\Rmo{5} \<c\baro,b\Bo\ra\Ro{6}\Ro{1}\Ro{3}\Ro{2}\>\<\bar S
c\barfo,b\Bfo\ra\Rmt{5}\Rmt{4}\Rmt{10}\Rmt{8}\>,}
where we use the  axiom for the braided-antipode $\bar S$ of $\bar C$. We also
cancelled $\CR_7\CR^{-1}_9$ from the final expression.

The resulting two expressions differ only by the order of the $\CR$ factors.
They are in fact equal on using the QYBE several times, i.e. equating the
action of two braids generated by the corresponding braidings $\Psi$.

This completes the proof that the map $\Delta:U\to U\tens U$ is an algebra
homomorphism. It is coassociative since this is known for $\bar C\lbiprod \bar
H$ and $H\rbiprod B$ separately. The tensor product of the counits of $\bar
C,H$ and $B$ clearly provides the counit for $\Delta$; hence we have a
bialgebra $U$. \endproof

We define the antipode of $U$ by $S(chb)=(Sb)\cdot(S(ch))$ or $(S(hb))\cdot Sc$
 in terms of the known antipodes $S$ of $H\rbiprod B$ or $\bar C\lbiprod \bar
H$. We do not need its exact form; it can be computed from the formulae in the
Preliminaries and the product (\ref{Uprod}).
\note{\eqn{Uantip}{ S(c\tens h\tens b)=\sum....}}

\begin{lemma} The antipodes of $\bar C\lbiprod \bar H$ and $H\rbiprod B$ extend
to an antipode $S:U\to U$.  \end{lemma}
\proof  The two extensions $S(chb)=(Sb)\cdot(S(ch))$ and $S(chb)=((Shb))\cdot
Sc$ are equal because both are equal (using associativity in $U$, as proven
above) to $(Sb)(Sh)(Sc)$. We use that the restriction $S(ch)$ is the antipode
of $\bar C\lbiprod\bar H$, which is known to be a Hopf algebra\cite{Ma:bos} so
that $S$ is antimultiplicative on $ch$. Likewise, the restriction $S(hb)$ is
the antipode of $H\rbiprod B$ and is therefore antimultiplicative as well.
Working in $U$, we now write $\Delta c=c\o\tens c\t$ (the bosonised coproduct)
where $c\o\in \bar C\lbiprod \bar H$ and $c\t\in \bar C$. Likewise, $\Delta b=
b\o\tens b\t$ where $b\o\in B$ and $b\t\in H\rbiprod B$. Since $U$ is a
bialgebra, as proven above, we have $(S(chb)\o)\cdot(chb)\t=(S((c\o h\o)\cdot
b\o))\cdot(c\t h\t)\cdot b\t=(Sb\o)\cdot((S(c\o h\o))c\t h\t)\cdot b\t=(\eps
c)(Sb\o)b\t=(\eps c)(\eps b)$. We used that the restriction of $S:U\to U$  to
$\bar C\lbiprod\bar H$ is its antipode, and then that the restriction of $S$ to
$H\rbiprod B$ is the antipode for that. Similarly for the proof of antipode
axiom on the other side. \endproof

This completes the proof of Theorem~3.2. When $B$ is finite-dimensional we take
$C=B^\star$ the categorical dual (i.e. $D=B^*$ the ordinary dual) and write
$U=U(B)$, the {\em double-bosonisation} of $B$. In this case we have a
canonical element or coevaluation for the duality pairing. As we have seen in
(\ref{BDexp}), when the braided coproducts are linear ones this coevaluation
plays the role of exponential. See \cite{KemMa:alg} where this point of view is
developed further in general (diagrammatic) terms as part of a braided Fourier
theory. With this in mind, we write
\eqn{bexp}{ \exp_B=\sum  e_a\tens f^a,\quad \bar{\exp}_B=\sum f^a\tens \und S
e_a}
where $\{e_a\}$ is a basis of $B$ with dual basis $\{f^a\}$, and $\und S$ is
the braided antipode of $B$. Here $\bar{\exp}_B=\exp_{B21}^{-1}$ is the
transposed inverse in the usual (unbraided) tensor product algebra $B^*\tens
B$.  A specific example of the same formalism is the braided Fourier
transform\cite{LyuMa:bra} which plays a role in conformal field theory.

\begin{propos} If $B$ is finite-dimensional then its double-bosonisation $U(B)$
is quasitriangular, with quasitriangular structure
\[ \CR_U=\bar{\exp}_B\cdot \CR=\sum (f^a\tens \Ro{2}\Ro{1}\tens 1)\tens
(\Rt{1}\tens \und Se_a\ra \Rt{2})\]
where the product $\cdot$ is in $U\tens U$.
\end{propos}
\proof We verify the quasitriangular structure directly. In the appendix we
introduce a (non-trivial) projection from a suitable Drinfeld quantum double,
which can also be used to obtain $\CR_U$.

{}From a categorical point of view, the exponential is the coevaluation for the
pairing between $B^*$ and $B$. In terms of the structure of $\bar
C=(B^*)^{\und{\rm cop}}$ and $B$, the pairing  (\ref{pairab}) corresponds to
the coevaluation property
\eqn{Dexp}{ (\bar\Delta\tens\id)\bar{\exp}_B=f^a\tens f^b\tens (\und S
e_a)(\und S e_b)=\bar{\exp}_B{}_{13}\bar{\exp}_B{}_{23}.}
Here the numerical suffices denote positions in the tensor product $B^*\tens
B\tens B$. Likewise, (\ref{paircd}) corresponds to the coevaluation property
\eqn{expD}{ (\id\tens\und\Delta)\bar{\exp}_B=(\Rt{}\la f^b)(\Ro{}\la f^a)\tens
\und S e_a\tens \und S e_b.}
This is equivalent to $(\id\tens\und\Delta)\exp_B=\exp_B{}_{12}\exp_B{}_{13}$
given the braided-antimultiplicativity of the braided antipode $\und S$. The
covariance of the pairing corresponds in terms of the coevaluation to $h\la
f^a\tens e_a=f^a\tens e_a\ra h$ for all $h\in H$.

Using (\ref{Dexp}), and since the coproduct of $U$ restricts on $H$ to its
coproduct and $\CR\in H\tens H$ already obeys (\ref{qua}), have
\align{&&\equad (\Delta\tens\id)\CR_U=(\Delta f^a\tens \und S
e_a)\CR_{13}\CR_{23}\\
&&=(f^a\baro\Rmo{}\tens\Rmt{}\la f^a\bart\tens \und S e_a)\CR_{13}\CR_{23}\\
&&=\bar{\exp}_B{}_{13}(\Rmo{}\tens  f^b\tens \und S
e_b\ra\Rmt{})\CR_{13}\CR_{23}\\
&&=\bar{\exp}_B{}_{13}\CR_{13}(\Rmo{1}\Rmo{2}\tens  f^b\tens \Rmt{1}(\und S
e_b\ra\Rmt{2}))\CR_{13}\CR_{23}\\
&&=\bar{\exp}_B{}_{13}\CR_{13}\bar{\exp}_B{}_{23}\CR_{23}
=\CR_U{}_{13}\CR_U{}_{23}.}
Products here are in $U$ or its tensor powers. The second equality applies the
coproduct of $U$ from (\ref{Ucoprod}) or $\bar C\lbiprod\bar H$. The third
equality is (\ref{Dexp}) and covariance of the coevaluation. The fourth inserts
$\CR_{13}\CR^{-1}_{13}$ and allows us (via (\ref{qua})) to recognise the
product in the middle section as $\bar{\exp}_B{}_{23}\CR^{-1}_{13}$ according
to the relations of $H\rbiprod B\subseteq U$.

Similarly, using (\ref{expD}), we have
\align{&&\equad (\id\tens\Delta)\CR_U=(f^a\tens \Delta \und S
e_a)\CR_{13}\CR_{12}\\
&&=(f^a \tens (\und S e_a)\Bo\ra\Ro{}\tens \Rt{}(\und S
e_a)\Bt)\CR_{13}\CR_{12}\\
&&=(f^b (\Ro{2}\Ro{1}\la f^a)\tens \und S e_a\tens \Rt{1}(\und S e_b\ra
\Rt{2})) \CR_{13}\CR_{12}\\
&&=\bar{\exp}_B{}_{13}((\Ro{1}\la f^a)\Ro{2}\tens \und S e_a\tens
\Rt{1}\Rt{2})\CR_{12}\\
&&=\bar{\exp}_B{}_{13}\CR_{13}\bar{\exp}_B{}_{12}\CR_{12}
=\CR_U{}_{13}\CR_U{}_{12}.}
The second equality is the coproduct of $U$ or $H\rbiprod B$, the third is
(\ref{expD}) and covariance of the coevaluation. The fourth recognises the
relations of $H\rbiprod B\subseteq U$ and the fifth recognises the relations of
$\bar C\lbiprod \bar H\subseteq U$.

This proves the first two parts of (\ref{qua}) for $\CR_U$.  Next, we compute
for all $b\in B\subseteq U$,
\align{&&\equad (\Delta^{\rm op} b)\CR_U=\Rt{2}b\Bt f^a\Ro{1}\tens
(b\Bo\ra\Ro{2})(\und S e_a)\Rt{1}\\
&&=\Rt{2}(\Rt{3}\la f^a\bart)\Rt{4}\Rmo{1}(b\Bt\Bt\ra\Rmo{2})\Ro{1}\tens
(b\Bo\ra\Ro{2})(\und S e_a)\Rt{1}\\
&&\ \<\Ro{3}\la f^a\baro, b\Bt\Bo\ra\Ro{4}\>\<\Rmt{1}\la\bar S f^a\barth,
b\Bt\Bth\ra\Rmt{2}\>\\
&&=\Rt{2}(\Rt{3}\la f^b)\Rt{4}\Rmo{1}(b\Bth\ra\Rmo{2})\Ro{1}\tens
(b\Bo\ra\Ro{2})(\und S e_a) (\und S e_b)(\und S e_c)\Rt{1}\\
&&\ \<f^a,b\Bt\ra \Ro{4}\Ro{3}\>\<\bar S f^c,b\Bfo\ra\Rmt{2}\Rmt{1}\>\\
&&= (\Rt{2}\Rt{3}\la f^b)\Rt{5}\Rt{4}\Rmo{1}(b\Bth\ra\Rmo{2})\Ro{1}\\
&&\ \tens (b\Bo\ra\Ro{5}\Ro{2})(\und S b\Bt\ra\Ro{4}\Ro{3})(\und S
e_b)(b\Bfo\ra \Rmt{2}\Rmt{1})\Rt{1})\\
&&= (\Rt{3}\la f^b)\Rt{4}\Rmo{1}(b\Bth\ra\Rmo{2})\Ro{1} \tens ((b\Bo(\und S
b\Bt)\ra\Ro{4}\Ro{3})(\und S e_b)(b\Bfo\ra \Rmt{2}\Rmt{1})\Rt{1})\\
&&= f^b \Rmo{1}(b\Bo\ra\Rmo{2})\Ro{1} \tens (\und S e_b)(b\Bt\ra
\Rmt{2}\Rmt{1})\Rt{1})\\
&&=f^b b\Bo \Rmo{1}\Ro{2}\Ro{1} \tens (\und S e_b)\Rt{1}(b\Bt\ra
\Rmt{1}\Rt{2})=f^b b\Bo\Ro{}\tens (\und S e_b)\Rt{} b\Bt\\
&&=f^b\Ro{1}(b\Bo\ra\Ro{2})\tens (\und S e_b)\Rt{1}\Rt{2} b\Bt=\CR_U\Delta b}
Here the first equality is the definition of the coproduct of $U$ or $H\rbiprod
B$ (transposed) and $\CR_U$. Products are in $U\tens U$. The second equality
uses the relations in Theorem~3.2 to reorder $b\Bt f^a$. The third equality
uses (\ref{Dexp}) in iterated form. The fourth equality evaluates the canonical
element as an identity mapping, and also uses the relations of $\bar C\lbiprod
\bar H\subseteq U$ to move $\Rt{2}$ to the right. The fifth equality uses
(\ref{qua}) in reverse and covariance of the product of $B$. The sixth equality
is the axioms for the braided-antipode $\und S$ of $B$. The seventh equality
recognises the relations of $H\rbiprod B$ as $b\Bth \Rmo{2}$, and also uses
these relations to move $\Rt{1}$ to the left. The result is $\CR_U\Delta b$ by
a further application of these relations of $H\rbiprod H$.

The proof of $(\Delta^{\rm op}c)\CR_U=\CR_U\Delta c$ for all $c\in\bar
C\subseteq U$ is strictly analogous. Finally,
\align{&&\equad (\Delta^{\rm op}h)\CR_U=(h\t f^a\tens h\o\und S e_a)\CR
=(h\t\la f^a)h\th\tens h\o\und S e_a)\CR\\
&&=(f^a h\th\tens h\o (\und S e_a\ra h\t))\CR =(f^a h\t\tens (\und
Se_a)h\o)\CR\\
&&=(f^a\tens Se_a)\CR(h\o\tens h\t)=\CR_U\Delta h}
for all $h\in H\subseteq U$, using in order: the relations in $\CH\subseteq U$,
the covariance of $\bar{\exp}_B$, the relation of $\HB\subseteq U$, and the
quasicocommutativity axiom (\ref{qua}) for $H$. Finally, the element $\CR_U$ is
manifestly invertible, with inverse $\CR^{-1}\exp_B$. \endproof

We see that the quasitriangular structure of $U$ is a product of the
quasitriangular structure of $H$ and the inverse of $\exp_B$. This is the
reason that `q-exponentials' of the root vectors appear in the quasitriangular
structure of $U_q(\cg)$; it is a rather general feature.

\begin{rem}{\rm Once we have constructed our algebra $U$, it is possible to use
the relations of $\CH$ and $\HB$ to write the cross relations in Theorem~3.2
more compactly as
\eqn{Urela}{bc=\Rt{} c\bart b\Bt \Rmo{} \<c\baro,b\Bo\ra\Ro{}\>\<\Rmt{}\la \bar
S c\barth, b\Bth\>}
for all $b\in B\subseteq U$ and $c\in \bar C\subseteq U$. Using the pairing
relation (\ref{paircd}) and axioms (\ref{qua}) we can also write these
relations as
\eqn{Urelb}{ b\Bo c\baro\Rmo{}\<\Rmt{}\la
c\bart,b\Bt\>=\<c\baro,b\Bo\ra\Ro{}\>\Rt{} c\bart b\Bt.}
Both these forms are useful. Note however, that we have {\em not} defined $U$
above as generated by $\bar C,H,B$
and relations between them, but rather we have built $U$ explicitly on the
tensor product space $\bar C\tens H\tens B$. For this purpose, the form shown
in Theorem~3.2 is more suited since it allows us to reorder $bc$ into the
canonical order in $\bar C\tens H\tens B$.  Alternatively, one could take
(\ref{Urela}) etc. as defining cross relations of a Hopf algebra generated by
$U$, but would then have to prove the `triangular decomposition' that the
product map $\bar C\tens H\tens B$ is a linear isomorphism. This triangular
decomposition is an intrinsic feature of our more explicit proofs above.}
\end{rem}

\begin{rem}{\rm We have written the formulae above in terms of $\bar{C}$ rather
a braided group $C$ paired in the more categorical way to $B$. When $H$ is a
(quasitriangular) Hopf algebra, the two are entirely equivalent. Their
underlying vector spaces and pairing maps with $B$ can be identified, their
products are opposite (in the usual sense) and their coproducts and $H$-module
structures are related by
\eqn{barca}{ c\baro\tens c\bart=c\Bo\ra\Rmo{} \tens c\Bt\ra\Rmt, \quad h\la
c=c\ra Sh,\quad \forall h\in H.}
This was explained at the start of the section. The explicit form uses the
braiding in ${}_H\CM$ from (\ref{lrbraid}) and the invariance of $\CR$ under
$S\tens S$. When viewed inside $U\tens U$, we have the further identity
\eqn{barcb}{c\baro\Rmo{}\tens \Rmt{}\la c\bart=\Ro{} c\Bo\tens c\Bt\ra\Rt{},}
which follows from the relations of $\CH$ and (\ref{qua}). Using such formulae,
 it is easy enough to rewrite all of the above in terms of $C,H,B$ rather than
in terms of $\bar C, H,B$. Explicitly, the structure of $U$ in this form is
\ceqn{Urelc}{ \Delta b=b\Bo\ra \Ro{}\tens \Rt{} b\Bt,\quad \Delta c=\Ro{}c\Bo
\tens c\Bt\ra\Rt{},\\
bh=h\o (b\ra h\t),\quad ch=h\t (c\ra h\o)\\
 b\Bo \Ro{} c\Bo \ev(c\Bt\ra \Rt{}, b\Bt)=\ev(c\Bo,b\Bo\ra\Ro{}) c\Bt\Rt{}
b\Bt.}}

{\rm The {\em right-covariant} formulation with $C$ has some advantages  (and
we will sometimes prefer it), because both $B,C$ are then in the same braided
category of right $H$-modules. The above {\em bicovariant} formulation in terms
of $\bar C$ (which is left-covariant) and $B$ (which is right-covariant) is
more natural from a purely algebraic point of view, allowing us the write $U$
in a symmetrical way on $\bar C\tens H\tens B$.}

{\rm In particular, we note that none of the above proofs of the
quasitriangular bialgebra structure of $U$ require an antipode for $H$; we
conclude that if $H$ is a quasitriangular bialgebra, $B$ a braided group in
$\CM_H$, $D\in {}_H\CM$ is dually paired with $B$ in the usual sense and $\bar
C=D^{\und{\rm cop}}$, then $U$ defined in the same way as in Theorem~3.2 is a
bialgebra. It is quasitriangular in the finite-type non-degenerately paired
case.}
\end{rem}

\begin{rem} {\rm As another immediate generalisation, we note that none of the
proofs of the Hopf algebra structure of $U$ actually used $\CR\in H\tens H$
itself but rather, in all expressions we find one component of $\CR$ or
$\CR^{-1}$ acting on $\bar C$ or $B$. Hence we are free to replace this
combination by a coaction of a Hopf algebra $A$ dual $H$ and  composition with
one of the maps $\CR,\bar{\CR}:A\to H$ or their convolution inverses. Hence we
conclude that if $(H,A)$ is weakly quasitriangular, $B\in {}^A\CM$ and $C\in
{}^A\CM$ are categorically dually-paired (or $D\in\CM^A$ is
ordinary-dually-paired) then there is a bialgebra or Hopf algebra $U(\bar
C,H,B)$ defined as above with the above changes. The relations of between these
sub-Hopf algebras becomes
\eqn{Uweakprod}{b\Bo c\baro
\bar{\CR}(c\baro\bt)\<c\bart\bo,b\Bt\>=\<c\baro,b\Bo\bt\>\CR(b\Bo\bo)c\bart
b\Bt}
and the coproduct is
\eqn{Uweakcoprod}{ \Delta b=b\Bo\bt\tens \CR(b\Bo\bo) b\Bt,\quad \Delta
c=c\baro\bar{\CR}(c\bart\bt)\tens c\bart\bo.}
Note that if $b,c$ are primitive elements (for their braided coproducts) then
these relations simplify to
\ceqn{Uweakprim}{[b,c]=\CR(b\bo)\<c,b\bt\>-\<c\bo,b\>\bar{\CR}(c\bt)\\
\Delta b=b\bt\tens \CR(b\bo)+1\tens b,\quad \Delta c=c\tens
1+\bar{\CR}(c\bt)\tens c\bo.}
}
\end{rem}

\begin{rem}{\rm Finally, we note that we also do not use directly the
braided-antipodes $B$ or $C$ or their inverses  in the proof of the bialgebra
structure of $U$. What was actually used in the proofs was the composition
$\<\bar S(\ ),\ \>$ and its properties expressed in terms of the pairing in
(\ref{pairab}) and (\ref{paircd}). Just as for the ordinary quantum
double\cite{Ma:mor}, one can extend the definitions to bialgebras by assuming
in place of $\<\bar S(\ ),\ \>$ a map
$\<\ ,\ \>^{-1}$ characterised as the inverse in the convolution algebra
$\hom(D\tens B,k)$ (the usual tensor product coalgebra here on $D\tens B$). In
terms of $\bar C\tens B$, this is equivalent to (\ref{pairab}) and
(\ref{paircd}). Likewise for the proof of the quasitriangular structure, we
need only assume the inverse of the coevaluation element $\exp_B$ rather than
the antipode of $B$.}
\end{rem}

These are all variations of the theory above; in our presentation we have
chosen the most easily accessible framework (based on modules, assuming the
antipodes etc.) for simplicity of presentation only, leaving the other cases as
routine variations along established lines. We proceed now on the same basis
for further general theory.
To study the representation theory of $U$, however, we do need to break the
left-right symmetry by working  with either left $U$-(co)modules or right
$U$-(co)modules. In the latter case (say) we work with $B,C\in \CM_H$ as
explained in Remark~3.8.

\begin{lemma} In the setting of Theorem~3.2, right $U$-module (algebras) can be
identified with right $H$-module algebras $V$ such that (i) $V$ is braided
$B$-module algebra in the braided category $\CM_H$ (ii) a left braided
$C$-module algebra in the category $\CM_H$ (iii) the actions $\ra,\la$ of $B,C$
are compatible in the sense
\[ c\Bo\la (v\ra b\Bo\Ro{})\ev(c\Bt\ra\Rt,b\Bt)=\ev(c\Bo,b\Bo\ra \Ro{})(c\Bt\la
v)\ra \Rt{}b\Bt\]
for all $v\in V$, $b\in B\subseteq U$, $c\in C\subseteq U$.
\end{lemma}
\proof Since $H\subseteq U$ is a sub-Hopf algebra, we require $V$ an $H$-module
algebra, i.e. an algebra in $\CM_H$. By the bosonisation
theorem\cite[Thm~4.2]{Ma:bos} the modules of $H\rbiprod B$ and their tensor
products coincide with the (right) $B$-modules in $\CM_H$ and their (braided)
tensor products computed categorically in $\CM_H$. This is part (i). For part
(ii) the same picture applies for left $\bar C$-modules in ${}_{\bar H}\CM$,
but this is not directly applicable now. Instead, we need to realise the
relations of $\CH$ acting from the right. We view an action of $\bar C$ from
the right in the trivial way as an action of $C$ from the left. For this to be
a morphism in $\CM_H$, we require $(c\la v)\ra h=(c\ra h\o)\la (v\ra h\t)$. In
terms of the corresponding action of $U$, the left hand side is $v\ra ch$ and
the right hand side is $v\ra(h\t (c\ra h\o))$, which is the cross relation in
(\ref{Urelc}) for $\CH$ when expressed in terms of $C,H$. The tensor product
actions also match; thus $c\la(v\tens w)=c\Bo\la (v\ra\Ro{})\tens
(c\Bt\ra\Rt{})\la w=v\ra(\Rmo{}c\Bo)\tens w\ra(c\Bt\ra\Rmt{})=(v\tens
w)\ra\Delta c$ which is the tensor product action according to the coproduct of
$U$.
The first equality here is the definition of the braided tensor product
$C$-module\cite{Ma:tra} computed in the category $\CM_H$ using the braiding
from (\ref{lrbraid}). The last equality is (\ref{barcb}). Finally, part (iii)
is manifestly the requirement that these module-algebra structures in $V$
respect the final cross relation in (\ref{Urelc}), which is equivalent to the
reordering relation in Theorem~3.2. \endproof

This digression into braided category theory is developed further in
Appendix~C. Here we give a purely algebraic consequence of it.

\begin{theorem} In the setting of Theorem~3.2, the algebra $V=B$ is a right
$U$-module algebra by the maps
\[ v\ra b=(\und S b\Bo\ra\Ro{})(v\ra\Rt{})b\Bt,\quad v\ra c=\<\bar S
c,v\Bt\ra\Rmt{}\>v\Bo\ra\Rmo{} \]
and the tautological action of $H$ (i.e. the same as for $B$ as an object in
$\CM_H$.) We call this the {\em fundamental} or {\em Schr\"odinger}
representation of $U$ on $B$.
\end{theorem}
\proof This is best done diagrammatically; see the appendix Proposition~C.4.
The action of $B$ is the right braided adjoint action of $B$ on itself. The
right action of $\bar C$ is the left action of $C$ from the proposition, which
is the standard braided coregular representation. Both are known braided group
constructions \cite{Ma:introm}\cite{Ma:lie}\cite{Ma:fre} and Proposition~C.4
checks that they are suitably compatible to form a representation of $U$. We
then convert over from $C$ to $\bar C$ as explained in Remark~3.8. The
restriction to $\CH$ is a version of the fundamental representation used
already in \cite{Ma:csta}. One can also verify directly by the same techniques
as in the proofs above that these actions on $V$ are compatible as required in
part (iii) of the preceding lemma. \endproof

Our previous Remarks~3.8 and~3.9 about only requiring $H$ to be a bialgebra and
only part of a weakly quasitriangular pair apply. The fundamental
representation in Theorem~3.12 becomes in the weak case
\eqn{weakfund}{ v\ra b=(\und S b\Bo\bt) (v\ra \CR(b\Bo\bo))b\Bt,\quad v\ra
c=\<\bar S c,v\Bt\bo\>v\Bo\ra\bar{\CR}(v\Bt\bt).}
If $b$ is braided-primitive, we have
\eqn{fundprim}{v\ra b=vb-b\bt (v\ra{\CR}(b\bo)),}
which is the form that we will use. The action $\ra c$ is by braided
differentiation\cite{Ma:fre} on the braided group $B^{\und{\rm cop}}$ with the
braided opposite coproduct.

To conclude our general theory, we note that $H\subseteq U$ as Hopf algebras.
Hence by (a right handed version of) the  transmutation theorem\cite{Ma:tra}
there is a braided group   $\und U=B(H,U)\in\CM_H$ which consists of $U$ as an
algebra but has
a modified coproduct. The action of $H$ on $\und U$ is by the right adjoint
action $h\ra g=(Sg\o)h g\t$, the given action on $B$ and the action $c\ra
g=(Sg\o)cg\t=(Sg\t\la c)(Sg\o)g\th$ on $\bar C$. Here $\und U\supseteq
B(H,H)\rcross B$ a (right handed) braided group cross product by the braided
group $B(H,H)$ associated to the identity mapping. Indeed, this is the original
construction of the bosonisation $\HB$ in \cite{Ma:bos} as such that its
transmutation is a cross product. Moreover, in the finite-dimensional
non-degenerately paired case we
know from \cite{Ma:tra} that $\und U$ is braided-quasitriangular (a quantum
braided group in the strict sense).
Explicitly, the braided coproducts on $H,B\subseteq \und U$ and its
braided-quasitriangular structure are
\eqn{undU}{ \Delta_{\und U} h=h\o\ra\Ro{}\tens (S\Rt{})h\t,\quad \Delta_{\und
U} b=b\Bo\ra\Ro{}\tens (S\Rt{})b\Bt,\quad \CR_{\und U}=\sum f^a\ra\Ro{}\tens
(S\Rt{})\und Se_a,}
while $\Delta_{\und U}$ restricted to $\bar C$ is more complicated. We can also
transmute by $\bar H\subseteq U$, in which case the restriction to $\bar C$ is
simple and the restriction to $B$ more complicated.

\section{Recovering Lusztig's construction of {$U_q(\cg)$}}

A {\em Cartan datum} in Lusztig's construction of $U_q(\cg)$ is a set $I=\{i\}$
and a symmetric bilinear form $\cdot$ on $\Z[I]$ such that
\eqn{cartan}{ i\cdot i\in\{2,4,6,\cdots\},\quad a_{ij}\equiv {2i\cdot j\over
i\cdot i}\in \{0,-1,-2,\cdots\},\quad \forall i\ne j.}
Thus, $a_{ij}$ is a symmetrizable Cartan matrix. Also defined is  a {\em root
datum}, which is two finitely generated free Abelian groups $Y,X$ with a
perfect pairing $\<\ ,\ \>:Y\times X\to \Z$ and inclusions (the second denoted
$i\mapsto i'$)
\eqn{root}{Y\supseteq I\subseteq X,\quad {\rm s.t.}\quad \<i,j'\>=a_{ij}.}

We show now that these data provide now a weakly quasitriangular dual pair in
the sense explained in the Preliminaries. We let $H=k Y$ with basis elements
$\{K_\mu:\mu\in Y\}$. This forms a Hopf algebra with $\Delta K_\mu=K_\mu\tens
K_\mu$ and product $K_\mu K_\nu=K_{\mu+\nu}$, extended linearly. Let
$A=k\Z[I]=k[g_i,g_i^{-1}]$ the group algebra of $\Z[I]$ with $\Delta
g_i=g_i\tens g_i$. Using the Cartan datum (symmetric or not) it is clear that
$A$ is dual quasitriangular, with $\CR(g_i,g_j)=q^{i\cdot j}$, for any $q\in
k^*$. This more general point of view is relevant to Appendix~B. For the
present we really need a weakly quasitriangular dual pair. We can work over
$k=\Q(q)$, for example.

\begin{lemma} Using the root datum, we define a pairing
\eqn{Luspair}{\<\ ,\ \>:H\tens A\to k,\quad \<K_\mu,g_i\>=q^{\<\mu,i'\>}}
and relative to it, we have a weak quasitriangular structure
\eqn{Lusweak}{ \CR,\bar{\CR}:A\to H,\quad \CR(g_i)=K_i^{i\cdot i\over 2},\quad
\bar{\CR}(g_i)=K_i^{-{i\cdot i\over 2}}.}
\end{lemma}
\proof The inclusion $I\subseteq X$ induces a homomorphism $\Z[I]\to X$ and
hence a homomorphism $A\to kX$ of Hopf algebras. The latter is dually paired
with $kY$ via the assumed group pairing. This gives the pairing between $H,A$
(it need no longer be non-degenerate, however). We also have well-define
algebra homomorphisms $\CR,\bar{\CR}$ as stated. Here $\CR^{-1}=\bar{\CR}$ as a
consequence of the symmetry of $\cdot$. We check that
\[ \<{\CR}^{-1}(g_j),g_i\>=q^{-{j\cdot j\over 2}\<j,i'\>}=q^{-j\cdot
i}=q^{-i\cdot j}=q^{-{i\cdot i\over 2}\<i,j'\>}=\<K_i^{-{i\cdot i\over
2}},g_j\>=\<\bar{\CR}(g_i),g_j\>\]
in virtue of the symmetry. \endproof

Next, we let $\tilde{B}=k\<e^i\>$ the free non-commutative algebra on $I$. This
lives
in the braided category of left $A$-comodules by $e^i\mapsto g_i\tens e^i$.
Hence it also lives in the category of right $H$-modules by
\eqn{eK}{ e^i\ra K_\mu=\<K_\mu,g_i\>e^i=q^{\<\mu,i'\>}e^i.}

\begin{lemma} The maps defined by $\und\Delta e^i=e^i\tens 1+1\tens e^i$,
$\und\eps e^i=0$, $\und Se^i=-e^i$ make $\tilde{B}$ a braided group
in the category of left $A$-comodules with braiding provided by the weak
quasitriangular structure from Lemma~4.1.
\end{lemma}
\proof  We check that the braiding as defined by the weak quasitriangular
structure is the desired one, namely
\[ \Psi(e^i\tens e^j)=e^j\tens e^i\ra \CR(g_j)=e^j\tens e^i\ra K_j^{j\cdot
j\over 2}=e^j\tens e^i q^{\<j,i'\> {j\cdot j\over 2}}=q^{j\cdot i}e^j\tens
e^i.\]
The rest is clear from \cite{Ma:any}\cite{Ma:csta}\cite{Ma:fre} or the
computations in Lusztig \cite{Lus}. Indeed, the theory of braided groups
ensures that we have natural braided tensor product algebras
$\tilde{B}\und\tens \tilde{B}$ (and higher braided tensor products as well).
The relations are
\[ (1\tens e^i)(e^j\tens 1)=q^{j\cdot i}(e^j\tens 1)(1\tens e^i).\]
We extend $\und\Delta:\tilde{B}\to \tilde{B}\und\tens\tilde{B}$ as an algebra
homomorphism, and $\und S:\tilde{B}\to \tilde{B}$ as a braided
anti-homomorphism using $\Psi$. \endproof

It is also clear from \cite{Ma:fre} or computations in \cite{Lus} that
$\tilde{B}$ has
ordinary-dual $\tilde{D}=k\<f_i\>$ in the category of right $A$-comodules by
$f_i\mapsto f_i\tens g_i$. It also lives in the category of left $H$-modules by
$K_\mu\la f_i=\<K_\mu,g_i\>f_i=q^{\<\mu,i'\>}f_i$ and forms a braided group
with $f_i$ braided-primitive and $\Psi(f_i\tens f_j)=q^{i\cdot j}f_j\tens f_i$.
We take the pairing with
$\tilde{B}$ to be
\eqn{tBDpair}{ \<f_i,e^j\>=(q_i-q_i^{-1})^{-1}\delta_i{}^j,}
where $q_i=q^{i\cdot i\over 2}$; we inserted here a choice of normalisation
factor for each $e^i$. Following Lusztig, we pass to the quotients $B,D$ by the
radical of the pairing, generated by the $q$-Serre-relations. We let $\bar
C=D^{\und{\rm cop}}$.

\begin{propos} The $U(\bar C,H,B)$ construction in Section~3 reduces in this
setting to Lusztig's construction of $U_q(\cg)$ in suitable (right-handed)
conventions.
\end{propos}
\proof The braided group $\bar C$  has the primitive braided coproduct
$\bar\Delta f_i=f_i\tens 1+1\tens f_i$ (because $\Psi(1\tens f_i)=f_i\tens 1$
etc.) But it extends to products with the inverse braiding to that of
$D$\cite{Ma:introp}, which in our case means
\[ \Psi(f_i\tens f_i)=\bar{\CR}(g_i)\la f_j\tens f_i=q^{-i\cdot j} f_j\tens
f_i.\]

The algebras $\bar C,H,B$ are all included in $U$. From the relations
$bh=h\o(b\ra h\t)$ and $hc=(h\o\la c)h\t$ in $\HB$ and $\CH$, we have
\[ e^iK_\mu=\<K_\mu,g_i\>K_\mu e^i=q^{\<\mu,i'\>}K_\mu e^i,\quad K_\mu f_i= f_i
K_\mu \<K_\mu,g_i\>=q^{\<\mu,i'\>}f_i K_\mu.\]
{}From the formulae (\ref{Uweakprim}) we have  the cross relations and
coproducts
\cmath{ [e^i,f_j]=(\CR(g_i) -\bar{\CR}(g_j))\<f_j,e^i\>={K_i^{i\cdot i\over
2}-K^{-{i\cdot i\over 2}}_i\over q_i-q_i^{-1}}\delta^i{}_j\\
\Delta e^i=e^i\tens K_i^{i\cdot i\over 2}+1\tens e^i,\quad \Delta f_i=f_i\tens
1+K_i^{-{i\cdot i\over 2}}\tens f_i.}
Hence we recover the  structure of $U_q(\cg)$ in \cite{Lus} in suitable
conventions. The identification is by $e^i=-E_i$, $f_i=F_i$ and an interchange
of $K_i$ with $K_i^{-1}$. \endproof

The relations between non-primitive elements of $\bar C,B$ and their coproducts
  follow just as easily from Theorem~3.2, as
\eqn{genUqg}{bc=\tilde{K}_{|b\Bo|} c\bart b\Bo
\tilde{K}^{-1}_{|c\barth|}\<c\baro,b\Bo\>\<\bar S c\barth, b\Bth\>,\quad
\Delta b=b\Bo\tens \tilde{K}_{|b\Bo|} b\Bt,\quad \Delta c=c\baro
\tilde{K}^{-1}_{|c\bart|}\tens c\bart}
for all $b\in B$, $c\in\bar C$ of homogeneous degree.  Here $|e^i|=|f_i|=i\in
\Z[I]$ corresponding to the coactions above, and $\tilde{K}_{\sum_i \nu_i
i}\equiv \prod_i K_i^{{i\cdot i\over 2}\nu_i}=\CR(\prod_i g_i^{\nu_i})$ for
$\sum_i\nu_i i\in\Z[I]$. We also deduce the triangular decomposition of
$U_q(\cg)$ into $\bar C,H,B$. These facts and formulae all require substantial
proof in \cite[Sec. 3.1.5, Prop. 3.1.7, Sec. 3.2]{Lus}, where $U_q(\cg)$ is
defined by generators and relations.

We are also in a position to apply general constructions for braided groups,
e.g. proven diagrammatically, to obtain results about $U_q(\cg)$ somewhat more
easily than by the usual direct calculation.

\begin{propos} The fundamental representation in Theorem~3.12 of $U$ on $B$,
with generators denoted $x^i$, is
\[ x^i\ra K_\mu=q^{\<\mu,i'\>}x^i,\quad v(x)\ra e^i={v(x)x^i-x^iv(x\ra
K_i^{i\cdot i\over 2})\over q_i-q^{-1}_i},\quad v(x)\ra f_i=-\del_i v(x),\]
where $\del_i$ is the braided differentiation
\[ \del_i\left( (x^1)^{\nu_1}\cdots
(x^r)^{\nu_r}\right)=q^{-i\cdot\sum_{j=1}^{i-1}\nu_j j}(x^1)^{\nu_1}\cdots
(x^{i-1})^{\nu_{i-1}}[\nu_i,q_i^{-2}](x^i)^{\nu_i-1}(x^{i+1})^{\nu_{i+1}}\cdots
(x^r)^{\nu_r}\]
where $I=\{1,\cdots,r\}$ and $[m,q]={1-q^m\over 1-q}$. This representation is
adjoint to a version of Lusztig's Verma module representation in \cite{Lus}.
Moreover, it makes $B$ into a $U$-module algebra.
\end{propos}
\proof The action of $H$ is the given action on $B$ from the right. The action
of $e^i$ is the right braided adjoint action (\ref{fundprim}). More precisely,
we have rescaled the $e^i$ in (\ref{tBDpair}) by a factor
$(q_i-q_i^{-1})^{-1}$, whereas we keep the $x^i$ in the usual normalisation
where $\<f_i,x^j\>=\delta_i{}^j$. For the action of $f_i$, we write the action
in (\ref{fundprim}) as
\[v\ra f_i=\<\bar Sf_i ,v\opBo\>v\opBt=-\del_i v\]
where $v\opBo \tens v\opBt\equiv\Psi^{-1}(v\Bo\tens v\Bt)$ is the coproduct of
$B^{\und{\rm cop}}$. This has the same linear form on the generators $x^i$
since $\Psi(x^i\tens 1)=1\tens x^i$ etc., but extends to products with the
opposite braiding according to \cite{Ma:introp}. Hence we can compute this as
shown, where $\del_i v$ is characterised by $v\opBo \tens v\Bt=x^i\tens \del_i
v+$terms not of the form $x^i\tens(\ )$. It is a version of the operator
${}_ir$ in \cite{Lus}. The braided-module algebra properties of the
braided-coregular representation correspond to the braided-Leibniz
rule\cite{Ma:poi}\cite{Ma:fre}, which in the present case takes the form
$\del_i(vw)=(\del_iv)w+q^{-i\cdot|v|}v\del_iw$ for $v,w\in B$ and $v$
homogeneous, as in \cite{Ma:poi}\cite{Lus}.   From this, we find easily the
explicit formula shown for $\del_i$. It can also be obtained from (\ref{bdiff})
with $R^i{}_j{}^k{}_l=\delta^i{}_j{}^k{}_l q^{-i\cdot k}$.

By contrast, the Verma representation in \cite[Sec. 3.4.5]{Lus} consists of
$f_i$ acting by multiplication on $\bar C$ (both $B$ and $\bar C$ are versions
of Lusztig's algebra $f$). The action of $e^i$ is a difference of the
differentiation operators ${}_ir$ and $r_i$ in \cite{Lus}. This is clearly the
adjoint of the fundamental representation above: the adjoint under the duality
pairing of braided differentiation is braided multiplication. Note, however,
that $U$ does not respect the algebra structure of $\bar C$ (it respects its
coalgebra). For this reason we consider the action on $B$ rather than its
adjoint to be more fundamental, since by Theorem~3.12 we know that $V=B$
becomes a $U$-module algebra. \endproof

There is also a  weak quasitriangular structure which can  be obtained from a
version of Proposition~3.6, at least for the classical ABCD series of
finite-dimensional semisimple Lie algebras. It requires, however, a description
of the coordinate algebra $G_q$ dual to $U_q(\cg)$, which we do not cover here.
We remark only that in a suitable setting where we work over formal power
series\cite{Dri}, we obtain the formula
\[  \CR_{U_q(\cg)} = \bar{\exp}_B\CR_H;\quad \CR_H=q^{\sum_{i,j}h_i\tens h_j
{i\cdot i\over 2} (a^{-1})_{ij}}.\]
Here we assume $K_i=q^{h_i}$ and $\bar{\exp}_B$ is obtained from the duality
pairing coevaluation. Because the duality pairing is non-degenerate we know
that $\bar{\exp}_B$ exists as a formal power series built from a basis and dual
basis of $B$. It corresponds to the element $\Theta$ studied in \cite[Sec.
4]{Lus}.

In the simplest case we  take for $B$ the braided line $k\<e\>$ and for
$\bar{C}$ the braided line $k\<f\>$. Then
$U=U_q(sl_2)$. Its fundamental representation on $k\<x\>$ is
\eqn{1conf}{ v(x)\ra K=v(q^2x),\quad v\ra e=-qx^2\del_{q^2}v,\quad v\ra
f=-\del_{q^{-2}}v,}
where $\del_qv(x)={v(x)-v(qx)\over (1-q)x)}$ is the usual 1-dimensional
$q$-derivative on polynomials $v\in k\<x\>$.
The fundamental representation in this case is a $q$-deformation of the action
on $k\<x\>$ of $sl_2$ as the degree 1,0,-1 subalgebra of the Virasoro algebra
in
physics. Meanwhile, the braided exponential part of the quasitriangular
structure comes out as
\eqn{1exp}{ \bar{\exp}_B=e_{q^{-2}}^{-(q-q^{-1})f\tens e},}
as in \cite{KirRes:rep}, where $e_q$ denotes the usual $q$-exponential defined
with $[m;q]!=[m;q]\cdots [2;q]$ in place of $m!$.

This completes our outline of how the explicit constructions in  \cite[Part
I]{Lus} and other results can be recovered as an application of the double
bosonisation and its properties in Section~3. The theory applies in other
related settings just as well:

\begin{example} If $q\ne 1$ is a primitive $r$-th root of $1$ with $r$ odd and
invertible in $k$, we take for $B$ the anyonic line $k\<e\>/e^r$ from
\cite{Ma:any}. This forms a braided group in the category of
$k\Z/r\Z$-comodules (or modules, since  $\Z/r\Z$ is self-dual). Then $U(B)$ is
the finite-dimensional reduced form $u_q(sl_2)$ with $K^r=1,e^r=0=f^r$. The
fundamental representation on $k\<x\>/x^r$ is as in (\ref{1conf}). The
quasitriangular structure from Proposition~3.6  is the product of
\[ \bar{\exp}_B=\sum_{m=0}^{r-1}f^m \tens (-e)^m (q-q^{-1})^m
([m;q^{-2}]!)^{-1},\quad \CR_H=r^{-1}\sum_{m,n=0}^{r-1} q^{-2mn} K^m\tens
K^n.\]
\end{example}
\proof We use the quasitriangular $\CR_H$ structure on $H=k\Z/r\Z$ introduced
in \cite{Ma:any}.  The braided group $B=k\<e\>/e^r$ is also introduced there,
with $|e|=1$ and  $\und\Delta e=e\tens1+1\tens e$. We take a slightly different
braiding (to fit with the conventions above), namely $\Psi(e^m\tens
e^n)=q^{2mn}e^n\tens e^m$ defined with $q^2$ in place of $q$. Similarly for
$k\<f\>/f^r$. The pairing $\<f^m,e^n\>=\delta^m_n[m;q^2]!$ is non-degenerate
between these finite-dimensional braided groups. Hence the coevaluation is
$\exp_B=\sum_{m=0}^{r-1} e^m\tens f^m (q-q^{-1})^m([m;q^2]!)^{-1}$. The
braided-antimultiplicativity of the braided antipode implies that $\und
S(e^m)=q^{m(m-1)}(-e)^m$, which gives $\bar{\exp}_B$ as stated in the example
(the same applies to the formal power series (\ref{1exp})). We also check that
the quasitriangular structure $\CR_H$  induces the correct weak quasitriangular
structure, by evaluation. Thus $A=k\Z/r\Z$ with pairing $\<K,g\>=q^2$ and
$\<\Ro{},g\>\Rt{}=r^{-1}\sum_{m,l=0}^{r-1}q^{-2mn}q^{2m}K^n=K$ since
$r^{-1}\sum_{m=0}^{r-1} q^{m(2(n-1))}=\delta_{2(n-1),0}=\delta_{n,1}$. The
Kronecker delta functions here are on $\Z/r\Z$. \endproof

This recovers $u_q(sl_2)$ as in \cite{Lus:roo} and (with the quasitriangular
structure) \cite{Ma:any}\cite{LyuMa:bra}. The braided version $\und U$ from
(\ref{undU}) recovers the anyonic quantum group $\und{u_q(sl_2)}$ in
\cite{Ma:any}.  In a different direction, we can include as well the case where
$\cdot$ on $\Z[I]$ is antisymmetric. For example, we can suppose
$\<i,j'\>=i\cdot j$ and define $\CR(g_i)=K_i=\bar{\CR}(g_i)$ as a weak
triangular structure. This case is not very interesting for us, however: The
category of comodules in this case is symmetric rather than braided and the
relations in (\ref{Uweakprim}) reduce to $[b,c]=0$ for all braided-primitive
$b\in B$ and $c\in\bar C$.

\section{New quantum group constructions}

We now apply our constructions to more general quantum groups $H$ in the role
of `Cartan' subalgebra. Here the main datum we need is a $R$-matrix, i.e. a
matrix $R\in M_n\tens M_n$ which is an (invertible) solution of the matrix
Yang-Baxter equations. We let $A(R)$ denote the usual matrix bialgebra with
generators $\vect=\{t^i{}_j\}$, relations $R\vect_1\vect_2=\vect_2\vect_1R$ and
coproduct $\Delta \vect=\vect\tens\vect$ in standard notation\cite{FRT:lie}. We
let $\widetilde{U(R)}=A(R)\dcross A(R)$ the double cross product bialgebra
constructed in \cite{Ma:mor}. It consists of two copies of $A(R)$ with
generators $\vecm^-,\vecm^+$, say, and the cross relations and coalgebra
\eqn{UR}{R\vecm^+_1\vecm^-_2=\vecm^-_2\vecm^+_1R,\quad \Delta
m^\pm{}^i{}_j=\sum_a m^\pm{}^a{}_j\tens m^\pm{}^i{}_a,\quad \eps
m^\pm{}^i{}_j=\delta^i{}_j.}
Combining results in \cite{Ma:mor} for the weak quasitriangular structure and
\cite{Ma:poi} for the braided plane $V(R)$, we have:

\begin{lemma} Let $R$ be an invertible matrix solution of the QYBE. Then
$(\widetilde{U(R)},A(R))$ is a weakly quasitriangular dual pair with
$\CR(\vect)=\vecm^+$ and $\bar{\CR}(\vect)=\vecm^-$ and  $V(R)=k\<e^i\>$ is a
braided group with $\und\Delta e^i=e^i\tens 1+1\tens e^i$ extended by the
corresponding braiding.
\end{lemma}
\proof Motivated by \cite{FRT:lie}, we considered in \cite{Ma:mor} the
bialgebras $A(R)$ and $\widetilde{U(R)}$ as dually paired by
$\<m^+{}^i{}_j,t^k{}_l\>=R^i{}_j{}^k{}_l$ and
$\<m^-{}^i{}_j,t^k{}_l\>=R^{-1}{}^k{}_l{}^i{}_j$, and verified that $\CR$
extends as a weak quasitriangular structure. The same applies for $\bar{\CR}$,
and the two are clearly inverse-transpose as required for a weak
quasitriangular structure in the sense (\ref{wqua}). There is left coaction of
$A(R)$ on $V(R)$ given by $e^i\mapsto t^i{}_a\tens e^a$, which induces a right
action of $\widetilde{U(R)}$ by
\eqn{emact}{ e^i\ra m^+{}^j{}_k=e^a\<m^+{}^j{}_k,t^i{}_a\>=R^j{}_k{}^i{}_a
e^a,\quad e^i\ra m^-{}^j{}_k=e^a\<m^-{}^j{}_k,t^i{}_a\>=R^{-1}{}^i{}_a{}^j{}_k
e^a.}
The induced braiding $\Psi(e^i\tens e^j)=e^a\tens e^i\ra\CR(t^j{}_a)$ is then
the correct one for $V(R)$ as required in (\ref{freevec}) in the Preliminaries.
We use it to define a braided tensor product algebra $V(R)\und\tens V(R)$ and
extend $\und\Delta$ as an algebra homomorphism to it\cite{Ma:poi}. \endproof

We take this for $B$ and we take $D=\Vhaj(R)=k\<f_i\>$ as described in the
Preliminaries. The braided group $\bar C$ has the same algebra as $D$ and the
same linear form of the braided-coproduct on the generators $f_i$, but extended
to products with the inverse braiding, i.e. it is $\Vhaj(R_{21}^{-1})$. The
bialgebra version of Theorem~3.2 yields:

\begin{propos} There is a bialgebra
$U=U(V(R),\widetilde{U(R)},\Vhaj(R_{21}^{-1}))$ generated by $\vecm^\pm,
\vece=\{e^i\},\vecf=\{f_i\}$ with the cross relations and coproduct
\cmath{ e^im^+{}^j{}_k=R^j{}_a{}^i{}_b m^+{}^a{}_k e^b ,\quad  m^-{}^i{}_j
e^k=R^k{}_a{}^i{}_b e^a m^-{}^b{}_j,\\
 m^+{}^i{}_jf_k=f_bm^+{}^i{}_a R^a{}_j{}^b{}_k,\quad  f_i
m^-{}^j{}_k=m^-{}^j{}_b f_a R^a{}_i{}^b{}_k,\quad [e^i,f_j]={m^+{}^i{}_j -
m^-{}^i{}_j\over q-q^{-1}}\\
\Delta e^i=e^a\tens m^+{}^i{}_a+ 1\tens e^a,\quad \Delta f_i= f_i\tens
1+m^-{}^a{}_i\tens f_a,\quad \eps e^i=\eps f_i=0.}
Here $\widetilde{U(R)}$ appears as a sub-bialgebra. The factor $q-q^{-1}\in
k^*$ is an arbitrary choice of normalisation for the $e^i$, chosen for
conventional purposes.
\end{propos}
\proof The action (\ref{emact}) and a similar computation for the action on
$f_i$, immediately give the relations
with $\vecm^\pm$ (using the matrix form of the coproduct of the latter). We
then use the formulae for the pairing and weak quasitriangular structure to
compute the cross relations and coproduct  from (\ref{Uweakprim}), giving the
results stated. In doing so, we introduce an overall normalisation factor
$(q-q^{-1})$ for the $e^i$, so that $\<f_i,e^j\>=(q-q^{-1})^{-1}\delta_i{}^j$
for some $q$ with $q^2\ne 0,1$. This is purely conventional to suit the
examples of interest; it can be any constant. \endproof

At this level, the formulae (\ref{bdiff})--(\ref{BDexp}) from \cite{Ma:fre}
provide us with the structure needed for $\bar{\exp}_B$ and for the fundamental
representation. Thus,

\begin{propos} In the setting of Proposition~5.2, the free algebra $V(R)$,
generated by $x^i$, is a $U$-module algebra by
\[ x^i\ra m^+{}^j{}_k= R^j{}_k{}^i{}_a x^a,\quad x^i\ra m^-{}^j{}_k =
R^{-1}{}^i{}_a{}^j{}_k x^a,\quad v\ra e^i={vx^i- x^a (v\ra m^+{}^i{}_a)\over
q-q^{-1}},\quad v\ra f_i=-\del_i v\]
for all $v\in V(R)$, where $\del_i$ is the braided differentiation\cite{Ma:fre}
on $V(R_{21}^{-1})$ from (\ref{bdiff}).
\end{propos}
\proof We action of $\vecm^\pm$ is the given action on $B$ in (\ref{emact}). We
use the right braided adjoint action (\ref{fundprim}), computed again from
Lemma~5.1, for the action of $e^i$. We take $x^i$ in their natural
normalisation where the pairing with $f_i$ is $\<f_i,x^j\>=\delta_i{}^j$. For
the action of $f_i$ we use the formulation of the coregular representation as
braided-differentiation introduced in \cite{Ma:fre}; as in the proof of
Proposition~4.4, we write it as evaluation against the coproduct of
$V(R)^{\und{\rm cop}}$, which is $V(R)$ with the opposite braiding, i.e.
$V(R_{21}^{-1})$. Hence $\del_i$ is given on monomials by the braided-integer
matrices $[m;R_{21}^{-1}]$ in the notation of the Preliminaries. \endproof

Next we consider the same construction at the level of quotient Hopf algebras
and quotient braided groups. These steps depend in fact on the normalisation of
$R$. In the framework of \cite{Ma:poi} where braided groups with quadratic
relations are constructed from $k\<e^i\>$, a necessary condition is that the
matrix $PR$ has an eigenvalue $-1$. The possible quadratic relations are
determined by a matrix $R'$ as in (\ref{BDrel}), obeying certain conditions. We
fix a choice of $R,R'$ (the braided plane data). Given these, we look for
quotients of the bialgebras $\widetilde{U(R)}$  and $A(R)$ which are Hopf
algebras, such that the pairing and weak quasitriangular structure descend.
Typically, this is possible provided we normalise $R$ (which enters into the
pairing) correctly,
e.g. provided we modify the pairing to $\<\vecm_1^+,\vect_2\>=\lambda R$ and
$\<\vecm_1^-,\vect_2\>=\lambda R_{21}^{-1}$, where $\lambda\in k^*$ is a
suitable constant. We say that such $R$ is {\em regular} and that $\lambda$ is
a {\em quantum group normalisation constant}. This framework has been
introduced (in an equivalent form) in \cite{Ma:poi}.

The $R$-matrices for the standard $ABCD$ series of Lie algebras are
known\cite{FRT:lie} and are regular in this sense when we work over $\C$, with
quotient weakly quasitriangular dual pair $(U_q(\cg),G_q)$ in suitable form
(which may be slightly different, however, from Lusztig's `minimal' form in
Section~4). Then we identify the image $\vecm^\pm =S\vecl^\mp$ in the
conventional notation of \cite{FRT:lie}, where $S$ is the antipode of
$U_q(\cg)$. Other examples of interest include the $q$-Lorentz group dual pair
\cite{Mey:new}. The quantum-braided planes $V(R)$ and $\Vhaj(R)$ and their
suitable quotients remain covariant under the quotients $G_q$ etc. of $A(R)$.
Lemma~5.1 no longer goes through, however.

\begin{lemma}cf.\cite{Ma:poi} Let $R$ be regular and $(H,A)$ a quantum group
quotient of $(\widetilde{U(R)},A(R))$
with associated normalisation constant $\lambda$. Let $\widetilde{H}=H\tens
k[c]$ and $\widetilde{A}=A\tens k[g]$ be
the centrally extended weakly quasitriangular dual pair defined by
\[ \Delta c=c\tens c,\quad\Delta g=g\tens g,\quad \<c,g\>=\lambda,\quad
\CR(g)=c^{-1},\quad \bar{\CR}(g)=c.\]
Then $V(R)$ and $\Vhaj(R)$ (and their covariant quotients) are braided groups
in the category of $\widetilde{A}$-comodules by the coactions $e^i\mapsto g
t^i{}_a\tens e^a$ and $f_i\mapsto f_a\tens g t^a{}_i$.
\end{lemma}
\proof We cast the construction in \cite{Ma:poi} into the present weakly
quasitriangular setting. Indeed, the induced action of $c$  is $e^i\ra
c=\lambda e^i$, hence $\Psi(e^i\tens e^j)=e^a\tens e^i\ra\CR(gt^j{}_a)=e^a\tens
e^i\ra c^{-1}m^+{}^j{}_a=R^j{}_a{}^i{}_b e^a\tens e^b$ as required. Similarly
$f_i\ra c=\lambda f_i$ adjusts correctly for the braiding. \endproof

We can then make the bosonisations $\widetilde{H}\rbiprod V(R',R)$ etc., which
is the general construction of inhomogeneous quantum introduced in
\cite{Ma:poi}. The elements $c,g$ are called in this context `dilaton
generators'.  The construction in \cite{Ma:poi} recovered  some of the specific
examples of inhomogeneous quantum groups obtained by other means. We make the
same extension when constructing our Hopf algebra $U$, whenever the
appropriate quantum group normalisation constant is not $1$. To be concrete, we
specialise to one of the standard weakly quasitriangular dual pairs
$(U_q(\cg),G_q)$ and $\A_q=V(R',R)$ a choice of quantum-braided plane covariant
as algebras under $G_q$ (i.e. of $G_q$ type) and
$\bar{\A^\star}_q=\Vhaj(R',R_{21}^{-1})$ its dual with braided-opposite
coproduct. The same formulae hold at the level of generality of Lemma~5.4.

\begin{corol} Let $\A_q$ be a quantum-braided plane of $G_q$ type and $\lambda$
the associated quantum group normalisation constant. Then
$U=U(\A_q,\widetilde{U_q(g)},\bar{\A^\star}_q)$ has the $m^\pm$ relations as in
Proposition~5.2 with $\lambda R$ in place of $R$, and the cross relations and
coalgebra
\cmath{  cf_i=\lambda f_i c, \quad e^ic=\lambda c e^i,\quad
[c,\vecm^\pm]=0,\quad [e^i,f_j]={m^+{}^i{}_j c^{-1}- c m^-{}^i{}_j\over
q-q^{-1}}\\
\Delta c=c\tens c,\quad \Delta e^i=e^a\tens  m^+{}^i{}_a c^{-1}+ 1\tens
e^a,\quad \Delta f_i= f_i\tens 1+ c m^-{}^a{}_i\tens f_a.}
\end{corol}
\proof We repeat the computations from (\ref{Uweakprim}), using this time the
coactions and weak quasitriangular structure of
$(\widetilde{U_q(\cg)},\widetilde{G_q})$ from Lemma~5.2. Because $\A_q$ and
$\bar{\A^\star_q}$ are well-defined braided groups in the corresponding braided
category, we know that the our previous calculations can be made at this level.
\endproof

The fundamental representation of $U$ on $\A_q$ also descends to this quotient
level. The formulae in Proposition~5.3 become
\eqn{extfund}{ v(x)\ra c=v(\lambda x),\quad v\ra e^i={vx^i- x^a (v\ra c^{-1}
m^+{}^i{}_a)\over q-q^{-1}},\quad v\ra f_i=-\del_i v}
and $\lambda R$ in place of $R$ for the action of $\vecm^\pm$. That $\del_i$
descend to the quotients (\ref{BDrel}) is shown in \cite{Ma:fre}. The first and
last actions in (\ref{extfund}) provide the fundamental representation of the
$q$-Poincar\'e algebra  in $q$-spacetime as a module algebra
c.f.\cite{Ma:poi}\cite{Ma:qsta}, for the appropriate regular $R$-matrix data
and quotients. Our double bosonisation construction extends this approach to
the $q$-conformal group defined by $U$, with the $e^i$ the additional
generators acting as in (\ref{extfund}). This geometrical picture of
(\ref{extfund}) will be developed elsewhere.

It should be clear that Corollary~5.5 leads to new quantum groups even when
$U_q(\cg)$ is one of the standard
$q$-deformations of ABCD type. For in these cases there is more than one
possible choice of quantum plane $\A_q$.
For the A series, there are two choices, namely of `fermionic' or `bosonic'
type. For generic $q$ the latter has the same dimensions at each degree as the
classical polynomial algebra in $n$ variables. There is such a standard choice
for each
of the ABCD series\cite{FRT:lie}. If we consider these  and work over  formal
power series $\C[[\hbar]]$ as in \cite{Dri}, then the pairing between $\A_q$
and $\bar{\A^\star_q}$ is non-degenerate since this is so `near' $q=1$, where
our algebras have a classical meaning (in this case the pairing via usual
differentiation). Hence we can expect that $\bar{\exp}_B$ exists as a formal
power series. In this deformation-theoretic setting we can also write
$c=\lambda^\xi$  and the weak quasitriangular structure in Lemma~5.4 becomes
$\CR_\xi=\lambda^{-\xi\tens\xi}$.  Then $U$ is quasitriangular with
\eqn{extR}{ \CR_U=\bar{\exp}_B\, \lambda^{-\xi\tens\xi}\,  \CR_{U_q(\cg)}.}
Moreover, from the relations in Corollary~5.5 we see that $U$ will also be a
deformation of a semisimple Lie algebra. This is because in the limit $q\to 1$
(in the sense of \cite{Dri}) we obtain $f_i$ and $e^i$ in the image of $[\xi,\
]$ and $\xi$ from the image of $[e^i,f_j]$. Hence we see that applied to
standard quantum groups with the standard `bosonic' choice of corresponding
quantum plane, the construction in Corollary~5.5
provides a way to construct quantum deformations of $U(\cg)$ by induction: the
induction step increases the rank
of $\cg$ by 1 and increases the dimension by $2n+1$, adjoining $\A_q$ to the
positive roots and $\bar{\A^\star_q}$ to the negative roots. Here $n$ is the
dimension of the defining representation $\cg\subseteq M_n$. From this, we
expect that the induction step takes a $q$-deformation of $U(sl_n)$ to one of
$U(sl_{n+1})$, a $q$-deformation of $U(so_n)$ to $U(so_{n+1})$ and a
$q$-deformation of $U(sp_n)$ to $U(sp_{n+1})$.

The same principle applies at our algebraic level. We demonstrate this now on a
concrete example, using $U_q(sl_n)$ in Lusztig's form in Section~4. For
technical reasons (order to use the known weak quasitriangular structure on
$U_q(sl_2)$ in the Drinfeld-Jimbo form) we adjoin the square roots $K^{\h}$ to
both input and output.

\begin{example} Let $\A_q$ be the standard bosonic quantum plane  of $sl_2$
type and suppose that $q$ has a square root in $k$. Let $\dot{U}_q(sl_2)$
denote
$U_q(sl_2)$ from Proposition~4.3 with $K^\h$ adjoined. This forms a weakly
quasitriangular dual pair with $SL_q(2)$ in a standard form. Then $U(\A_q,
\widetilde{\dot{U}_q(sl_2)}, \bar{\A_q^{\star}})$ is $U_q(sl_3)$ from
Proposition~4.3 with $K^\h$ adjoined.
\end{example}
\proof We start with $R$-matrix datum
\[ R=\pmatrix{q^2&0&0&0\cr0&q&q^2-1&0 \cr 0&0&q&0\cr0&0&0&q^2},\]
where the entry at row  $(ik)$ and column $(jl)$ (taken in the order
$11,12,21,22$) is $R^i{}_j{}^k{}_l$. We take  $R'=q^{-2}R$ in (\ref{BDrel}).
This gives the quantum-braided plane $\A_q^2$
with the usual relations $e^2e^1=qe^1e^2$ and the correct braiding
$\Psi$\cite{Ma:fre}. The dual $\bar{\A^\star_q}$ from (\ref{BDrel}) is similar,
with $f_2f_1=qf_1f_2$. The quantum group normalisation for $R$ needed   for the
(weak) quasitriangular structure on $\dot{U}_q(sl_2)$ in \cite{Ma:qua} is given
by $\lambda=q^{-3\over 2}$. The form of $\vecl^\pm$ in \cite{FRT:lie} provides
the weak quasitriangular structure explicitly as
\[ \vecm^+=\pmatrix{K^\h&-q^{-\h}(q-q^{-1})eK^{-\h}\cr 0& K^{-\h}},\quad
\vecm^-=\pmatrix{K^{-\h}&0\cr -q^\h(q-q^{-1})K^\h f&K^\h}\]
in our present conventions, where $e,f,K$ are the generators of $U_q(sl_2)$ in
Proposition~4.3.

We first compute the `Borel' relations between $\vecm^\pm$ (in the role of
`Cartan') and the $e^i$. Only five of the entries of $R$ are non-zero; in
particular the only off-diagonal entry is $R^1{}_2{}^2{}_1$. So we have
\cmath{ e^1m^+{}^1{}_1=\lambda R^1{}_1{}^1{}_1 m^+{}^1{}_1 e^1,\quad
e^1m^+{}^1{}_2=\lambda R^1{}_1{}^1{}_1 m^+{}^1{}_2 e^1,\quad
e^2m^+{}^2{}_2=\lambda R^2{}_2{}^2{}_2 m^+{}^2{}_2 e^2,\\
e^2m^+{}^1{}_2=\lambda R^1{}_1{}^2{}_2 m^+{}^1{}_2 e^2+ \lambda R^1{}_2{}^2{}_1
m^+{}^2{}_2 e^1,\\
m^-{}^2{}_1e^2=\lambda R^2{}_2{}^2{}_2 e^2m^-{}^2{}_1,\quad
m^-{}^2{}_1e^1=\lambda R^1{}_1{}^2{}_2 e^1m^-{}^2{}_1+  \lambda R^1{}_2{}^2{}_1
e^2m^-{}^1{}_1,}
with the other relations empty or redundant. From the form of $\vecm^\pm$ we
obtain
\cmath{ e^1K^\h=q^\h K^\h e^1,\quad e^1e=qee^1,\quad e^2K^{-\h}=q^\h
K^{-\h}e^2,\quad qee^2-e^2e=q^{-\h}e^1\\
{}[f,e^2]=0,\quad [f,e^1]=-q^{-\h}K^{-1}e^2}
correspondingly. The calculation for the $\vecm^\pm$ relations with $f_1,f_2$
is similar and the results analogous.
The remaining relations from Corollary~5.5 are clearly
\cmath{e^ic=q^{-3\over 2}ce^i,\quad cf_i=q^{-3\over 2}f_ic,\quad
[e^1,f_1]={K^\h c^{-1}-cK^{-\h}\over q-q^{-1}},\quad
[e^2,f_2]={K^{-\h}c^{-1}-cK^\h\over q-q^{-1}}\\
{}[e^1,f_2]=-q^{-\h}eK^{-\h}c^{-1},\quad [e^2,f_1]=q^\h cK^\h f\\
\Delta e^1=e^1\tens K^\h c^{-1}-q^{-\h}(q-q^{-1})e^2\tens eK^{-\h}c^{-1}+1\tens
e^1,\quad \Delta e^2=e^2\tens K^{-\h}c^{-1}+1\tens e^2 \\
\Delta f_1=f_1\tens 1+cK^{-\h}\tens f_1-q^\h(q-q)^{-1}cK^\h f\tens f_2,\quad
\Delta f_2=f_2\tens 1+cK^\h\tens f_2.}
Comparing with $U_q(sl_3)$ in Proposition~4.3 with Cartan matrix
$\pmatrix{2&-1\cr-1&2}$ we see that we can identify $e,f,K$ as the copy of
$U_q(sl_2)$ associated to  $i=1$ there, and we can identify
$e^2,f_2,K_2=K^{-\h}c^{-1}$ as the copy associated to $i=2$ there. We identify
the above elements $e^1,f_1$ from the quantum-braided plane as non-simple roots
generated by $q$-commutators with $e,f$.
We construct here $U_q(sl_3)$ with $K^\h$ adjoined.

We recover, in fact, more than just $U_q(sl_3)$ defined by generators and
relations (such as the $q$-Serre relations contained in Lusztig's algebra $f$
in Section~4) -- we are explicitly adjoining the non-simple root generators as
well. Their relations with the other root generators (which is the content of
the $q$-Serre relations) appear in our inductive approach from the quantum
plane relations and the `Cartan' relations with $\vecm^\pm$. These are all
provided by the inductive construction, along with their explicit coproducts
and cross relations.

It is possible to define a sub-Hopf algebra of $\widetilde{\dot{U}_q(sl_2)}$
generated by $\vecm^\pm c^{\mp 1}$, and a sub-Hopf algebra of
$\widetilde{SL_q(2)}$ generated by $\vect g$; our braided groups live more
precisely in the braided category generated by this weakly quasitriangular dual
pair. In this way, one can obtain precisely $U_q(sl_3)$ rather than its
extension by $K^\h$. \endproof

We see that our inductive procedure takes us from the standard quantisation to
the standard quantisation, at least in this example. In fact, it is clear from
another approach\cite{MaMar:glu}  that $U_q(sl_n)$ indeed becomes
$U_q(sl_{n+1})$ when we adjoin its fundamental bosonic quantum plane. This
approach\cite{MaMar:glu} provided an inductive construction of the $R$-matrix
datum and associated quantum matrices  by `gluing' quantum planes, but was
limited to the $A$ series or similar (Hecke type) R-matrices. Our present
approach is much more powerful and not limited in this way. From
(\ref{extfund}) it is also clear that $U_q(so_n)$ becomes $U_q(so_{n+2})$, i.e.
the rotation group is turned into the conformal group in the same dimension.
The same can be expected for the $C$ series in its standard $q$-deformation.

As a consequence of this inductive approach, we see that the quasitriangular
structure of $U_q(\cg)$ is then built up by iteration of (\ref{extR}) as a
product of braided exponentials in the roots on the left and `Gaussian' factors
of the form $\lambda^{-\xi\tens\xi}$  which can be collected to the right as
the Cartan part.
For example, in the deformation-theoretic setting of Example~5.6 we obtain a
quasitriangular structure
\eqn{expsl3}{ \CR_{sl_3}=
(\exp_B)_{21}^{-1}\lambda^{-\xi\tens\xi}\CR_{sl_2};\quad
\exp_B=\sum_{m=0}^\infty e^{i_1}\cdots e^{i_m}\tens f_{i_1}\cdots f_{i_m}
([m;q^2]!)^{-1},}
where the expression for $\exp_B$ follows immediately from the free form
(\ref{BDexp}) for any Hecke type R-matrix\cite{KemMa:alg}. Moreover, if each of
the quantum planes being adjoined has a natural basis, we build up by iteration
a natural `geometrical basis' of $U_q(\cg)$  using the triangular decomposition
in Section~3 at each step. Full details of this inductive construction of
$U_q(\cg)$ at least for the ABCD series will be developed in a sequel.

\appendix

\section{Appendix: Relation with the quantum double of a bosonisation}

Drinfeld in \cite{Dri} introduced a construction for a quasitriangular Hopf
algebra $D(H)$ associated to a (say, finite-dimensional)
Hopf algebra $H$. Here we use a generalised version $D(H,A)$ of
this\cite{Ma:mor} associated to a dual pair of Hopf algebras, i.e. two Hopf
algebras equipped with a (not necessarily non-degenerate) pairing. They can
even be bialgebras so long as the pairing is
convolution-invertible\cite{Ma:mor}. Here $D(H,A)$ is built on $H\tens A$ with
the tensor product coalgebra and unit, and the product
\eqn{DHA}{ (h\tens a)(g\tens b)=g\t h \tens a\t b \<Sa\o, g\o\>\<a\th,g\th\>.}
\note{There is a dual construction*** ${\rm co}D(A,H)$ with respect to which
$D(H,A)$ is weakly quasitriangular.}

If $H$ is a (weakly)  quasitriangular Hopf algebra with dual $A$ and
$B\in\CM^A, C\in {}_H\CM$ are dually paired braided groups in the sense of
a morphism $\ev:C\tens B\to k$ as explained in the Preliminaries, we have a
dual bosonisation $A\rbiprod B$ and bosonisation $C\lbiprod H$,
which are dually paired Hopf algebras\cite{Ma:mec}\cite{Ma:qsta} as in
(\ref{rcobos}). We can therefore form their generalised double.

\begin{lemma} Let $B\in \CM^A$ and $C\in{}_H\CM$ be categorically dually paired
braided groups, where $H$ is quasitriangular and dually paired with $A$. The
generalised quantum double $D(C\lbiprod H,A\rbiprod B)$ is a Hopf algebra
structure on $C\tens H\tens A\tens B$  with structure
\align{ &&\equad (c\tens h\tens a\tens b)\cdot(c'\tens h'\tens a'\tens b')=
(\Ro{1}\la c'\Bt)(\Rt{3}h'\t\la c)\tens \Rt{4}h'\th h\tens a\t (\Rt{6}h'\fiv\la
b\Bo)\bt a'\o\\
&&\qquad\tens  (\Ro{8}\Rt{7}h'\six\la
b\Bt)b'\<\Rt{1}\Rt{2}h'\o,Sa\o\>\<\Rt{8},a'\t\>\<\Rt{5}h'\fo,a\th\>   \\
&&\qquad \ev(c'\Bo,(\Rt{6}h'\fiv\la
b\Bo)\bo)\ev(\Ro{7}\Ro{6}\Ro{5}\Ro{4}\Ro{3}\Ro{2}\la \und S^{-1}
c'\Bth,h'\sev\la b\Bth),\\
&&\equad\Delta(c\tens h\tens a\tens b)= c\o\tens \Rt{} h\o\tens a\o\tens
b\Bo\bo\tens\Ro{}\la c\Bt\tens h\t\tens a\t b\Bo\bt\tens b\Bt.}
Similarly when $H,A$ are a weakly quasitriangular dual pair. In the
finite-dimensional non-degenerately paired case we have Drinfeld's
quasitriangular structure as
\[ \CR_D=(f^\alpha\cu\Ro{}\la\und S f^a\tens f^\beta S\Rt{}\tens 1\tens
1)\tens(1\tens 1\tens e_a\tens e_\alpha e_\beta)\]
where $\{e_a\}$ is a basis of $B$ with dual $\{f^a\}$ and $\{e_\alpha\}$ a
basis of $A$ with dual $\{f^\alpha\}$.
\end{lemma}
\proof This is a straightforward but tedious calculation. We emphasise the case
where $H$ is quasitriangular because the notation is more familiar; as
explained in Remark~3.9, our proofs and results convert immediately over to the
weakly quasitriangular case. The structure of $C\lbiprod H$ is in (\ref{lbos})
and its dual construction $A\rbiprod B$ is in (\ref{rcobos}). We need their
iterated coproducts as
\cmath{(\id\tens\Delta)\circ\Delta(c\tens h)=c\Bo\tens\Rt{1}\Rt{2}h\o\tens
\Ro{1}\la c\Bt\tens\Rt{3}h\t\tens\Ro{3}\Ro{2}\la c\Bth\tens h\th\\
(\Delta\tens\id)\circ\Delta(a\tens b)=a\o\tens b\Bo\bo\tens a\t b\Bo\bt\o\tens
b\Bt\bo\tens a\th b\Bo\bt\t b\Bt\bt\tens b\Bth.}
Then
\align{&&\equad (c\tens h\tens a\tens b)\cdot(c'\tens h'\tens a'\tens b')=\\
&&=(c'\tens h')\t\cdot (c\tens h)\tens (a\tens b)\t\cdot (a'\tens b')\<S(a\tens
b)\o,(c'\tens h')\o\>\<(a\tens b)\th,(c'\tens h')\th\>\\
&&=(\Ro{1}\la c'\Bt\tens \Rt{3}h'\t)\cdot(c\tens h)\tens (a\t\tens
b\Bo\bt\o\tens b\Bt\bo)\cdot(a'\tens b')\\
&&\qquad\<S(\Rt{1}\Rt{2}h'\o),a\o\>\<\Rt{4}h'\th,a\th b\Bo\bt\t b\Bt\bt\>\\
&&\qquad\ev(c'\Bo,b\Bo\bo)\ev(\Ro{4}\Ro{3}\Ro{2}\la \und S^{-1}c'\Bth,h'\fo\la
b\Bth)\\
&&=(\Ro{1}\la c'\Bt)(\Rt{3}h'\t\la c)\tens\Rt{4} h'\th h\tens a\t b\Bo\bt\o
a'\o\tens b\Bt\bo b'\\
&&\qquad\<\Rt{1}\Rt{2}h'\o,Sa\o\>\<\Rt{5}h'\fo,a\th b\Bo\bt\t
b\Bt\bt\t\>\<\Ro{6},b\Bt\bt\o\>\<\Rt{6},a'\t\>\\
&&\qquad\ev(c'\Bo,b\Bo\bo)\ev(\Ro{5}\Ro{4}\Ro{3}\Ro{2}\la \und S^{-1}c'\Bth,
h'\fiv\la b\Bth)}
using the definition (\ref{DHA}) applied to our case, the triple coproducts
above and then the products from (\ref{lbos}) and (\ref{rcobos}). We also used
invariance of $\ev$ in the form $\ev(c,h\la b)=\ev(S^{-1}h\la c,b)$ and the
properties of actions and coactions. Further rewriting the right coaction as a
left action $b\bo \<h,b\bt\>=h\la b$, gives the formula stated. The coproduct
of the quantum double is the tensor product one from
(\ref{lbos})--(\ref{rcobos}).

For the quasitriangular structure in the case $A=H^*, C=B^*$, we see from the
duality pairing (\ref{rcobos}) that if $e_{\alpha,a}=e_\alpha\tens e_b$ is a
basis of $A\tens B$ then
\align{&&\equad f^{\alpha,a}=S(f^a\tens
S^{-1}f^\alpha)=(S(S^{-1}f^\alpha)\t)\cu\Ro{}\la \und S f^a\tens
S(\Rt{}(S^{-1}f^\alpha)\o)\\
&&=f^\alpha\o\cu \Ro{}\la\und S f^a\tens f^\alpha\t S\Rt{}}
is a dual basis under the pairing. We computed the antipode (\ref{lbos}) of
$C\lbiprod H$ from \cite{Ma:bos}. Then Drinfeld's quasitriangular structure
$\CR=f^{\alpha,a}\tens 1\tens 1\tens e_{\alpha,a}$ giving the formula shown. We
have emphasised this case. When $A,H$ and $C,B\in \CM^A$ are merely dually
paired then the quantum double here is weakly quasitriangular with respect to
the natural codouble construction made along the same lines as above. \endproof

We now apply this general construction to the specific setting of Section~3.
There we had a braided group $B\in {}^A\CM$ which we consider as $B\in
\CM^{A^{\rm cop}}$, and categorically dually paired with it a braided group $C$
in ${}\CM_H$ which we consider as $C\in {}_{H^{\rm op}}\CM$.

\begin{theorem} In Section~3, there is a Hopf algebra surjection
$\pi:D(C\lbiprod H^{\rm op},A^{\rm cop}\rbiprod B)\to U(\bar C,H,B)$ from the
generalised quantum double, given by
\[ \pi(c\tens h\tens a\tens b)=c\ra Sh\o\tens h\t\Rt{}\<a,\Ro{}\>\tens b.\]
In the finite-dimensional non-degenerately paired case this is a surjection of
quasitriangular Hopf algebras.
\end{theorem}
\proof First, we write out the product from Lemma~A.1 applied to our particular
dual pair. We equip $H^{\rm op}$ with antipode $S^{\rm op}=S^{-1}$ and
quasitriangular structure $\CR^{\rm op}=\CR_{21}$. Then we write everything in
terms of the usual structure of $H$ and the usual coproduct of $A$. We write
the left action of $H^{\rm op}$ on $C$ in its original form as a right action
of $H$ (as in Section~3), and the right coaction of $A^{\rm cop}$ in its
original form as a left coaction of $A$ or (by evaluation) a right action of
$H$. These are simple manipulations on the  formulae in Lemma~A.1.  The
resulting product between various combinations of the tensor factors (embedded
in $C\tens H\tens A\tens B$ in the trivial way by tensoring with $1$) becomes
\ceqn{Drels}{ (h\tens a)(h'\tens a')=hh'\t\tens a\t
a'\<h'\th,a\o\>\<h'\o,S^{-1}a\th\>\\
(c\tens h)(c'\tens h')=c'(c\ra h'\o)\tens hh'\t,\quad (a\tens b)(a'\tens b')=a
a'\t\tens (b\ra\Rt{}) b'\<\Ro{},a'\o\>\\
(c\tens a)(c'\tens a')=(c'\ra\Rt{})c\tens a\o a'\<\Ro{},S^{-1}a\t\>,\quad
(h\tens b)(h'\tens b')=hh'\o\tens (b\ra h'\t)b'\\
(c\tens b)(c'\tens b')=c'\Bt(c\ra\Ro{1})\tens\Ro{2}\tens
(b\Bo\ra\Ro{3})\bo\tens (b\Bt\ra\Ro{4})b'\\
\qqquad\qqquad\qqquad\quad \ev(c'\Bo,(b\Bo\ra\Ro{3})\bt)\ev(\und
S^{-1}c'\Bth\ra\Rt{1}\Rt{2}\Rt{3}\Rt{4},b\Bth),}
while the coproduct becomes
\eqn{Dcop}{\Delta(c\tens h\tens a\tens b)=c\Bo\tens h\o\Ro{}\tens a\t\tens
b\Bo\bt\tens c\Bt\ra\Rt{}\tens h\t\tens a\o b\Bo\bo\tens b\Bt.}
We recognise a version of the generalised quantum double built on $H\tens
A^{\rm cop}$  as a sub-Hopf algebra. We also recognise sub-Hopf algebras
$(C\lbiprod H^{\rm op})^{\rm op}$ and $A^{\rm cop}\rbiprod B$ as expected.
Finally, we recognise that $\HB\subseteq U$ appears here as a subalgebra, as
does $\CH$ up to an elementary isomorphism.

We verify now that $\pi$ is a Hopf algebra map to $U$ from Section~3 (it is
clearly surjective). The restriction  to the generalised quantum double of
$H,A$ is $\pi(h\tens a)=h\Ro{}\<\Rt{},a\>$ which is a version of the Hopf
algebra projection $D(H)\to H$ introduced in \cite{Ma:dou} when $H$ is
quasitriangular. The restriction to $C\tens H$ is $\pi(c\tens h)=c\ra Sh\o\tens
h\t=h\o\la c\tens h\t$ and gives an isomorphism to the algebra $\CH\subseteq
U$. We use elementary properties of Hopf algebras and the conversion of the
action on $\bar C$ as in (\ref{barca}). The  restriction of $\pi$  to $H\tens
B$ is the identity and hence immediately an algebra map to $U$. For the
remaining restrictions, we have
\align{ \pi((a\tens b)(a'\tens b'))\equad&& =\<aa'\t,\Ro{1}\>\Rt{1}\tens
(b\ra\Rt{3})b'\<\Ro{3},a'\o\>\\
&&=\<a,\Ro{1}\>\<a',\Ro{3}\Ro{2}\>\Rt{1}\Rt{2}\tens (b\ra\Rt{3})b'\\
&&=\<a,\Ro{1}\>\<a',\Ro{2}\>(\Rt{1}\tens b)\cdot(\Rt{2}\tens b')=(\pi(a\tens
b))(\pi(a'\tens b')),}
using the elementary properties of the quasitriangular structure and the
relations of $\HB\subseteq U$. Likewise, we have
\align{\pi((c\tens a)(c'\tens a'))\equad&& =(c'\ra\Rt{1})\cdot_C c\tens
\Rt{2}\<a\o a',\Ro{2}\>\<\Ro{1},S^{-1}a\t\>\\
&&=c(\Rt{1}\la c')\tens \Rt{2}\Rt{3}\<a,\Ro{2}\Ro{1}\>\<a',\Ro{3}\>\\
&&=(c\tens\Rt{1})\cdot(c'\tens\Rt{3})\<\Ro{1},a\>\<\Ro{3},a'\>=(\pi(c\tens
a))(\pi(c'\tens a')),}
where the latter products are in $\bar C$ (which has the opposite algebra
structure to $C$) and $H$, or, finally,  in $\CH\subseteq U$. We used
invariance of $\CR$ under $S\tens S$. For the final restriction, we need the
iteration of the relation (\ref{barca}) between the coproducts of $C,\bar C$ as
\eqn{barcd}{(\bar \Delta\tens\id)\circ\bar\Delta c= c\Bo \ra
\Rmo{3}\Rmo{8}\tens c\Bt\ra\Rmo{7}\Rmt{8}\tens c\Bth\ra\Rmt{7}\Rmt{3},}
where the numbering of the copies of $\CR^{-1}$ is to keep them distinct from
other copies used in the proof.
Then
\align{&&\equad \pi((c\tens b)(c'\tens b'))=(c'\Bt\cdot_C (c\ra\Ro{1}))\ra
S\Ro{2}\tens \Ro{5}\Rt{6}\<\Ro{6},(b\Bo\ra\Ro{3})\bo\>\tens (b\Bt\ra\Ro{4})b'\\
&&\qquad\ev(c'\Bo,(b\Bo\ra\Ro{3})\bt)\ev(\und
S^{-1}c'\Bth\ra\Rt{1}\Rt{2}\Rt{5}\Rt{3}\Rt{4},b\Bth)\\
&&=(c'\Bt\ra\Rmo{7})\cdot_Cc\tens \Ro{6}\Ro{5}\tens (b\Bt\ra\Ro{4})b'\\
&&\qquad\ev(c'\Bo,b\Bo\ra\Ro{6}\Ro{3})\ev(\und
S^{-1}c\Bth\ra\Rmt{7}\Rt{3}\Rt{5}\Rt{4},b\Bth)\\
&&=c(c'\bart\ra\Rt{8})\tens \Rt{6}\Ro{5}\tens (b\Bt\ra\Ro{4})b'\\
&&\qquad \ev(c'\baro\ra\Ro{8},b\Bo\ra\Ro{6})\ev(\bar S
c'\barth\ra\Rt{5}\Rt{4},b\Bth)\\
&&=(\pi(c\tens b))(\pi(c'\tens b')).}
For the second equality we evaluated the coaction on $b\Bo\ra\Ro{3}$ as an
action on it of $\Ro{6}$, cancelled $\CR_1\CR_2^{-1}$ and then used the QYBE
applied to $\CR_3,\CR_6,\CR_5$. The third equality is invariance of $\ev$ and
(\ref{barcd}). We also adopt the product of $\bar C$. Making the further
notational change from $\ev$ to $\<\ ,\ \>$ as in Section~3, using invariance
under $\Rt{4}$ and writing the action of $\Ro{8}$ from the right as $S\Ro{8}$
from the left by (\ref{barca}) allows us to recognise the product in $U$ as the
final step.

To see that $\pi$ is a coalgebra map, we compute
\align{&&\equad (\pi\tens\pi)\Delta(c\tens h\tens a\tens b)=c\Bo\ra S\Ro{1}
Sh\o\tens h\t\Ro{2} h\t \Rt{3}\tens b\Bo\bt \<a\t,\Ro{3}\>\\
&&\qqquad\qqquad \qqquad\tens c\Bt\ra\Rt{1}\Rt{2}Sh\th\tens h\fo\Rt{4}\tens
b\Bt\<a\o b\Bo\bo,\Ro{4}\>\\
&&=h\o\la c\baro\tens \Ro{2}h\th \Rt{3}\tens b\Bo\ra
\Ro{5}\<a,\Ro{4}\Ro{3}\>\tens h\t S^{-1}\Rt{2}\la c\bart\tens
h\fo\Rt{4}\Rt{5}\tens b\Bt\\
&&=\Delta_U(h\o\la c\tens h\t\Rt{}\<a,\Ro{}\>\tens b)=\Delta_U\circ\pi(c\tens
h\tens a\tens b)}
as required. The first equality is the definitions and elementary properties
(\ref{qua}). The second equality is (\ref{barca}). We also write the action of
$Sh\o$ on $C$ from the right as $h\o$ from the left, etc. In addition, we use
the quasicocommutativity axiom (\ref{qua}) to $h\t\tens h\th$. We then identify
the result in terms of the coproduct of $U$.

Finally, in the finite-dimensional non-degenerately paired case we have the
quasitriangular structure converted from Lemma~A.1 (as explained above) as
\[ \CR_D=(\und S f^a\ra (S\Rt{})\cu f^\alpha\tens \Ro{} f^\beta\tens 1\tens
1)\tens(1\tens 1\tens e_\alpha e_\beta\tens e_a).\]
We used invariance of $\CR$ under $S\tens S$. Then
\align{&&\equad (\pi\tens\pi)(\CR_D)=\und S f^a\ra (S\Rt{1})\Rt{2}\cu f^\alpha
(Sf^\beta) \Ro{2}\tens\Ro{1}f^\gamma\tens 1\tens 1\tens\Rt{3}\<e_\alpha e_\beta
e_\gamma,\Ro{3}\>\tens e_a\\
&&=\Rt{1}\la \bar S^{-1}f^a\tens\Ro{1}\Ro{3}\tens 1\tens 1\tens\Rt{3}\tens
e_a=\CR_U}
using the definition of $\pi$, the antipode axioms followed by
$\Rt{}\cu\Ro{}=1$. We write the action of $S\Rt{1}$ from the right as an action
from the left. Using invariance of the coevaluation $f^a\tens e_a$ and changing
to a new basis $f'{}^a=\bar S^{-1}f^a$ with dual $\bar Se_a'$ gives $\CR_U$
from Proposition~3.6. \endproof

\section{Appendix: Double biproducts}

This appendix introduces a `double biproduct' construction which generalises
the  double bosonisation in Section~3.  We recall that Radford in
\cite{Rad:str} considered as `biproducts' the general class of Hopf algebras
which are a smash (or cross) product and
coproduct by the simultaneous action and coaction of a Hopf algebra $H$; it is
clear from \cite{Ma:skl}\cite{Ma:dou} that bosonisations can be viewed as
examples of this general form (see the Preliminaries).  We extend this
observation now to our `doubled' setting by introducing the required `double'
construction. This provides an alternative point of view because the focus is
now on actions and coactions (rather than on quasitriangular and
dual-quasitriangular structures), which should be more accessible to the
algebraically minded reader. The natural setting in the present section is with
$H$ a general Hopf algebra with bijective antipode. The antipode and its
inverse are not actually needed in the main construction, i.e. the result is
slightly more general.

The correct setting for biproducts was identified in \cite{Ma:dou}\cite{Ma:skl}
as the braided category of crossed modules $\CM^H_H$, as discussed
independently in \cite{Yet:rep} in another context. This category is in fact
nothing other than a version of the braided category of $D(H)$-modules
introduced through the work of V.G. Drinfeld. Simply, one casts the action of
$H^*\subseteq D(H)$ as a coaction of $H$; this is then well-defined even for
$H$ infinite-dimensional. An object in $\CM^H_H$ is a vector space on which $H$
acts and coacts from the right such that
\eqn{rcrossmod}{v\bo\ra h\o\tens v\bt h\t=(v\ra h\t)\bo\tens h\o(v\ra h\t)\bt.}
A morphism is a linear map intertwining both the action and coaction of $H$.
There is also a left-handed version ${}^H_H\CM$ with
\eqn{lcrossmod}{h\o v\bo\tens h\t\la v\bt=(h\o\la v)\bo h\t\tens (h\o\la
v)\bt.}
Both categories are braided, with
\eqn{lrcrossbraid}{\Psi_{V,W}(v\tens w)=w\bo\tens v\ra w\bt,\quad
\Psi_{V,W}(v\tens w)=v\bo\la w\tens v\bt}
for the two cases. These are just the braidings corresponding to Drinfeld's
quasitriangular structure on the quantum double. As explained by the author in
\cite{Ma:dou}\cite{Ma:skl}, braided groups in such categories exactly satisfy
the conditions in \cite{Rad:str} to make a simultaneous cross product and cross
coproduct by the given action and coaction, and obtain an ordinary Hopf
algebra. We consider braided groups $B\in \CM^H_H$ and $\bar C\in {}^H_H\CM$
and denote the corresponding biproducts by $H\rbiprod B$ and $\bar C\lbiprod H$
as an extension of our previous notation. We assume, moreover, that there is a
`pairing' $\<\ ,\ \>:\bar C\tens B\to k$ such that
\ceqn{bippaira}{ \<h\la c,b\>=\<c,b\ra h\>\\
\<c,ab\>=\<c\bart,a\ra b\bt\>\<c\baro,b\bo\>,\quad
\<cd,b\>=\<c,b\Bo\>\<d,b\Bt\>,\quad \<\bar S c,b\>=\<c,\und S^{-1}b\>}
for all $h\in H$, $a,b\in B$ and $c,d\in \bar C$. We deduce (from braided
antimultiplicativity of braided antipodes\cite{Ma:tra}) that
\eqn{bippairb}{\<\bar S c,ab\>=\<\bar S c\baro,a\>\<\bar S c\bart,b\>,\quad
\<\bar S(cd),b\>=\<c\bo\la \bar S d,b\Bo\>\<\bar S c\bt,b\Bt\>}
hold as well. The first condition in (\ref{bippaira}) expresses bicovariance of
the pairing under the action. In place of bicovariance under the coaction,
however, we adopt a compatibility condition
\eqn{bipcond}{b\bo\ra c\bo\tens b\bt\la c\bt=b\tens c,\qquad \forall b\in B,\
c\in\bar C}
between the two crossed module structures. This is the data needed for the
following construction.

\begin{theorem} Let $H$ be a Hopf algebra with invertible antipode,
$B\in\CM^H_H$, $\bar C\in{}^H_H\CM$  braided groups obeying (\ref{bipcond}),
and $\<\ , \>$ a braided  skew-pairing as above. Then there is a unique
ordinary Hopf algebra $\CU(\bar C,H,B)$ built on $\bar C\tens H\tens B$, the
{\em double biproduct}, containing $H\rbiprod B$ and $\bar C\lbiprod H$ as
sub-Hopf algebras with cross relations
\[ bc=b\Bo\bt c\bart b\Bt c\barth\bo \<c\baro,b\Bo\bo\>\<\bar S
c\barth\bt,b\Bth\>\]
\end{theorem}
\proof We provide an outline, following the same strategy as in the proofs of
Theorem~3.2. As before, we  work with the cross relations in a more canonical
form
\eqn{bipord}{bc=(b\Bo\bt\o\la c\bart)b\Bo\bt\t c\barth\bo\o (b\Bt\ra
c\barth\bo\t)\<c\baro,b\Bo\bo\>\<\bar S c\barth\bt,b\Bth\>,}
in view of the relations of $\HB$ and $\bar C\lbiprod H$. This has the form
$bc=\sum c_i R_i b_i$ for $c_i\in\bar C$, $R_i\in H$ and $b_i\in B$. The
general product is defined from this as in Section~3, and comes out as
\eqn{bipprod}{(chb)\cdot(dga)= c(h\o b\Bo\bt\o\la d\bart)(h\t b\Bo\bt\t
d\barth\bo\o g\o)(b\Bt\ra d\barth\bo\t g\t)a\<d\baro,b\Bo\bo\>\<\bar S
d\barth\bt,b\barth\>.}
The argument that it is enough to prove associativity of the products
$(a\cdot(chb))\cdot d=a\cdot((chb)\cdot d)$ goes through in the same way  (it
requires only covariance of the algebras of $B,\bar C$ under $H$). To prove
this special case we compute both sides from (\ref{bipprod}). The left hand
side is
\align{&&\equad (a\cdot(chb))\cdot d=\left((a\Bo\bt\o\la c\bart)(a\Bo\bt\t
c\barth\bo\o h\o)((a\Bt\ra c\barth\bo\t h\t)b))\right)\cdot d\<c\baro,
a\Bo\bo\>\<\bar S c\barth\bt,a\Bth\>\\
&&=(a\Bo\bt\o\la c\bart)(a\Bo\bt\t c\barth\bo\o h\o ((a\Bt\ra c\barth\bo\th
h\th)b)\Bo\bt\o\la d\bart)\\
&&\quad\tens a\Bo\bt\th c\barth\bo\t h\t ((a\Bt\ra c\barth\bo\th
h\th)b)\Bo\bt\t d\barth\bo\o\tens ((a\Bt\ra c\barth\bo\th h\th)b)\Bt\ra
d\barth\bo\t\\
&&\quad \<c\baro, a\Bo\bo\>\<\bar S c\barth\bt,a\Bth\>\<d\baro,((a\Bt\ra
c\barth\bo\th h\th)b)\Bo\bo\>\<\bar S d\barth\bt,((a\Bt\ra c\barth\bo\th
h\th)b)\barth\>\\
&&=(a\Bo\bt\o\la c\bart)(a\Bo\bt\t c\barth\bo\o h\o (a\Bt\ra c\barth\bo\th
h\th)\bt\o b\Bo\bt\t \la d\barth)\\
&&\quad\tens a\Bo\bt\th c\barth\bo\t h\t (a\Bt\ra c\barth\bo\th h\th)\bt\t
b\Bo\bt\th d\barfo\bo\o d\barfiv\bo\o\\
&&\quad \tens (a\Bt\ra c\barth\bo\fo h\fo b\Bo\bt\fo d\barfo\bo\t
d\barfiv\bo\t)(b\Bt\bo\ra d\barfo\bo\th d\barfiv\bo\th)\<c\baro,
a\Bo\bo\>\<\bar S c\barth\bt,a\Bfiv\>\\
&&\quad \<d\bart,(a\Bt\ra c\barth\bo\th h\th)\bo\ra
b\Bo\bt\o\>\<d\baro,b\Bo\bo\>\<\bar S d\barfo\bt, a\Bfo\ra b\Bo\bt\fiv
b\Bt\bt\>\<\bar S d\barfiv\bt,b\Bth\>,}
where the first two equalities are two applications of the product from
(\ref{bipprod}). The third equality then puts in the iterated braided coproduct
of a product in $B$ from the first half of
\ceqn{triplecoprod}{(\und\Delta\tens\id)\circ\und\Delta(ab)=a\Bo b\Bo\bo\tens
(a\Bt\ra b\Bo\bt\o)b\Bt\bo\tens (a\Bth\ra b\Bo\bt\t b\Bt\bt)b\Bth\\
(\id\tens\bar \Delta)\circ\bar\Delta(cd)=c\baro(c\bart\bo c\barth\bo\o\la
d\baro)\tens c\bart\bt (c\barth\bo\t\la d\bart)\tens c\barth\bt d\barth.}
We also use covariance of the products and coproducts under the action and
coaction of $H$, and break down the pairing with products from $B$ using
(\ref{bippaira})--(\ref{bippairb}). This gives the expression above. We now use
bicovariance of the pairing to move $\ra b\Bt\bt$ to act on $d\barfo\bt$ and
(\ref{bipcond}) to cancel it, and we use the crossing condition
(\ref{rcrossmod}) on $\Delta(c\barth\bo\o h\o (a\Bt\ra c\bart\bo\th
h\t)\bt)\tens (a\Bt\ra c\bart\bo\th h\t)\bo$. These steps give
\align{&&\equad (a\cdot(chb))\cdot d=(a\Bo\bt\o\la c\bart)(a\Bo\bt\t a\Bt\bt\o
c\barth\bo\t h\t b\Bo\bt\t \la d\barth)\\
&&\quad\tens a\Bo\bt\th  a\Bt\bt\t c\barth\bo\th h\th b\Bo\bt\th d\barfo\bo\o
d\barfiv\bo\o\\
&&\quad \tens (a\Bt\ra c\barth\bo\fo h\fo b\Bo\bt\fo d\barfo\bo\t
d\barfiv\bo\t)(b\Bt\ra  d\barfiv\bo\th)\<c\baro, a\Bo\bo\>\<\bar S
c\barth\bt,a\Bfiv\>\\
&&\quad \<d\bart,a\Bt\bo\ra c\barth\bo\o h\o b\Bo\bt\o\>
\<d\baro,b\Bo\bo\>\<\bar S d\barfo\bt, a\Bfo\ra b\Bo\bt\fiv\>\<\bar S
d\barfiv\bt,b\Bth\>.}
The calculation for $a\cdot((chb)\cdot d)$ is strictly similar: we use the
second line of (\ref{triplecoprod}) and the other halves of
(\ref{bippaira})--(\ref{bippairb}). Indeed, there is a strict symmetry
involving: reversal of all products and tensor products (reflection in a
mirror), interchange of $a,b,\la$ with $d,c,\ra$, interchange of $\<\ ,\ \>$
with $\<\bar S\ ,\ \>$ and reversal of the numbering of the coactions and all
coproducts. On the other hand, our final expression for $(a\cdot(chb))\cdot d$
is self-symmetric under this operation. Hence it coincides with the result for
$a\cdot((chb)\cdot d)$. Hence associativity is proven.

For the coalgebra structure, our requirement that $\HB$ and $\bar C\rbiprod H$
are sub-bialgebras forces the general coproduct to be
\eqn{bipcoprod}{ \Delta(chb)=c\baro\tens c\bart\bo h\o\tens b\Bo\bo\tens
c\bart\bt\tens h\t b\Bo\bt\tens b\Bt.}
As before, we only need to prove the bialgebra homomorphism property for the
special case $\Delta(b\cdot c)$. For brevity, we outline the proof using the
cross relations stated in the theorem rather than the more explicit ordering
relations (\ref{bipord}). This is equivalent, although less direct. Thus,
working in the algebra $\CU$, we have
\align{&&\equad (\Delta b)(\Delta c)=b\Bo\bo c\baro c\bart\bo\tens b\Bo\bt b\Bt
c\bart\bt\\
&&=b\Bo\bt\o c\bart b\Bt\bo c\barth\bo c\barfo\bo c\barfiv\bo
c\barsix\bo\o\tens b\Bo\bt\t b\Bt\bt b\Bth\bt b\Bfo\bt c\barfiv\bt b\Bfiv
c\barsix\bo\t\\
&&\quad \<c\baro,b\Bo\bo\>\<\bar S
c\barth\bt,b\Bth\bo\>\<c\barfo\bt,b\Bfo\bo\>\<\bar S c\barsix\bt,b\Bsix\> \\
&&=b\Bo\bt\o c\bart b\Bt\bo c\barth\bo c\barfo\bo c\barfiv\bo\o\tens b\Bo\bt\t
b\Bt\bt b\Bth\bt  c\barfo\bt b\Bfo c\barfiv\bo\t\\
&&\quad \<c\baro,b\Bo\bo\>\<(\bar S
c\barth\bt\baro)c\barth\bt\bart,b\Bth\bo\>\<\bar S c\barfiv\bt,b\Bfiv\>\\
&&=b\Bo\bt\o c\bart b\Bt\bo c\barth\bo c\barfo\bo\o  \tens b\Bo\bt\t b\Bt\bt
c\barth\bt b\Bth c\barfo\bo\t\<c\baro,b\Bo\bo\> \<\bar S c\barfo\bt,b\Bfo\>\\
&&=\Delta(bc),}
where the second equality uses the cross relations in each factor of
$\CU\tens\CU$ and covariance of the braided coproducts under the coaction of
$H$. The third equality uses the pairing axiom (\ref{bippaira}) with
$b\Bth\bo\tens b\Bfo\bo\tens b\Bth\bt b\Bfo\bt=\und\Delta (b\Bth\bo)\tens
b\Bth\bt$. We then use the braided-antipode property in $\bar C$. Finally we
recognise $\Delta(bc)$ in using the further identity
\eqn{bipcom}{ b\bo c\bo\tens b\bt c\bt=c\bo b\bo\tens c\bt b\bt,\quad \forall
b\in B,\ c\in\bar C,}
which follows directly from (\ref{bipcond}). As before, the antipode on $\CU$
exists and is uniquely determined by the antipodes of $\HB$ and $\bar C\lbiprod
H$. \endproof

We remark that we do not need the antipode or inverse antipode of $H$ for this
construction, though this is the natural setting for the input data. Without
the inverse antipodes, the braidings $\Psi$ in (\ref{lrcrossbraid}) are not
invertible, but the braided-homomorphism properties as in (\ref{triplecoprod})
still make sense. We do not need braided antipodes on $B,\bar C$ either, but
only need to assume a `convolution inverse' to $\<\ ,\ \>$ in place of $\<\bar
S(\ ),\ \>$, as characterised by (\ref{bippaira})--(\ref{bippairb}). In this
way, the above construction lifts entirely to the bialgebra level. We can also
write the relations in the theorem as
\eqn{biprelb}{b\Bo c\baro c\bart\bo\<c\bart\bt,b\Bt\>=\<c\baro,b\Bo\bo\>b\Bo\bt
c\bart b\Bt.}

This is an extension of  the construction of Section~3. The precise inclusion
of one construction in the other is provided by the functors
${}_H\CM\hookrightarrow {}^H_H\CM$  and $\CM_H\hookrightarrow \CM^H_H$ when $H$
is quasitriangular. These functors (introduced by the author in\cite{Ma:dou})
use the quasitriangular structure to induce from an action a compatible
coaction, forming a crossed module. In the specific setting of Section~3 we
have $B\in\CM_H\hookrightarrow \CM^H_H$ and $\bar C\in{}_{\bar
H}\CM\hookrightarrow {}_H^H\CM$ by induced right, left coactions
\eqn{quafunc}{ b\bo\tens b\bt=b\ra \Ro{}\tens \Rt{} ,\quad c\bo\tens
c\bt=\Rmo{}\tens\Rmt{}\la c}
respectively. In the weak quasitriangular case we have $B\in{}^A\CM$ and define
both the action of $H$ and the induced coaction of $H$ from this data as
explained in Remark~3.9 (i.e. the action is by evaluation against the coaction
of $A$, and the induced coaction is $b\mapsto b\bt\tens\CR(b\bo)$. Similarly
for $\bar C$. These inclusions are the standard way that bosonisations can be
viewed as examples of biproducts. It is easy to see that the condition
(\ref{bipcond}) also holds by cancellation of $\CR,\CR^{-1}$, so the
double-bosonisation $U$ in Section~3 can be viewed  as an example of the more
general $\CU$ above.

The double biproduct construction includes other examples as well. Thus, the
dual bosonisation construction is obviously also a biproduct (see the
Preliminaries), this time using the dual quasitriangular structure $\CR:H\tens
H\to k$ to induce an action from a coaction. In the present case we can use
these to map braided groups  $B\in\CM^H\hookrightarrow \CM^H_H$ and $\bar C\in
{}^{\bar H}\CM\hookrightarrow {}^H_H\CM$ by
\eqn{dquafunc}{b\ra h=b\bo\CR(b\bt\tens h),\quad h\la c=\CR^{-1}(h\tens
c\bo)c\bt}
respectively. Here $\bar H$ denotes $H$ with the inverse-transpose dual
quasitriangular structure.

\begin{example} Let $G$ be an Abelian group, $\beta:G\times G\to k$ an
invertible symmetric bicharacter and $B,D$
two $G$-graded Hopf algebras dually-paired in the usual way (respecting the
grading). Let $\bar C=D^{\und{\rm cop}}$. Then there is a
double-biproduct $\CU=\CU(\bar C,kG,B)$. For any $g\in G$ and $b,c$
braided-primitive of homogeneous $G$-degree $|\ |$, $\CU$ has the structure
\cmath{ g^{-1}cg=\beta(g,|c|)c,\quad g^{-1}bg=\beta(|b|,g)b,\quad
[b,c]=(|b|-|b|^{-1})\<c,b\>\\
\Delta g=g\tens g,\quad \Delta b=b\tens |b|+1\tens b,\quad \Delta c=c\tens
1+|c|\tens c\\
\eps g=1,\quad \eps b=\eps c=0,\quad Sg=g^{-1},\quad Sb=-b|b|^{-1},\quad
Sc=-|c|^{-1}c.}
\end{example}
\proof It is well known that group-graded algebras can be considered as
$kG$-comodule algebras, where $kG$ is the group algebra of $G$, see e.g.
\cite{Coh:hop}. $G$-graded Hopf algebras can be considered in the same way if
there is a bicharacter $\beta$, which we extend by linearity to a
dual-quasitriangular structure on $H=kG$ cf.\cite{Ma:csta}\cite{Ma:any}.
$G$-graded Hopf algebras thereby become Hopf algebras in the braided category
of comodules of a dual-quasitriangular Hopf algebra as introduced (by the
author) in \cite{Ma:bg}. We suppose that $B,D$ are such $G$-graded braided
groups, equipped with an ordinary pairing $\<\ ,\ \>:D\tens B\to k$ preserving
the degree in the sense $\<c,b\>|b|=|c|^{-1}\<c,b\>$ on all homogeneous
elements $b\in B$, $c\in D$. We view $B$ in the category of right
$kG$-comodules by $b\mapsto b\tens |b|$, and $D$ in the category of left
$kG$-comodules by $c\mapsto |c|\tens c$. Finally, we let $\bar C$ be the same
braided group as $D$ but with braided-opposite coproduct $\bar\Delta$. Then
$\bar C$ lives in the braided category of left $kG$-comodules with braiding
determined by the inverse-transposed bicharacter. Explicitly,
\[\und\Delta (ab)=\beta(|a\Bt|,|b\Bo|)a\Bo b\Bo\tens a\Bt b\Bt,\quad \bar\Delta
(cd)=\beta^{-1}(|c\bart|,|d\baro|)c\baro d\baro\tens c\bart d\bart\]
on homogeneously decomposed coproducts.  As in (\ref{dquafunc}) we view $B,\bar
C$ in the right/left crossed $kG$-module categories via the induced actions
$b\ra g=b\beta(|b|,g)$ and $g\la c=\beta^{-1}(g,|c|)c$. It is easy to see that
(\ref{bippaira}) are satisfied; the bicovariance of the pairing requiring the
symmetry of $\beta$. We then apply Theorem~B.1 and compute the structure of
$\CU$ as follows: the $\HB$ and $\bar C\lbiprod H$ algebras are the
corresponding right and left cross products as shown, while the cross relations
(\ref{biprelb}) and coproduct (\ref{bipcoprod}) become
\eqn{Gbip}{ b\Bo c\baro |c\bart|\<c\bart,b\Bt\>=\<c\baro,b\Bo\> |b\Bo|c\bart
b\bart,\quad \Delta b=b\Bo\tens |b\Bo|b\Bt,\quad \Delta c=c\baro |c\bart|\tens
c\bart,}
which computes as stated on braided-primitive elements. \endproof

The same construction works if $\beta$ is skew-symmetric as in \cite{Sch:gen}.
The only difference is that in this case we suppose that $\<\ ,\ \>$ preserves
the grading in the sense $\<c,b\>|b|=|c|\<c,b\>$ (no inversion). This case is
not so interesting, however, because on braided-primitive elements one then has
$[b,c]=0$. Both cases work more generally; the
skew case works with $kG$ replaced by any triangular Hopf algebra $H$. We just
consider $D\in{}^H\CM$ dually paired with $B\in \CM^H$ in a covariant manner
(in the sense $c\bo\<c\bt,b\>=\<c,b\bo\>b\bt$), and set $\bar C=D^{\und{\rm
cop}}\in {}^{\bar H}\CM$. However, we again have $[b,c]=0$ on braided-primitive
elements.

It seems likely that there is variant of Theorem~B.1 which works for general
dual quasitriangular Hopf algebras too. The present version is enough to
include Example~B.2, which generalises Lusztig's construction in Section~4 to
the case of a symmetric bicharacter on any (not-necessarily free) Abelian group
$G$ and a pair of suitably dual $G$-graded braided groups. Note however, that
there is in general no surjection to $\CU$ from the quantum double as in the
preceding Appendix~A, and hence no corresponding quasitriangular structure in
the finite-dimensional non-degenerately paired case. As with single
bosonisations and their duals, the key properties of $U$ in Section~3 do not
come from this biproduct point of view, though it may be a useful as a
complement.

\section{Appendix: Diagrammatic construction of $U$}

In this Appendix we provide a more category-theoretic point of view on the
$U(B)$ construction in Section~3: we construct a braided category and recover
$U$ by Tannaka-Krein reconstruction as its generating Hopf algebra. The
construction is more general than the algebraic one in Section~3, though this
remains the main example.

We work in a general braided monoidal or quasitensor category $\CC$. Let $B$ be
a braided group (Hopf algebra) in $\CC$ and let $C$ be another braided group in
$\CC$ which is dually paired with $B$ from the left, i.e. there is a morphism
$\ev:C\tens B\to \und 1$ obeying the categorical duality axioms as explained in
the Preliminaries. Here $\und 1$ denotes the identity object in $\CC$. We use a
diagrammatic notation in which braidings are denoted
$\Psi=\epsfbox{braid.eps}$, inverse braidings by the reversed braid crossing,
$\ev=\cup$ and $\und 1$ is omitted. Other morphisms are represented by nodes
with the appropriate valency, pointing generally downwards. Unless otherwise
labelled, $\cdot=\epsfbox{prodfrag.eps}$ a product and
$\und\Delta=\epsfbox{deltafrag.eps}$ a coproduct. We suppress the associativity
morphisms in the category. Functoriality of the braiding says that morphisms on
nodes can be pulled through braid crossings on either side. Algebra and Hopf
algebra in braided categories using this notation was introduced by the author,
in \cite{Ma:bra}\cite{Ma:bos}\cite{Ma:tra}\cite{Ma:introp}\cite{Ma:introm}.

\begin{defin} Let $C,B$ be dually paired braided groups in $\CC$. We define the
category ${}_C\CC_B$ of {\em braided crossed $C-B$-bimodules} to consist of
objects $(V,\la,\ra)$ where $V$ is an object of $\CC$, $\la$ is a left action
in the category of $C$ on $V$ and $\ra$ is a right action of $B$ on $V$ such
that
\ceqn{crosscond}{
(\la\tens\ev)\circ(\id\tens\Psi_{C,V}\tens\id)\circ
(\id\tens\id\tens\ra\tens\id)\circ(\Delta_C\tens\id\tens\Delta_B)\qquad\\
\qquad=(\ev\tens\ra)\circ(\id\tens\Psi_{V,B}\tens\id)\circ
(\id\tens\la\tens\id\tens\id)\circ(\Delta_C\tens\id\tens\Delta_B)}
as morphisms $C\tens V\tens B\to V$. The condition is shown in Figure~1(a).
Morphisms in ${}_C\CC_B$ are morphisms in $\CC$ which intertwine the actions of
$C,B$.
\end{defin}

\begin{propos} The category ${}_C\CC_B$ is monoidal, where the tensor product
is the usual tensor product of $C$ and $B$ modules in $\CC$. The forgetful
functor ${}_C\CC_B\to \CC$ is monoidal.
\end{propos}
\proof The proof is shown in Figure~1(c). The left hand diagram is the
condition (\ref{crosscond}) from Figure~1(a) applied to the module $V\tens W$.
The latter is a right $B$ module via the coproduct of $B$ and the braiding
$\Psi$ of $\CC$, and a left $C$-module via the coproduct of $C$ in the usual
way\cite{Ma:introp}. Iterated coproducts
$(\id\tens\und\Delta)\circ\und\Delta=(\und\Delta\tens\id)\circ\und\Delta$  are
depicted by nodes with 1 input line and 3 output lines. The first identity
applies the assumed condition (\ref{crosscond}) for $W$ (the upper layer in the
figure). The second identity then applies the assumed condition
(\ref{crosscond}) for $V$ (the lower layer). The result is the right hand side
of Figure~1(a) in for the tensor product module $V\tens W$. The identity object
$\und 1$ equipped with the trivial module structures via the counits of $C,B$
provides the identity object of ${}_C\CC_B$. Since we build the tensor product
on $V\tens W$, the forgetful functor of $\CC$ is monoidal. \endproof

\begin{figure}

\[ \epsfbox{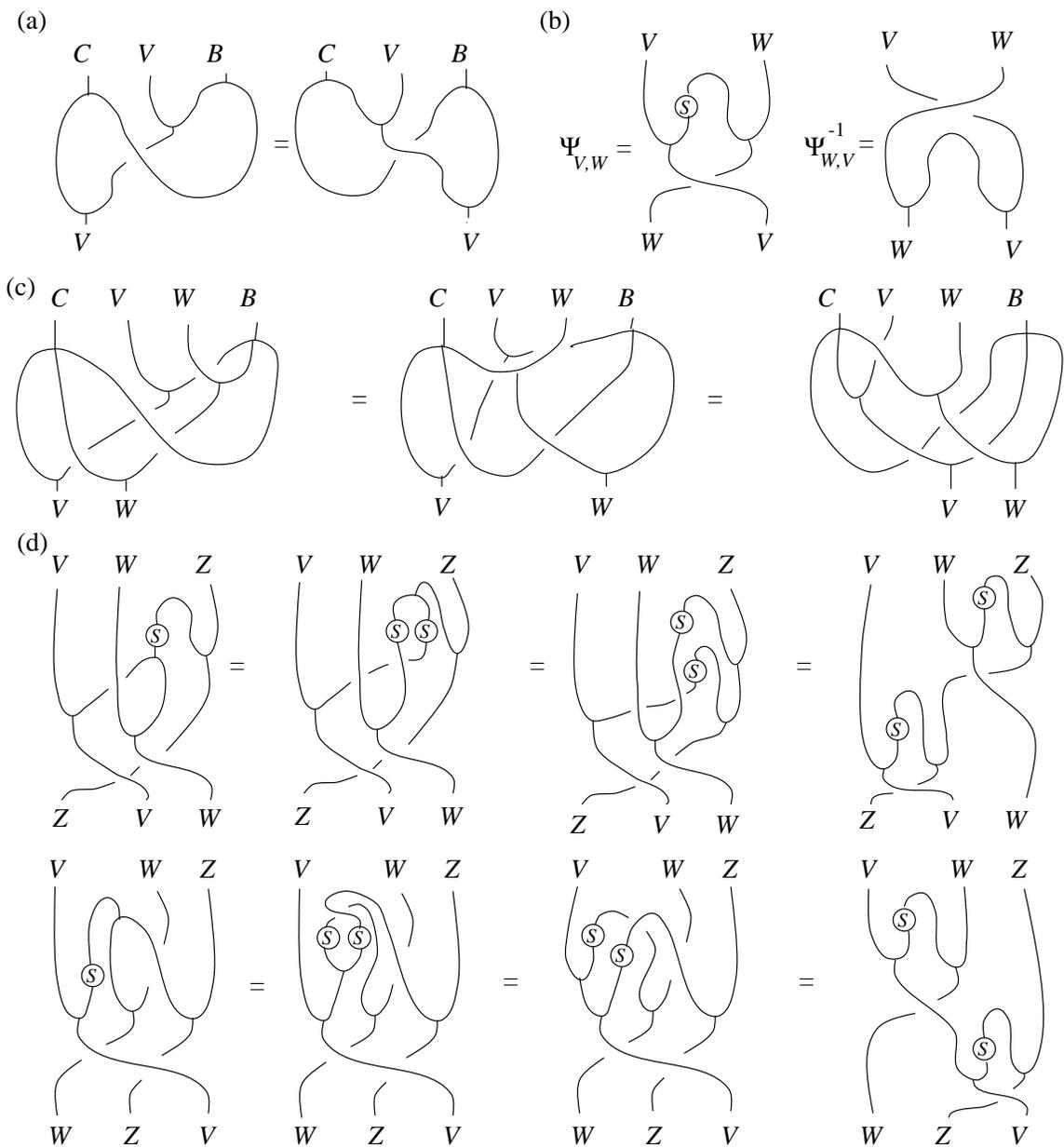}\]
 \caption{Compatibility condition (a) defining the category of $C-B$-crossed
bimodules. Its braiding (b) in the rigid case. Proof (c) that the category is
monoidal and proof (d) that the braiding obeys hexagon identities}
\end{figure}

If there is also a coevaluation $\coev:\und 1\to B\tens C$ making $B$ a rigid
object in the category, we write  $C=B^*$ and $\coev=\cap$. The
`bend-straightening axioms' pertain, see \cite{Ma:introm}. In this case we can
turn right $B^*$-modules into right $B$-comodules in $\CC$. Then
${}_{B^*}\CC_B$ is equivalent to the category of crossed $B$-modules
${}^B\CC_B$ as in \cite{Bes:cro}\cite{Dra:bos}, where the latter category is
shown to be braided. In our formulation the same observation is:

\begin{propos} cf.\cite{Bes:cro}\cite{Dra:bos} ${}_{B^*}\CC_B$ is braided by
\[ \Psi_{(V,\ra,\la),(W,\ra,\la)}=\Psi_{V,W}\circ(\ra\tens\la)\circ (\id\tens
S\tens\id\tens\id)\circ(\id\tens\coev\tens\id)\]
as shown in Figure~1(b).
\end{propos}
\proof This is shown in Figure~1(d). The upper left hand side shows
$\Psi_{(V\tens W,\ra,\la),(Z,\ra,\la)}$ where we use the tensor product action
on $V\tens W$ from Proposition~C.2. The first equality is the
braided-antimultiplicativity of $S$ proven in \cite{Ma:tra}. The second
identity dualises the coproduct of $B$ as a product in $B^*$ and then writes
the action by this product as two applications of $\la$. The third identity
uses functoriality to write the result as
$\Psi_{(V,\ra,\la),(Z,\ra,\la)}\circ\Psi_{(W,\ra,\la),(Z,\ra,\la)}$.
This verifies one of the so-called hexagon identities for the braiding in
${}_{B^*}\CC_B$. The second line in
Figure~1(d) is the proof of the other hexagon, and is similar. The first
identity dualises the coproduct in $B^*$ as a product in $B$ and applies
braided-antimultiplicativity of the braided antipode. The second identity uses
that $\ra$ is an action. Functoriality then provides the required right hand
side. The inverse braiding shown in Figure~1(b) is clearly the inverse morphism
after we use that $\ra,\la$ are actions to write the composition of
$\Psi,\Psi^{-1}$ as the action of products in $B^*,B$. The product in $B^*$ can
then be written as the coproduct in $B$, providing an antipode loop, which we
cancel (by the axiom of a braided antipode). \endproof

These steps are similar to the study of the representations of the usual
quantum double in \cite{Ma:dou}, except that now all modules are objects in a
background braided category. Just as the quantum double has a canonical
`Schr\"odinger representation' by the coregular and adjoint actions, we show
now that the same holds in our braided setting.

\begin{propos} Let $C,B$ be dually paired braided groups in $\CC$. Then $V=B$
is an algebra in ${}_C\CC_B$ where $B$ acts on itself by the right adjoint
action\cite{Ma:lie} (upper box in Figure~2) and $C$ acts on $B$ by the left
coregular representation\cite{Ma:introm}\cite{Ma:qsta} (lower box in Figure~2).
\end{propos}
\begin{figure}
 \[ \epsfbox{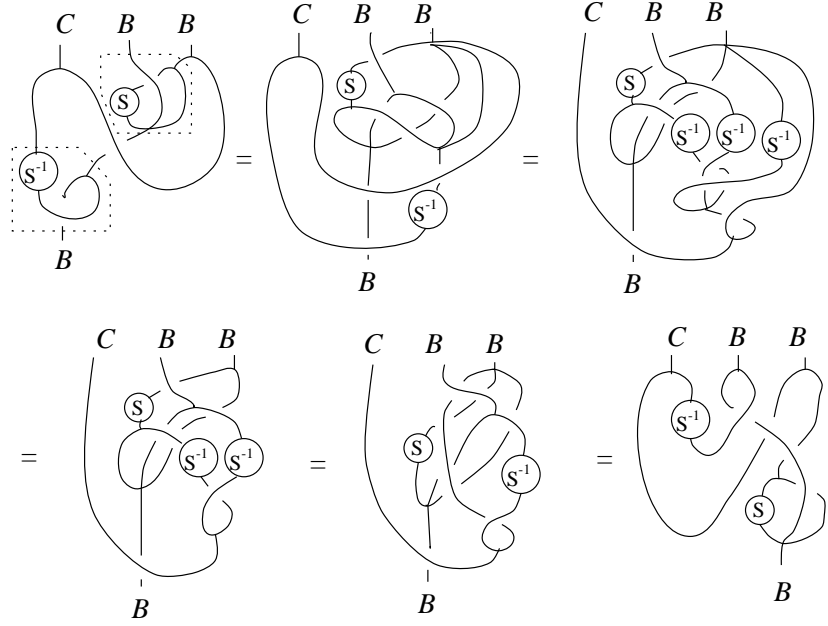}\]
 \caption{Proof that $B$ itself is an object in the category of $C-B$-crossed
bimodules}
\end{figure}
\proof This is shown in Figure~2. The right adjoint action shown is the mirror
image of the left adjoint action in \cite{Ma:lie}. More precisely, consider the
proofs for the left adjoint action of a braided group on itself, reflect in a
mirror about a vertical axis and reverse all braid crossings. This gives the
proof that the $B$ is a right $B$-module algebra by the right-handed $\Ad$ as
shown. The left coregular representation shown in the lower box is the right
coaction provided by the  coproduct, converted to a left action of $C$ via the
evaluation pairing (a coevaluation is not required)\cite{Ma:introp}.
Equivalently, it is the left braided differentiation in \cite{Ma:fre} provided
by evaluation against the left output of the coproduct, but applied to
$B^{\und{\rm cop}}$ with coproduct $\Psi^{-1}\circ\und\Delta$. The braided
antipode converts   $B$ to a braided left $C$-module algebra according to
\cite{Ma:introp}. Hence $V=B$ is a right $B$-module algebra in $\CC$ and a left
$C$-module algebra in $\CC$.

It remains to show that the left and right actions obey (\ref{crosscond}) in
Figure~1(a). The first expression in Figure~2 is the left hand side of this
condition. The first equality is uses the coproduct homomorphism property to
compute the coproduct of a triple-product. We combine iterated coproducts into
a multiple node. The second identity is the braided-antimultiplicativity of the
inverse braided antipode\cite{Ma:introm}, iterated. We also write the coproduct
of $C$ as a product in $B$, in view of their duality pairing. The third
equality cancels a loop involving the inverse braided antipode. The latter is
the antipode for $B^{\und{\rm op}}$\cite{Ma:introm}.
The fourth equality uses the braided-antimultiplicativity of $\und S$ to
simplify further. The last equality writes a product and inverse braided
antipode of $B$ back in terms of the coproduct and inverse braided-antipode of
$C$. We obtain the right hand side of the condition (\ref{crosscond}).
\endproof

The corresponding result for braided crossed modules is also new: the braided
group $B$ provides a canonical algebra in ${}^B\CC_B$ by the right adjoint
action and the coproduct viewed as a left $B$-comodule in the category.

These are general categorical constructions. We now apply them in the case
$\CC=\CM_H$ where $H$ is a quasitriangular Hopf algebra with quasitriangular
structure $\CR$. Similarly in ${}^A\CM$ for a weakly quasitriangular Hopf
algebra dual pair $(H,A,\CR)$. We consider the {\em completely forgetful
functor} which is the forgetful functor $\CM_H\to\Vec$ composed with the
forgetful functor ${}_C\CC_B\to \CM_H$. It is manifestly monoidal and hence, by
general Tannaka-Krein arguments we deduce the existence of a Hopf algebra
through the (say, right) modules of which it factors. In the present setting we
have an equivalence of categories:

\begin{propos} Let $C,B$ be dually paired braided groups in $\CM_H$. Then the
category ${}_C\CC_B$ of crossed $C-B$-modules is monoidally equivalent (in a
way compatible with the forgetful functors to $\Vec$) to the category of right
$U$-modules. If $B$ is rigid then the equivalence becomes one of braided
categories.
\end{propos}
\proof The calculation has been done in the proof of Lemma~3.11, where we used
the concrete version of the diagrammatic condition (\ref{crosscond}), obtained
from the form of braiding (\ref{lrbraid}) in the category $\CM_H$. We have seen
that the left action of $C$ can be viewed as a right action of $\bar C$, in
which case the condition becomes the relation (\ref{Urela}) in $U$. The $\CH$
relations and $\HB$ relations are just the condition that $\la,\ra$ are
morphisms in $\CM_H$, as explained in the course of Lemma~3.11. If we know that
these relations fully characterise $U$, as we know from Section~3 (see
Remark~3.7), we conclude the result. \endproof

The fundamental representation in Proposition~C.5 provides in the case of
$\CM_H$, the concrete representation of $U(\bar C,H,B)$ on $B$ in Theorem~3.12.
We read off the formulae there from the braiding (\ref{lrbraid}) in $\CM_H$. It
is also possible to understand in categorical terms the projection $\pi$ from
the quantum double in Appendix~A. \note{Here the bosonisation $H\rbiprod B$ is
the representing object for the completely forgetful functor from the category
$\CC_B$ of $B$-modules in $\CC$\cite{Ma:bos}. Hence the modules of the quantum
double of $H\rbiprod B$ generates a braided category which is the double or
dual category of $\CC_B$. This consists of objects which are $B$-modules in
$\CC$ equipped with}
%\bibliographystyle{unsrt}
%\bibliography{biblio}

\begin{thebibliography}{10}

\bibitem{Ma:bra}
S.~Majid.
\newblock Braided groups and algebraic quantum field theories.
\newblock {\em Lett. Math. Phys.}, 22:167--176, 1991.

\bibitem{Ma:bg}
S.~Majid.
\newblock Braided groups.
\newblock {\em J. Pure and Applied Algebra}, 86:187--221, 1993.

\bibitem{Ma:introm}
S.~Majid.
\newblock Algebras and {H}opf algebras in braided categories.
\newblock volume 158 of {\em Lec. Notes in Pure and Appl. Math}, pages 55--105.
  Marcel Dekker, 1994.

\bibitem{Ma:introp}
S.~Majid.
\newblock Beyond supersymmetry and quantum symmetry (an introduction to braided
  groups and braided matrices).
\newblock In M-L. Ge and H.J. de~Vega, editors, {\em Quantum Groups, Integrable
  Statistical Models and Knot Theory}, pages 231--282. World Sci., 1993.

\bibitem{Ma:poi}
S.~Majid.
\newblock Braided momentum in the {$q$}-{P}oincar{\'e} group.
\newblock {\em J. Math. Phys.}, 34:2045--2058, 1993.

\bibitem{Ma:tra}
S.~Majid.
\newblock Transmutation theory and rank for quantum braided groups.
\newblock {\em Math. Proc. Camb. Phil. Soc.}, 113:45--70, 1993.

\bibitem{Ma:bos}
S.~Majid.
\newblock Cross products by braided groups and bosonization.
\newblock {\em J. Algebra}, 163:165--190, 1994.

\bibitem{Ma:qsta}
S.~Majid.
\newblock Quasi-{$*$} structure on {$q$}-{P}oincar{\'e} algebras.
\newblock {\em Preprint}, Damtp/95-11, 1995.

\bibitem{Ma:skl}
S.~Majid.
\newblock Braided matrix structure of the {S}klyanin algebra and of the quantum
  {L}orentz group.
\newblock {\em Commun. Math. Phys.}, 156:607--638, 1993.

\bibitem{Lus}
G.~Lusztig.
\newblock {\em Introduction to Quantum groups}.
\newblock Birkhauser, 1993.

\bibitem{Dri}
V.G. Drinfeld.
\newblock Quantum groups.
\newblock In A.~Gleason, editor, {\em Proceedings of the {ICM}}, pages
  798--820, Rhode Island, 1987. AMS.

\bibitem{Jim:dif}
M.~Jimbo.
\newblock A {$q$}-difference analog of {$U(g)$} and the {Y}ang-{B}axter
  equation.
\newblock {\em Lett. Math. Phys.}, 10:63--69, 1985.

\bibitem{Ma:csta}
S.~Majid.
\newblock {$\C$}-statistical quantum groups and {W}eyl algebras.
\newblock {\em J. Math. Phys.}, 33:3431--3444, 1992.

\bibitem{Ma:fre}
S.~Majid.
\newblock Free braided differential calculus, braided binomial theorem and the
  braided exponential map.
\newblock {\em J. Math. Phys.}, 34:4843--4856, 1993.

\bibitem{Ma:lie}
S.~Majid.
\newblock Quantum and braided {L}ie algebras.
\newblock {\em J. Geom. Phys.}, 13:307--356, 1994.

\bibitem{Rad:str}
D.~Radford.
\newblock The structure of {H}opf algebras with a projection.
\newblock {\em J. Algebra}, 92:322--347, 1985.

\bibitem{Bes:cro}
Yu.~N. Bespalov.
\newblock Crossed modules and quantum groups in braided categories i,ii.
\newblock {\em Preprints}, 1994.

\bibitem{Dra:bos}
B.~Drabant.
\newblock Braided bosonization and inhomogeneous quantum groups.
\newblock {\em Preprint}, 1994.

\bibitem{Swe:hop}
M.E. Sweedler.
\newblock {\em {H}opf Algebras}.
\newblock Benjamin, 1969.

\bibitem{Ma:pro}
S.~Majid.
\newblock Quantum groups and quantum probability.
\newblock In {\em Quantum Probability and Related Topics VI (Proc. Trento,
  1989)}, pages 333--358. World Sci.

\bibitem{Ma:tan}
S.~Majid.
\newblock Tannaka-{K}rein theorem for quasi{H}opf algebras and other results.
\newblock {\em Contemp. Maths}, 134:219--232, 1992.

\bibitem{Ma:mor}
S.~Majid.
\newblock More examples of bicrossproduct and double cross product {H}opf
  algebras.
\newblock {\em Isr. J. Math}, 72:133--148, 1990.

\bibitem{Ma:qua}
S.~Majid.
\newblock Quasitriangular {H}opf algebras and {Y}ang-{B}axter equations.
\newblock {\em Int. J. Modern Physics A}, 5(1):1--91, 1990.

\bibitem{Ma:dou}
S.~Majid.
\newblock Doubles of quasitriangular {H}opf algebras.
\newblock {\em Commun. Algebra}, 19(11):3061--3073, 1991.

\bibitem{Ma:mec}
S.~Majid.
\newblock The quantum double as quantum mechanics.
\newblock {\em J. Geom. Phys.}, 13:169--202, 1994.

\bibitem{KemMa:alg}
A.~Kempf and S.~Majid.
\newblock Algebraic $q$-integration and {F}ourier theory on quantum and braided
  spaces.
\newblock {\em J. Math. Phys.}, 35:6802--6837, 1994.

\bibitem{LyuMa:bra}
V.V. Lyubashenko and S.~Majid.
\newblock Braided groups and quantum {F}ourier transform.
\newblock {\em J. Algebra}, 166:506--528, 1994.

\bibitem{Ma:any}
S.~Majid.
\newblock Anyonic quantum groups.
\newblock In Z.~Oziewicz et~al, editor, {\em Spinors, Twistors, Clifford
  Algebras and Quantum Deformations (Proc. of 2nd Max Born Symposium, Wroclaw,
  Poland, 1992)}, pages 327--336. Kluwer.

\bibitem{KirRes:rep}
A.N. Kirillov and N.Yu. Reshetikhin.
\newblock Representations of the algebra {$U_q(sl(2))$}, {$q$}-orthogonal
  polynomials and invariants of links.
\newblock In V.G. Kac, editor, {\em Infinite-Dimensional {L}ie Algebras and
  Groups}, pages 285--339. World Sci., 1989.

\bibitem{Lus:roo}
G.~Lusztig.
\newblock Quantum groups at roots of {$1$}.
\newblock {\em Geometricae Dedicata}, 35:89--114, 1990.

\bibitem{FRT:lie}
L.D. Faddeev, N.Yu. Reshetikhin, and L.A. Takhtajan.
\newblock Quantization of {L}ie groups and {L}ie algebras.
\newblock {\em Leningrad Math. J.}, 1:193--225, 1990.

\bibitem{Mey:new}
U.~Meyer.
\newblock {$q$}-{L}orentz group and braided coaddition on $q$-{M}inkowski
  space.
\newblock {\em Commun. Math. Phys.}, 168:249--264, 1995.

\bibitem{MaMar:glu}
S.~Majid and M.~Markl.
\newblock Glueing operation for {$R$}-matrices, quantum groups and link
  invariants of {H}ecke type.
\newblock {\em Math. Proc. Camb. Phil. Soc.}, 118, 1995.

\bibitem{Yet:rep}
D.N. Yetter.
\newblock Quantum groups and representations of monoidal categories.
\newblock {\em Math. Proc. Camb. Phil. Soc.}, 108:261--290, 1990.

\bibitem{Coh:hop}
M.~Cohen.
\newblock {H}opf algebras acting on semiprime algebras.
\newblock {\em Contemp. Math.}, 43:49--61, 1985.

\bibitem{Sch:gen}
M.~Scheunert.
\newblock Generalized {L}ie algebras.
\newblock {\em J. Math. Phys.}, 20:712--720, 1979.

\end{thebibliography}

\end{document}